\documentclass[letterpaper,12pt]{article}
\usepackage{mathtools} 
\usepackage{latexsym,amssymb,amsmath,amscd,amsfonts}
\usepackage{graphicx} 
\usepackage{url}
\usepackage{color}
\usepackage{hyperref}
\hypersetup{colorlinks=true}

\newcommand{\units}{\sim}
\newcommand{\mystar}{{\boldsymbol{\star}}}

\newcommand{\bdot}{\boldsymbol{\cdot}}

\newcommand{\real}{\mathbb{R}}

\newcommand{\lbilin}{\left<\left<}
\newcommand{\rbilin}{\right>\right>}

\newcommand\mystack[2]{\genfrac{}{}{0pt}{}{#1}{#2}}

\newcommand{\half}{{\frac{1}{2}}}
\newcommand{\thalf}{{\frac{3}{2}}}

\newcommand{\dx}{\triangle x}
\newcommand{\dy}{\triangle y}
\newcommand{\dz}{\triangle z}
\newcommand{\dt}{\triangle t}

\newcommand{\grad}{{\vec{\nabla}}}
\newcommand{\divg}{{\vec{\nabla} \bdot }}
\newcommand{\curl}{{\vec{\nabla} \times }}
\newcommand{\lap}{{\Delta}}
\newcommand{\lapv}{{\mathbf \Delta}}

\newcommand{\Nodes}{\mathcal N}
\newcommand{\Edges}{\mathcal E}
\newcommand{\Faces}{\mathcal F}
\newcommand{\Cells}{\mathcal C}
\newcommand{\StarNodes}{{{\mathcal N}^\mystar}}
\newcommand{\StarEdges}{{{\mathcal E}^\mystar}}
\newcommand{\StarFaces}{{{\mathcal F}^\mystar}}
\newcommand{\StarCells}{{{\mathcal C}^\mystar}}

\newcommand{\SpaceN}{{S_\Nodes}}
\newcommand{\SpaceE}{{V_\Edges}}
\newcommand{\SpaceF}{{V_\Faces}}
\newcommand{\SpaceC}{{S_\Cells}}
\newcommand{\StarSpaceN}{{S_{\Nodes^\mystar}}}
\newcommand{\StarSpaceE}{{V_{\Edges^\mystar}}}
\newcommand{\StarSpaceF}{{V_{\Faces^\mystar}}}
\newcommand{\StarSpaceC}{{S_{\Cells^\mystar}}}

\newcommand{\GRAD}{{\mathcal G}}
\newcommand{\CURL}{{\mathcal R}}
\newcommand{\DIVG}{{\mathcal D}}

\newcommand{\GRADstar}{{{\mathcal G}^\mystar}}
\newcommand{\CURLstar}{{{\mathcal R}^\mystar }}
\newcommand{\DIVGstar}{{{\mathcal D}^\mystar }}

\newcommand{\Dx}{{\Delta x}}
\newcommand{\Dy}{{\Delta y}}
\newcommand{\Dz}{{\Delta z}}

  \newcommand{\norm}[1]{\left|\left|#1\right|\right|}
\usepackage[mathscr]{euscript}
\numberwithin{figure}{section}
\numberwithin{table}{section}

\begin{document}
\title{Explicit Time Mimetic Discretiztions Of Wave Equations \thanks{ }}
\author{
    Stanly L. Steinberg \\
    Department of Mathematics and Statistics \\
    University of New Mexico, 
    Albuquerque NM 87131-1141 USA
    }

\maketitle
\newpage \clearpage
\tableofcontents
\listoftables
\listoffigures
\newpage
\begin{abstract}

Updated versions of this paper will appear at arXiv.
This version has be revised thru Section 4.

This paper is part of a program to combine a staggered time and staggered
spatial discretization of continuum wave equations so that important properties
of the continuum that are proved using vector calculus can be proven in an
analogous way for the discretized system. The spatial discretizations are
second order accurate and mimetic.  The time discretizations second are
order accurate explicit leapfrog schemes.  The discretizations have a
conserved quantity that guarantees stability for a reasonable constraint on
the time step.  The conserved quantities are closely related to the energy
for the continuum equations.  The well known Yee grid discretization of
Maxwell's equations \cite{Yee1966} is the same as discretization described
here and is an early example of using a staggered space and time grids.

Motivation of the discussion begins by studying the discretization of the
harmonic oscillator and then use this to introduce the modification of the
usual discrete energy that is conserved.  Next the discretization and creation
of a conserved quantity for possibly infinite liner systems of constant
coefficient wave equations is used to show how to create conserved quantities
for partial differential equations.

As an introduction to discretizing higher dimensional wave equations the one
dimensional wave equation with variable coefficients is used to motivate the
rest of this paper. The discussion begins with the constant coefficient wave
equation and is then extended to variable coefficients. The conserved quantity
contains a term proportional to the square of the time step.  Requiring that
the conserved quantity is positive is the well known Courant-Friedrichs-Lewy
(CFL) condition for stability. For constant material properties and some
variable material properties the simulation codes are second order accurate
but for other materials they are only close to second order accurate.
The cause and importance of this is under study.

To study partial differential equations in higher dimensions the notion of
exact sequences and diagram chasing from differential geometry are used to
produce second order variable coefficient differential operators that are self
adjoint and either positive or negative and thus can be used to define many of
the known wave equations. These operators are then used to define second
order wave equations and their equivalent first order systems.

Next mimetic spatial finite difference approximations and leap-frog time
discretization are use to discretize the first order systems of wave equations
in one, two and three dimensional spaces. The discretization also provide
solutions of the second order wave equation that are equivalent to the first
order system. The discretizations are at least second order accurate and have
a conserved quantity that guarantees stability for a modest restriction on the
time step. The conserved quantities imply the conservation of energy.

MatLab simulation codes are provided for all of equations studied.

\end{abstract}

\noindent \textit{Key Words: mimetic discretization, leapfrog,
	energy conservation, stability}

\newpage \clearpage
\setcounter{equation}{0}
\section{Introduction}

The main goal is to show that second-order accurate mimetic finite difference
spatial discretizations can be combined with an {\em explicit} second order
accurate leapfrog finite difference time discretizations to discretize wave
equations to produce stable second-order accurate simulations and then to
create a general method for extending mimetic discretizations to model waves in
inhomogeneous and anisotropic materials. Most importantly, discrete
conserved quantities that have simple relationships to the energy are
introduced.  These discrete conserved quantities are positive quadratic
forms for modest restrictions on the time step and consequently imply the
discretizations are stable.  It is important that the material properties do not
depend on the solution of the differential equation so that the equation are
linear and additionally that the material properties do not depend on time.
The techniques depend on writing second order wave equations as systems of
two first order equations.  The spatial and temporal units of variables play
and important role in the discussion. Simulation codes are provided for all
the examples discussed.  It is an open question if these ideas can be extended
to higher order discretizations or to nonlinear differential equations.
Simulation programs are provided for the differential equations.  These codes
confirm that the that the conserved quantities are constant to within a small
multiple of machine epsilon {\tt eps}. For constant and some variable
material properties the codes confirm the second order accuracy. For other
variable material properties the codes are only close to second order accurate.
The cause of this is currently under investigation.

In Section 2 {\em The Harmonic Oscillator and Conserved Quantities} the
standard second-order harmonic oscillator ordinary differential equation is
written as a first order system which is discretized using the standard leap
frog or staggered in time method. For this discretization a conserved
quantity is created that is closely related to the oscillator energy.
This quantity depends on the time step and is positive for small time steps
thus guaranteeing stability of the discretization. The constraint on the time
step for reasonable accuracy is much smaller that the constraint for stability.

In Section 3 {\em Systems of Ordinary Differential Equations} the results for
the harmonic oscillator are extended to possibly infinite systems of first
order linear and constant coefficient ordinary differential equations that are
wave equations. For the discrete system conserved quantities that are
generalization of those for the harmonic oscillator are derived.  Simple
generalizations of these conserved quantities will provide conserved
quantities for variable material properties wave equations.

In Section 4 {\em Discretizing the One Dimensional Wave Equation} the
discussion of spatial discretizations of partial differential equations begins
by discretizing wave equations with variable material properties in one
spatial dimension.  First the second order one dimensional wave equation is
derived from the continuum conservation of energy to make clear the correct
form of a wave equation with variable material properties. The second order
equation is then written as a first order system and then this provides
another second order wave equation. The critical point is that weighted inner
products can be introduced so that the system has the same form as the
equations in Section 3. This immediately gives a conserved quantity that
implies the conservation of energy.

First the constant coefficient system is discretized followed by the variable
coefficient case.  These discretizations use the same grids that are staggered
in space and time.  In space, as in standard in mimetic discretizations, two
grids are used where a point in one of the grids is at the midpoint of the
cells in the other grid.  The time discretization is the standard leapfrog
scheme.  This results in each first order equation being discretized on
separate space and time grids.  The discussion in Section 3 then easily
provides a conserved quantities that contains a term that is proportional to
minus the square of the time step times a positive quantity.

In the constant coefficient discretization the wave speed $c$ is introduced so
that it will be easy to compare the mimetic discretization to the usual
discretizations found in the literature.  Importantly, requiring the conserved
quantity to be positive produces the standard Courant-Friedrichs-Lewy (CFL)
condition for stability. For the variable coefficient case an estimate of
the maximum wave speed can be introduced that provides an analog of the
CFL condition for stability.

Two codes {\tt Wave1DCMP.m} and {\tt Wave1DVMP.m} that implement the constant
material properties {\tt CMP} and variable material properties {\tt VMP}
discretizations are described. For the constant materials case both codes
produce the same results. Both codes keep a their conserved quantity constant
to within a small multiple of machine epsilon {\tt eps}.  Simulations confirm
that the {\tt CMP} code is second order accurate. A significant puzzle is that
in some cases the {\tt VMP} code solution has oscillatory errors that start
near the boundary and propagate to the interior of the spatial domain. The
convergence rate is reduce but is still greater than 1. In a sense made clear
in this section the {\tt VMP} discrete solutions are close to second order
accurate. On the other hand when the material properties are given by
$1+\left(2 \, x \, (1-x)\right)^2$ the VMP code is second order accurate and
there are no small oscillatory errors.

The first 4 sections have been revised. Still working on the following sections.
The plan is to first create the three dimensional theory and codes and then
use this to implement the two dimensional codes because mimetic differential
operators are easier to understand in three dimensions than in two dimensions.
In particular the mimetic three dimensional method naturally produces the
gradient, curl and divergence operators.

In Section 5 {\em Three Dimensional Variable Coefficient Differential Operators}
the notions of exact sequences and diagram chasing used in differential 
geometry are used to create twelve second order differential operators with
variable coefficients.  These operators are built using the gradient, curl and
divergence and have coefficients that are positive functions or real symmetric
3 by 3 positive matrices that are smooth functions of the spatial variables.
These functions correspond to the star operators in differential geometry.
An additional 4 operators that are a linear combination of two of the above
operators are also be defined.

The second order operators are next written as a product of two first order
operators with variable coefficients.  Weighted inner products are introduced
and used to compute the adjoint operators of the first order operators and the
show that the second order operators are self adjoint and either positive or
negative.  Thus the negative operators and minus the positive operators can
be used to define second order wave equations.  It is not clear if all of the
16 general equations are useful but they do provide flexibility in dealing
with spatially dependent material properties.  When the scalar coefficients
are constant and the matrix coefficient are a constant multiple of the
identity matrix the second order operators are reduced to 4 unique operators.

Revisions done to here.

In Section 6 {\em 3D Wave Equations With Variable Material Properties} second
order variable coefficient wave equations are defined. Diagram chasing is used
to write the second order equations as first order systems. These systems have
the same form as those in Section \ref{ODEs} which then easily provides
conserved quantities for all of the equations. This section ends by showing
the relationship of the diagram chasing wave equations to the standard wave
equations in the literature.

In Section 7 {\em Mimetic 3D Discretizations} three dimensional staggered
grids are introduced and used to define mimetic discretizations of the
gradient, curl and divergence \cite{RobidouxSteinberg2011}.  This type of
discretization is motivated by the use of exact sequences from differential
geometry, except here there are two grids so the exact sequence becomes a
double exact sequence.  This is combined with a leapfrog discretization.
The spatial discretization uses two grids where the corners of
cells in one grid are the same points as the centers of the cells in the other
grid and the centers of the edges on one grid are the cell face centers in the
other grid.  An early example of this type of spatial and time discretization
was given by Yee \cite{Yee1966} for Maxwell's equations.  Mimetic spatial
discretizations have been used extensively to create simulation programs for
problems in continuum mechanics, see \cite{LipnikovMS14} and the volume
\cite{Koren2014} in which this work appeared. They have also been used
to model inhomogeneous and anisotropic materials in two dimensions
\cite{HymanMorel02,HymanSS97}. For a comparison of mimetic finite difference,
finite volume and finite element discretizations see \cite{BochevHyman06}.

The mimetic spatial discretization studied here are extensions of those
described in \cite{RobidouxSteinberg2011} where it is shown that mimetic
discretizations of the gradient, curl and divergence satisfy important
properties of the continuum operators. For example the discrete divergence
of the discrete curl operator is identically zero and the adjoint of the
discrete gradient is the minus the discrete divergence.  To extend this work
to anisotropic materials it is critical that in three dimensions the
anisotropic properties of many important materials are describe by a
$3 \times 3$ symmetric positive definite matrix \cite{Nye04}, see Chapter
4, Sections 1 and 4, for the permittivity.  Additionally,
\cite{LipnikovMS14,BrezziLS05,HymanMSS02,HymanShashkov01} discuss the
incorporation of material properties into mimetic discretizations using
symmetric positive definite matrices.

In Section 8 the results in Section 7 are used to discretize the scalar wave
equation and Maxwell's equations for general material proprieties. Three
case will be considered: trivial material properties where the scalar material
property is 1 and the matrix material property is the identity matrix;
constant material proprieties where the scalar material properties is a positive
constant and the matrix material properties is a multiple of the identity
matrix and finally general material proprieties where the scalar material
property is a positive smooth function and the matrix material proprieties are
given by positive definite matrix and both material properties are smooth
functions of the spatial variable.

The code {\tt Wave3DTMP.m} shows that for the scalar wave equation with
trivial material proprieties the solution is forth order accurate and
the conserved quantity is conserved quantity is constant to a small 
multiple of machine epsilon {\tt eps}. Note that the one spatial dimension
scalar wave equation could also produce forth order accurate solutions
in some simple cases.

Too be done are the scalar wave equation with general material properties
and Maxwell's equations with general material proprieties.

It would be interesting to use ghost points to implement general Robin
boundary conditions as was done in \cite{CocoRusso13}.

\newpage \clearpage
\setcounter{equation}{0}

\section{The Harmonic Oscillator \label{Harmonic Oscillator}}

The goal is to use the time discretization of the harmonic oscillator to
motivate time discretizations of special systems of ordinary differential
equation that are wave equations, and then this is used to discretize wave
equations in 1, 2 and 3 dimensional spaces.  The discretizations are second
order accurate, explicit and conserve a quantity that is a second order
approximation of a constant multiple of the energy. Requiring this quantity
to be positive gives the standard constraint on the time step for stability
and is far less restrictive that a reasonable accuracy constraint.

First the explicit second order discretization of the second order continuum
oscillator equation is described and the novel conserved discrete quantity is
introduced. This quantity contains a term that is {\em minus} the square of
the time step times a positive quantity. The space and time units of the
conserved quantity is that of energy.  The simulation program
{\tt Oscillator2ndOrder.m} confirms that the discretization is second
order accurate and that the conserved quantity is constant within a small
multiple of machine epsilon {\tt eps}.

Next the second order equation is written as a first order system and then
this system is discretized using staggered grids which is typically called a
leapfrog discretization.  This discretization gives solutions that are
equivalent to solutions of the discretized second order equation.  Next a
conserved quantity for the discrete system is described.  As before the
simulation code {\tt OscillatorSystem.m} confirms that the method is second
order accurate and that conserved quantity gives a second order accurate
approximation of a constant multiple continuum energy.  This system was used
to make the phase plane plot \ref{OscillatorPhasePlane} for the harmonic
oscillator that illustrates how well the new conserved quantity approximates
the continuum conserved quantity.

The constraint on the time step to keep the quantity positive becomes the
Courant-Friedrichs-Lewy (CFL) condition for spatially dependent wave equations.
Appendix \ref{C-N} reviews the implicit Crank-Nicholson discretization of the
oscillator which conserves a natural second order accurate discretization of
the continuum energy.
The paper \cite{WanBihloNave2015} uses a natural discretization of the
classical energy to derive a discretization of the oscillator equation that
is equivalent to the Crank-Nicholson discretization. This discretization
is implicit.

\subsection{The Harmonic Oscillator and Conserved Quantities}

The linear harmonic oscillator equation is given by 
\[
u'' + \omega^2 u = 0 \,,
\]
where $u = u(t)$ is a smooth function of time $t$ and $u' = du/dt$,
$u'' = d^2u/dt^2$ and $\omega$ is a positive constant.  The total energy is a
multiple of the average of the kinetic and potential energies which is 
\[
E = \frac{(u')^2 + (\omega \, u)^2}{2} \,.
\]
This is conserved quantity because
\[
E' =
u'' \, u' + \omega^2 u \, u' = 
\left(u'' + \omega^2 u \right) \, u' = 0 \,.
\]
In general conserved quantities are labeled $C$ but if the conserved quantity
is a {\em constant multiple} of the energy it will usually be labeled with $E$.

In this paper time and spatial units are used extensively to help
construct appropriate models.
For the energy note that $u$ has units of distance, $t$ has units of time,
$u'$ has units distance over time and $\omega$ has units of reciprocal time:
\begin{align*}
 u & \sim d \,; \\
 t & \sim t \,;  \\
 u' & \sim d/t \,;  \\
 \omega & \sim 1/t \, \\
 E & \sim d^2/t^2 \,.
\end{align*}
So both the differential equation and the conserved quantity are
dimensionally consistent.

A natural way to write the oscillator equation as a first order system is
by introducing $v = v(t)$ and requiring
\begin{equation}\label{First Order System}
u' = - \omega \, v \,,\quad v' = \omega \, u \,.
\end{equation}
The minus sign can be put in either equation. The consequence is that
now $v$ must have units $d$ and thus is not a velocity. This could
be fixed by setting $v = u'$ but then this is inconsistent for what
must be done when the material properties are variable and the density
and tension are use to describe the physical properties of the material
the waves are traveling in as is done in Section \ref{1D Wave}. In
any case this simple example will be helpful for understanding more
complex problems.

For the system, set
\begin{equation}\label{Conserved Quantity Simple}
C = \frac{1}{2} \left( u^2 + v^2 \right)\,.
\end{equation}
This quantity is conserved because
\[
C' = u\,u' + v\,v' = u \, \omega \, v - v \, \omega \, u = 0 \,,
\]
Use \eqref{First Order System} to remove $v$ to get
\begin{equation}
C  = \frac{E}{\omega^2} \,.
\end{equation}
For partial differential equations the relationship between energy and
the conservation laws is similar but not so simple.
Note that the phase plane plots using $(u(t),v(t))$ are circles where the
radius $r$ is given by $r^2 = 2 \, C$.

Because the second order equation is linear with constant coefficients
the time derivative $u'$ satisfies the same equation a $u$. Similarly
the time derivatives $u'$ and $v'$ also satisfy the system thus creating
an infinity of conserved quantities.

The condition that $\omega > 0$ and not that $\omega \geq 0 $ is important
because for $\omega = 0$ the second order equation with $u(0)= 0$ and
$u'(0) = 1$ has the solution $u(t) = t$ for which the energy is unbounded.
However, for $\omega = 0$ the system only has only constant solutions
which have a bounded conserved quantity. So the second order equation and the
system are not consistent for $\omega = 0$.

\subsection{Discretizing the Second Order Oscillator Equation}

If $\Delta t > 0$ then a standard explicit discretization of the second order
oscillator equation using the discrete times
$t_n = n \, \Delta t$, $0 \leq n < \infty$ is where $n$ is and integer and
\[
\frac{u^{n+1} -2 \, u^n + u^{n-1}}{\Delta t ^2} + \omega^2 \, u^n = 0 ,
\quad n \geq 1 \,.
\]
Given the two initial conditions $u(0)$ and $u'(0)$ set
\begin{align}
u^0 & = u(0) \,  \nonumber \\
u^1 & = u(\Delta t) = u(0) + \Delta t \, u'(0)
- \frac{1}{2} \, \Delta t^2 \, \omega^2 \, u(0) \,. \nonumber
\end{align}
The discrete equation is then 
\[
u^{n+1} = (2 -(\omega \Delta t)^2) u^n - u^{n-1} , \quad n \geq 1 \,.
\]

A natural proposal for a second-order accurate discrete conserved quantity is 
\[
C^n =
(u^n)^2 + \left( \frac{u^{n+1}-u^{n-1}}{2 \, \omega \, \Delta t} \right)^2 \,.
\]
A little algebra shows that $C^n$ is not conserved. However this computation
shows that
\[
C^n = \left(1-\left(\frac{\omega \, \Delta t}{2}\right)^2 \right) (u^n)^2 +
\left( \frac{u^{n+1}-u^{n-1}}{2 \, \omega \, \Delta t} \right)^2 
\]
is conserved. Consequently for $0 < \Delta t /\omega < 2$ the discretization
is stable.  It is important that this constraint is less restrictive than
requiring an accurate solution.

The initial condition given above
is only first order accurate but a Taylor series expansion can
be used to make the order of accuracy higher:
\begin{align}
u^1 & = u(0)
+ \Delta t \,  u'(0)
+ \frac{1}{2} \, \Delta t^2 \, u''(0)
+ \frac{1}{6} \, \Delta t^3 \, u'''(0)
+ \cdots \nonumber \\
    & = u(0) 
+ \Delta t \,  u'(0)
- \frac{1}{2} \, \Delta t^2 \, \omega^2 \, u(0)
- \frac{1}{6} \, \Delta t^3 \, \omega^2 \, u'(0)
+  \cdots \label{HigherOrder} \\
\nonumber
\end{align}

The program {\tt Oscillator2ndOrder.m} confirms that the algorithm is stable
for $0 < \Delta t /\omega < 2$ and is second order accurate.  The discrete
conserved quantity is second order accurate and constant to within a small
multiple of machine epsilon {\tt eps}.  As an exact solutions is known for
the system the simulations codes uses the value $u(\Delta t)$ as the initial
condition. 

\subsection{Staggering the Time Discretization}

A time staggered grid is used to discretize the first order system
\eqref{First Order System} which is given by a primal grid
$t^n = n \, \Delta t$ and a dual grid $t^{n+1/2} = (n+1/2) \, \Delta t$,
where $\Delta t > 0$ and $0 \leq n < \infty$ is an integer.
The staggered or leapfrog discretization of the harmonic oscillator is
then given by
\begin{equation}\label{Main Difference Equations}
\frac{u^{n+1}-u^n}{\Delta t} = - \omega \, v^{n+1/2} , \quad  n \geq 0
\,,\quad
\frac{v^{n+1/2}-v^{n-1/2}}{\Delta t} = \omega \, u^{n} , \quad n \geq 1
\,.
\end{equation}
The minus sign can be put in either equation, but it is important for
consistent time units to have an $\omega$ in both equations.  As before,
the initial conditions $u(0)$ and $u'(0)$ are given and then $u^0 = u(0)$ and
\begin{align}
v^{1/2} & = v\left(\frac{\Delta t}{2}\right) \nonumber \\
& = v(0) +
\frac{\Delta t }{2} \, v'(0) +
\frac{1}{2} \, \left(\frac{\Delta t }{2} \right)^2 \, v''(0) + \cdots 
	\nonumber \\
& = v(0)+
\frac{\Delta t }{2} \, \omega \, u(0) 
- \frac{1}{2} \, \left(\frac{\Delta t }{2} \right)^2 \, \omega^2 \, v(0) 
+ \cdots \label{Higher Order IC} \\
\nonumber
\end{align}
The update algorithm starts with $u^0$ and $v^{1/2}$ and then for
$n \geq 0$
\begin{align*}
u^{n+1} = u^n - \Delta t \,  \omega \, v^{n+1/2} \,,\quad
v^{n+3/2} = v^{n+1/2} + \Delta t \, \omega \, u^{n+1} \,.
\end{align*}
Note that the second equation depends on the update in the first equation,
so the order of evaluation is critical.
This staggered grid discretization gives two standard single grid
discretization of the second order oscillator equation:
\[
\frac{u^{n+2} - 2 \, u^{n+1} +u^n }{\Delta t^2} + \omega^2 u^{n+1} = 0
\,;\quad
\frac{v^{n+3/2} - 2 \, v^{n+1/2} +v^{n-1/2} }{\Delta t^2} +
	\omega^2 v^{n+1/2} = 0 \,.
\]
So the solution $u$ of the fractional step methods is identical to the 
solution of the second order equations.

Again a simple proposed conserved quantity for
\eqref{Main Difference Equations} is
\begin{equation}\label{Simple Conserved}
C^n = \frac{1}{2}
\left( (u^n)^2 +\left( \frac{v^{n+1/2}+v^{n-1/2)}}{2}\right)^2 \right) \,.
\end{equation}
A little algebra gives
\[
C^{n+1} - C^n = \frac{\omega^2 \, \Delta t^2 }{4} \,
    \left((u^{n+1})^2 - (u^n)^2 \right) \,.
\]
So $C^n$ is {\em not} conserved. However, set
\[
\alpha = \frac{\omega \, \Delta t}{2} \,,
\]
and then the following two quantities are conserved:
\begin{equation}\label{Scaler Conserved n}
C^n = \frac{1}{2}
\left(
\left(1 - \alpha ^2 \right)
(u^n)^2 +\left( \frac{v^{n+1/2}+v^{n-1/2)}}{2}
\right)^2 \right) \,;
\end{equation}
\begin{equation}\label{Scaler Conserved n+1/2}
C^{n+1/2} = \frac{1}{2} \left(
\left(\frac{u^{n+1}+u^{n}}{2}\right)^2 + 
\left(1 - \alpha^2 \right) (v^{n+1/2})^2 \right) \,.
\end{equation}
The important properties for the staggered scheme are that it is explicit,
second order accurate and stable for $ \alpha = \omega \, \Delta t / 2 < 1$. 
By modifying the discretization, a similar result was obtained in
\cite{SchuhmannWeiland2001}, Equation 45, for the Yee time discretization
of Maxwell's equations.

The simulation program {\tt OscillatorSystems.m} confirms that the numerical 
solutions are second order accurate and that the two conserved quantities
are constant to within a small multiple of {\tt eps}.

\begin{figure}
\begin{center}
\includegraphics[width = 4.00in, trim = 50 0 90 240, clip]{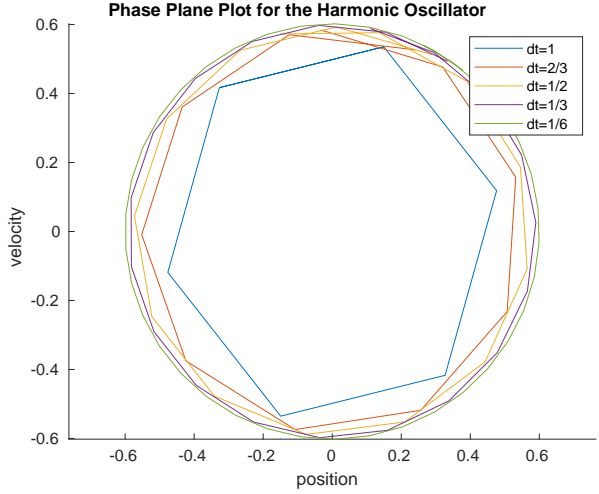}
\caption{Phase Plane Plots For The Harmonic Oscillator}
\label{OscillatorPhasePlane}
\end{center}
\end{figure}

\subsection{Summary}

For the second order harmonic oscillator and the first order harmonic
oscillator conserved quantities are introduce that are second order accurate
approximations of a multiple of the continuum energy but contain a term that
is minus the square of the time step times a positive quantity. Requiring this
conserved quantity to be positive gives the standard stability condition on
the time step for the discrete oscillator.  The restriction on $\Delta t$ to
keep the conserved quantity positive is far less stringent than the
restriction for reasonably accurate solutions.  Simulations confirm that the
discrete conserved quantities are constant up to a small multiple of machine
epsilon $\tt eps$ and are second order accurate. For the simulations code it
is a bit awkward to check but the systems code is stable for
$\Delta t/\omega = 1.99$ and unstable for $\Delta t /\omega = 2.30$ supporting
the stability condition $0 < \Delta t /\omega < 2$.  A modified version of the
systems code produced the phase plane plot in figure \ref{OscillatorPhasePlane}
using rather large time steps which confirms that the conserved quantity
converges to the continuum conserved quantity rapidly.

\newpage \clearpage
\setcounter{equation}{0}
\section{Systems of Ordinary Differential Equations \label{ODEs}}

The next task is to study a special class of systems of linear ordinary
differential equations that are wave equations. Discrete conservation laws are
easy to find by following the harmonic oscillator example. Most importantly,
the discussion here will show how to find conservation laws for discretized
spatially dependent wave equations. Consequently the notation here is set up
to conveniently apply the results to these equations.  To see how this is done
note that the standard wave equation in 3D is given by the second time
derivative of a scalar function equals the divergence of the gradient of the
function. Critical points are that the adjoint operator of the divergence is
minus the gradient and vice versa and the divergence of the gradient is a
symmetric negative operator.

\subsection{Continuous Time} \label{Continuous Time}

Let $X$ and $Y$ be linear spaces (finite or infinite dimensional).  It is
important {\em not} to assume that $X$ and $Y$ have the same dimension.
If $f1$ and $f2$ are in $X$ then their inner product is
$\left< f1 , f2 \right>_X$ and the norm of $f1$ is given by
$\norm{f1}_X^2 = \left< f1 , f1 \right>_X$ with a similar notation for $Y$.

Let $A$ be a linear operator mapping $X$ to $Y$ with adjoint $A^*$
that maps $Y$ to $X$: 
\[
X \overset{A} \rightarrow Y \,;
\quad
Y \overset{A^*} \rightarrow X \,.
\]
Consequently, if $f \in X$ and $g \in Y$ then
\[
\left< A \, f , g \right>_Y = \left< f, A^* \, g \right>_X \,.
\]

Next, if $f = f(t) \in X$ and $g = g(t) \in Y$ then the generalization of
the harmonic oscillator system that was studied in the previous section
is given by
\begin{equation} \label{Wave System}
f' = - A^* \, g \,,\quad g' = A \, f  \,.
\end{equation}
Consequently
\[
f'' = - A^* \, A f \,,\quad g'' = - A \, A^* g \,.
\]
Note that $A^* \, A$ and $A \, A^*$ are symmetric positive operators.
If they are positive definite then the solutions will be oscillatory.
When studying spatially dependent wave equations $A$ will correspond to
the gradient and $-A^*$ will correspond to the divergence so that
$- A \, A^*$ corresponds to the divergence of the gradient that is the
Laplacian in the scalar wave equation.  Also  $- A^* \, A$ corresponds
to the gradient of the divergence that occurs in the second order vector
wave equation.

The matrix form of the system is
\[
\left[ \begin{matrix}
f'  \\
g'
\end{matrix} \right]
= 
\left[ \begin{matrix}
0 & - A^* \\
A & 0 
\end{matrix} \right]
\left[ \begin{matrix}
f \\
g
\end{matrix} \right] \,.
\]
The coefficient matrix is skew adjoint, so it must have purely imaginary
spectra so the solutions of this system must be made up of waves and constant
functions.  Additionally, the adjoint of the product of two operators $A$ and
$B$ is $(A \, B)^* = B^* \, A^*$ so both $A \, A^*$ and $A^* \, A$ are self
adjoint and positive. If $A^* \, A$ is positive then the second order
equations have solution that are made up of waves.  However if $A^* \, A$ has
a zero eigenvalue then the second order equation can have solutions that grow
linearly in time.

There are three natural initial conditions: for the system specify
$f(0)$ and $g(0)$; for the second order equation in $f$ specify,
$f(0)$ and $f'(0)$; and for the second order equation in $g$ specify,
$g(0)$ and $g'(0)$. Note that if $f(0)$ and $f'(0)$ are given then
to use the first order system, Equation \eqref{Wave System} must be
solved for $g(0)$ which may not be possible. So the system and the
second order equations may not be equivalent.

If the dimensions of $X$ and $Y$ are the same so that it make sense to
assume that $A$ is inevitable then the system \eqref{Wave System} will
have properties similar to the harmonic oscillator system
\eqref{First Order System} when $\omega > 0$. The most interesting case
is when the dimensions of the spaces are different which provides insight
in to the discretization of the scalar and vector wave equation and also
the Maxwell equations.  Also the case when $A$ is self adjoint, $A^* = A$,
provides insight into the discretization of Maxwell's equations.

\subsection{Continuous Time Conserved Quantities}

An interesting point here is that there is a conserved quantity that is not
the energy but implies that an analog of the energy is conserved.  This
fundamental conserved quantity is
\[
C(t) = \frac{1}{2} \left( \norm{ f(t) }_X^2 + \norm {g(t)}_Y^2  \right) \,,
\]
which is analogous to \eqref{Conserved Quantity Simple}.  Because
\begin{align*}
C'(t)
 & = \left< f'(t), f(t) \right>_X + \left< g'(t), g(t) \right>_Y \\
 & = \left< - A^* \, g(t), f(t) \right>_Y + \left< A \, f(t), g(t) \right>_Y \\
 & = - \left< g(t), A \, f(t) \right>_X + \left< A \, f(t), g(t) \right>_Y \\
 & = 0 \,,
\\
\end{align*}
this quantity is conserved.

Additionally because the system of equations is linear and constant
coefficient, if $f$ and $g$ are solutions then so are $f'$ and $g'$.
Consequently
\begin{align}
\label{Continuum Total Energy}
E(t)
   & = \frac{1}{2} \left(\norm{f'(t)}_X^2 + \norm {g'(t)}_Y^2  \right) \,,
        \nonumber \\
   & = \frac{1}{2} \left(\norm{f'(t)}_X ^2 + \norm {A^* \, f(t) }_Y^2 \right) \,, 
        \nonumber \\
   & = \frac{1}{2} \left(\norm{A \, g(t)}_X^2 + \norm {g'(t)}_Y^2  \right) 
        \nonumber \\
\end{align}
is constant.  For the wave equations to be studied later this is essentially
the total energy that is the sum of the kinetic and potential energies.
The $C(t$) type conserved quantities will used from now on.

\subsection{Staggered Time Discretization\label{Staggered Time Discretization}}

A second order centered leapfrog discretization for the first order system is  
\begin{equation} \label{Staggered Time Equations}
\frac{f^{n+1} - f^n}{\Delta t} = -A^* \, g^{n+1/2}  \,,\quad
\frac{g^{n+1/2} - g^{n-1/2}}{\Delta t} = A \, f^{n} \,.
\end{equation}
Assuming that $f^0$ and $g^{1/2}$ are given then for $n \geq 0$ the
leapfrog time stepping scheme is 
\[
f^{n+1} = f^n - \dt \, A^* \, g^{n+1/2} \,,\quad
g^{n+3/2} = g^{n+1/2} + \dt \, A \, f^{n+1} \,.
\]
The order of evaluation is important.

The initial conditions for the discretized system require $f^0 = f(0)$
and $g^{1/2} = g(\Delta t/2)$.
If $f(0)$ and $g(0)$ are given then
\[ g^{1/2}
 \approx g(0) + \frac{\Delta t}{2} \, g'(0) + \frac{\Delta t^2}{2} g''(0)
 = g(0) + \frac{\Delta t}{2} \, A \, f(0) - \frac{\Delta t^2}{2} A \, A^* g(0) \,.
\]

The remaining cases require the solution of a system of equations, that is one
of the operators $A \, A^*$ or $A^* \, A$ needs to inverted. These operators
are symmetric and positive but under our assumptions cannot be guaranteed to
be definite.
In the other cases,
if $g(0)$ and $g'(0)$ are given then solve
\[
A^* \, g'(0) =  A^* \, A \, f(0)
\]
for $f(0)$ and if $f(0)$ and $f'(0)$ are given solve
\[
A \, f'(0) = -A \, A^* \, g(0)
\]
for $g(0)$.

Both $f$ and $g$ satisfy a second order difference equation:
\begin{align} \label{Second Order Difference}
\frac{f^{n+1} - 2 f^n + f^{n-1}}{\Delta t^2}
    & = - A^* \, A f^{n} \,;  \nonumber \\
\frac{g^{n+3/2} - 2 g^{n+1/2} + g^{n-1/2}}{\Delta t^2}
    & = - A \, A^* \, g^{n+1/2} \,. \\
\nonumber
\end{align}
Additionally a second order average is needed for computing
conserved quantities. So multiply the above equations by $\Delta t^2/4$
and then add $f^n$ to the first equation and $g^{n+1/2}$ to the 
second equation to get:
\begin{align} \label{Second Order Average}
\frac{f^{n+1} + 2 \, f^{n} + f^{n-1}}{4}
& = f^{n}  - \frac{\Delta t^2}{4} A^* \, A f^{n} \,; \nonumber \\
\frac{g^{n+3/2} + 2 \, g^{n + 1/2} + g^{n-1/2}}{4} 
& = g^{n+1/2} - \frac{\Delta t^2}{4} A \, A^* \, g^{n+1/2} \,. \\
\nonumber
\end{align}

When comparing this discretization to the simple oscillator
discretization it is important that $\omega > 0$, while here
the operators $A$ and $A^*$ may not be invertible which is typically
the case when studying spatially dependent partial differential wave equations.

\subsection{Discrete Time Conserved Quantities \label{DTCQ}}

To show that $A$ not being invertible is not serious problem a detailed
derivation of the conservation laws that are analogs of
\eqref{Scaler Conserved n} and \eqref{Scaler Conserved n+1/2} are given.
So let
\begin{align*}
C_1^{n} & = \norm{\frac{g^{n+1/2} + g^{n-1/2}}{2}}_Y^2 \,, \\
C_2^{n} & = \norm{f^{n}}_X^2 \,, \\
C_3^{n} & = \norm{ A \, f^{n} }_Y^2 \,. 
\end{align*}
As above compute the changes in $C_1$, $C_2$ and $C_3$ over a time step:
\begin{align*}
C_1^{n+1}- C_1^{n}
& = \left<
\frac{g^{n+3/2} + 2 \, g^{n+1/2} + g^{n-1/2}}{4} \,,\,
   g^{n+3/2} -g^{n-1/2} \right>_Y ; \\
& = \left<
g^{n+1/2}  - \frac{\Delta t^2}{4} A \, A^* \, g^{n+1/2} \,,\,
   g^{n+3/2} -g^{n-1/2} \right>_Y ;  \\
& = \left< g^{n+1/2} \,, g^{n+3/2} -g^{n-1/2}  \right>_Y 
- \frac{\Delta t^2}{4} \left< A \, A^* \, g^{n+1/2} \,,\,
   g^{n+3/2} -g^{n-1/2} \right>_Y \,.
\end{align*}
\begin{align*}
C_2^{n+1}- C_2^{n} 
& = \left< f^{n+1}-f^{n} \,,\, f^{n+1}+f^{n} \right>_X ; \\
& = \left< \Delta t \,  A^* \, g^{n+1/2} \,,\, f^{n+1}+f^{n} \right>_X ; \\
& = \Delta t \left< g^{n+1/2} \,,\, A f^{n+1}+A f^{n} \right>_Y ; \\
& = \Delta t \left< g^{n+1/2} \,,\,
	-\frac{g^{n+3/2}- g^{n-1/2}}{\Delta t} \right>_Y ; \\
& = - \left< g^{n+1/2} \,,\,
	g^{n+3/2}- g^{n-1/2} \right>_Y  .
\end{align*}
\begin{align*}
C_3^{n+1}- C_3^{n}
& = \left<
   A f^{n+1}- A f^n  \,,\,  A f^{n+1}+ A f^n
 \right>_Y ; \\
& = \left<
   A f^{n+1}- A f^n  \,,\, - \frac{g^{n+3/2}-g^{n-1/2}}{\Delta t}
 \right>_Y ; \\
& = - \left<
   \Delta t \, A A^* \, g^{n+1/2} \,,\, \frac{g^{n+3/2}-g^{n-1/2}}{\Delta t}
 \right>_Y ; \\
& = - \left<
   A A^* \, g^{n+1/2} \,,\, g^{n+3/2}-g^{n-1/2} \right> \,.
\end{align*}

Consequently
\begin{equation} \label{Conserved A full}
C^n = C_1^n + C_2^n - \left( \frac{\Delta t}{2}\right) C_3^n \,.
\end{equation}
is a conserved quantity. Moreover, if $\norm{A}$ is the operator norm of $A$
then 
\begin{equation}
\norm{C^n} \geq \left(1 - \frac{\Delta t^2}{4} \norm{A}^2\right) \norm{f^n}^2
	+ \norm{\frac{g^{n+1/2} + g^{n-1/2}}{2}}^2 \,,
\end{equation}
so $C^{n}$ is positive for $\Delta t \leq 2/\norm{A}$.

Next let
\begin{align*}
C_1^{n+1/2} & = \norm{\frac{f^{n+1} + f^{n}}{2}}_X^2 \,, \\
C_2^{n+1/2} & = \norm{ g^{n+1/2}}_Y^2 \,, \\
C_3^{n+1/2} & = \norm{A^* \, g^{n+1/2} }_X^2  \,.
\end{align*}
Now compute the changes in $C_1$, $C_2$ and $C_3$ over a time step:
\begin{align*}
C_1^{n+1/2} - C_1^{n-1/2}
& = \left< \frac{f^{n+1}+f^{n}}{2} , \frac{f^{n+1}+f^{n}}{2} \right>_X
  - \left< \frac{f^{n}+f^{n-1}}{2} , \frac{f^{n}+f^{n-2}}{2} \right>_X ; \\
& = \left<
\frac{f^{n+1} + 2 \, f^{n}  + f^{n-1}}{4} , f^{n+1} - f^{n-1}
\right>_X ; \\
& = \left<
f^{n}  - \frac{\Delta t^2}{4} A^* \, A f^{n} , f^{n+1} - f^{n-1}
\right>_X ; \\
& = \left<
f^{n} , f^{n+1} - f^{n-1} \right>_X -
\frac{\Delta t^2}{4}
\left< A^* \, A f^{n} , f^{n+1} - f^{n-1} \right>_X .
\end{align*}
\begin{align*}
C_2^{n+1/2}- C_2^{n-1/2}
& = \left< g^{n+1/2}, g^{n+1/2} \right>_Y
  - \left< g^{n-1/2}, g^{n-1/2} \right>_Y \\
& = \left< g^{n+1/2} + g^{n-1/2} , g^{n+1/2} - g^{n-1/2}  \right>_Y ; \\
& = \left< g^{n+1/2} + g^{n-1/2} , \Delta t \, A f^n  \right>_Y ; \\
& = -\Delta t \, \left< -A^* \, g^{n+1/2} - A^* \, g^{n-1/2} , f^n \right>_X ; \\
& = -\Delta t \, \left< \frac{f^{n+1}-f^{n-1}}{\Delta t} , f^n \right>_X ; \\
& = - \left< f^n , f^{n+1}-f^{n-1} \right>_X \,.
\end{align*}
\begin{align*}
C_3^{n+1/2}- C_3^{n-1/2}
& = \left<
A^* \, g^{n+1/2} - A^* \, g^{n-1/2} \,,\, A^* \, g^{n+1/2} + A^* \, g^{n-1/2}
\right>_X ; \\
& = \left<
A^* \left(g^{n+1/2} - g^{n-1/2}\right) \,,\,
-\frac{f^{n+1}-f^n}{\Delta t}-\frac{f^n-f^{n-1}}{\Delta t}
\right>_X ; \\
& = \left< \Delta t \, A^* \, A f^{n} \,,\,
	- \frac{f^{n+1}-f^{n-1}}{\Delta t} \right>_X ; \\
& = - \left< A^* \, A f^{n} \,,\, f^{n+1}-f^{n-1} \right>_X \,.
\end{align*}
Consequently
\begin{equation} \label{Conserved A* half}
C^{n+1/2} = C^{n+1/2}_1 + C^{n+1/2}_2 -
\left(\frac{\Delta t}{2}\right)^2 \, C^{n+1/2}_3 
\end{equation}
is a conserved quantity. Moreover, if $\norm{A^*}$ is the operator norm
of $A^*$ then
\[
C^{n+1/2} \geq \norm{\frac{f^{n+1} + f^{n}}{2}}_X^2 +
\left( 1 - \frac{\Delta t^2}{4} \norm{A^*}^2 \right)
\norm{ g^{n+1/2}}_Y^2 , 
\]
so $C^{n+1/2}$ is positive for $\Delta t \leq 2/\norm{A^*}$.

The program {\tt SystemODEs.m} confirms that the solutions of the system are
second order accurate and that the conserved quantities are second order
accurate and are constant to a small multiple of {\tt eps}.

\subsection{Notes}

Both $A\,A^*$ and $A^* \, A$ are self-adjoint positive operators,
but may not be positive definite.  So if $h \neq 0$ and $A \, h = 0$
then $ g(t) = t\, h$ is an unbounded solution of the second second order
equation while if $A^*h = 0$ then $f(t) = t \, h$ is an unbounded solution
of the first second order equation. For this $f(t)$ the system becomes
$h = A \, g(t) \,,\, g'(t) = 0  \,.  $ So $g(t) = k$ a
constant and then
$ \left< h , h \right>
= \left< h , A \, k \right>
= \left< A^* h ,  k \right> =  0 $,
that is $h=0$ and then $f(t) = 0$ and $g(t) = k$ and $A \, k = 0$.  So the
unbounded solution of the second order equation is not a solution of the
system.  So the second order equations and the system are not exactly
consistent.

There is also a problem with the initial conditions for the system
and the second order equations. If $A$ is an $n$ by $m$ matrix, then
$A^*$ is $m$ by $n$ matrix and consequently $A \, A^*$ is an $n$ by
$n$ matrix and $ A^* \, A$ is an $m$ by $m$ matrix.  So the first of
the second order equation needs $2\,n$ initial conditions, and the second
of the second-order equations needs $2\,m$ initial conditions. The
system needs $n+m$ initial conditions.
However, for example, if one knows $f(t)$ then $g(t)$ can be found using
simple integration and the initial condition for $g(t)$ and conversely for
knowing $g(t)$.

If $n=m$ then the number of initial conditions are the same for all three
variants of the ordinary differential equations.  If $A$ is inevitable then
the system and second order equations are consistent. Unfortunately, the
$n \neq m$ is far more analogous to the situation for equations modeling
waves in spaces of dimension 1 or larger.

\newpage \clearpage
\setcounter{equation}{0}
\section{Discretizing the One Dimensional Wave Equation\label{1D Wave}}

\begin{table}[ht]
\begin{center}
\begin{tabular}{|l|c|l|}
\hline
quantity & units & name \\
\hline
$x$    & $d$       & spatial position  \\
$dx$   & $d$       & spatial increment \\
$L$    & $d$       & end of spatial interval \\
$t$    & $t$       & time              \\
$u$    & $d$       & displacement      \\
$v$    & $d/t$     & velocity          \\
$u_t$  & $d/t$     & velocity          \\
$u_x$  & $1$       & slope             \\
$\rho$ & $1/d$     & 1D density           \\
$\tau$ & $d/t^2$   & 1D Young's Modulus\\
$E$    & $d^2/t^2$ & energy            \\
\hline
\end{tabular}
\caption{Quantities and their one dimensional space-time units.}
\label{Units}
\end{center}
\end{table}

The physical model used to motivate this discussion is piece of string or wire
pinned at both ends that is under tension and whose mass and strength varies
along the string or wire.  To better understand the correct form of the one
dimensional wave equation in such materials, three forms of this equation will
be derived from the conservation of energy. The three forms are two second
order equations and most importantly a first order system of two equations.
An important tool for obtaining correct equations is to use the one dimensional
units for distance $d$ and time $t$ units of the variables in the equations
that are given in Table \ref{Units} so that the equations are dimensionally
consistent.

The discussion of the continuum problem begins by introducing the energy of
the string using this to derive the second order wave equation the describes
the motion of the string. The energy then suggests how to write the second
order equation as a system of two first order equations and also how to
introduce two weighted inner products so the system has the same
form as the system of ODEs studied in Section \ref{ODEs}. The results
in Section \ref{ODEs} then motivate the discretization of the system
so that the discrete system has conserved quantities.

For spatially dependent differential equations two case will be studied:
CMP standing for constant material properties and VMP standing for
variable material properties. The simulation codes are created first
for the CMP case and then extended to the VMP case.  For the VMP case two
functions are introduced to represent the material properties, one is the one
dimensional density $\rho(x)$ while the other is the one dimensional Young's
modulus $\tau(x)$ that depends on the strength of the material and the force
applied to the string.  Two inner products weighted by $\rho$ and $\tau$ and
their related norms are also introduced and used to show that this system of
equations again has the same form as the system studied in Section \ref{ODEs}.
This is then used to introduce conserved quantities that are defined in terms
of the weighted norms.  The conserved quantities have a simple relationship to
the energy that was used to derive the differential equations and thus imply
the conservation of energy. For the CMP case it is assumed that $\rho$ and
$\tau$ are constant and then the wave speed $c^2 = \tau/\rho$ is
introduced.

The discussion of the discretizations begins with constant material properties
CMP followed by variable material properties VMP.  Two reason: the CMP case is
what is commonly studied in the literature and is relatively simple;
the VMP case requires the introduction of weighted inner products and is more
complicated.   More importantly, the CMP as done in the literature doesn't
provide a good motivation for VMP case.  In both cases a variable $v$ is
introduced to define the first order system, however the units of $v$ in the
two cases are different.  In the CMP case the wave speed is introduced which
is commonly used in text books. In this case $v$ has units distance
$v \sim d$.  In the VMP case the material is characterized by the one
dimensions density $\rho$ and the one dimensional Young's modulus $\tau$.
In this case a natural first order system is obtained by introducing a
velocity $v$ with units $v \sim d/t$.

The discretization of the system begins by introducing a grid staggered in
both space and time and the using central differences in both space and time
so that the discretizations are second order accurate.  The time discretization
is the same as leapfrog discretization used in Section \ref{ODEs} while the
spatial discretization is the 3D mimetic discretization
\cite{RobidouxSteinberg2011} specialized to one dimension. In this case
the spatial discretization is also a simple staggered discretization.
The discretization of the system also gives a standard discretization of the
two second order wave equations.  The resulting discrete system can be put in
the form used to study system of ODEs studied in the section \ref{ODEs} and
consequently the results in that section provide two discrete conserved
quantities.

The simulation codes {\tt Wave1DCMP.m} and {\tt Wave1DVMP.m} test the
developed theory for the constant and variable material properties cases.
For both cases the simulation programs show that the discrete solution
converges at least order 2 to the continuum solution and that conserved
quantities are second order accurate and constant to within small multiple
of machine epsilon ${\tt eps}$ that increases slowly with decreasing time and
space step sizes.  Surprisingly, for the CMP code if the final time is a
multiple of half of the period of the solution then the convergence rate is
order 4.  For the CMP code it is important that the condition that
guarantees the conserved quantities are positive is the same as the
Courant-Friedrichs-Lewy (CFL) condition for stability.
This condition guarantees that the conserved quantity is positive and
consequently the discretization is stable.

Next the VMP equations are discretized in the same manner as the CMP equations
and again the results of Section \ref{ODEs} can be applied.  As in the
continuum, weighted inner products and norms of discrete functions are
introduced and used to define discrete weighted conservation laws. Again,
these conserved quantities are positive for constraints on the time step that
are a generalization of the CFL condition. The discretization is stable
when the conserved quantities are positive.

For constant material properties the simulation program {\tt Wave1DVMP.m} will
produce the same results as the constant properties code {\tt Wave1DCMP.m} if
$\rho = 1/c$ and $\tau = c$. Both codes show that two conserved quantities
are constant and second order accurate. For the CMP case it is easy to find
exact solutions for wave equation and thus estimate the error in the solution
and the convergence rate.  For the general variable coefficient case it is not
possible to find exact solutions of the wave equation so an alternative method
of estimating the convergence rate is used.

For the VMP case the convergence rates is computed by comparing the discrete
solution on a given grid to the solution on a grid where each cell is cut in
half and the time step is also cut in half. The accuracy of the solution on
the original grid is estimated by the difference of the two solutions on
the course grid. The convergence rate is easily estimated but is a bit noisy.
It is far more informative divide the error estimate by the square of the size
of the cells in the course grid and plot the results. If the plots overlap
with decreasing step size then this shows that the solutions converge at least
order 2. Most example are order 2 but for some examples the plots increase
indicating that the convergence rate order is less than 2. In all cases the
convergence rate order is at least 1.  There are a several examples given in
{\tt Wave1DVMP}.
The "method of manufactured solutions" could be used to compute the order
of convergence of trivially modified equations but that was not done here.

\subsection{Derivation of the VMP 1D Wave Equation}
Consider a piece of string that is pinned at both ends and whose displacement
from equilibrium is given by $u(x,t)$ for $0 \leq x \leq 1$ and $t \geq 0$.
The boundary conditions are then $u(0,t) = u(1,t) = 0$.
The derivatives of functions will be given by subscripts:
$u_t = \partial u / \partial t$;
$u_x = \partial u / \partial x$;
$u_{tt} = \partial^2 u / \partial t^2$;
$u_{xx} = \partial^2 u / \partial x^2$.
The derivation of the variable material properties wave equation will start
with the conservation of the energy of a vibrating string which is the average
of the kinetic and potential energies,
\[
E = \frac{1}{2} \int_{0}^{1} \, \rho \,  u_t^2 \, dx
+ \frac{1}{2} \int_{0}^{1} \tau u_x^2 \, dx  \,.
\]
Here $\rho = \rho(x)$ is the one dimensional density and $\tau = \tau(x)$ is
the one dimensional Young's modulus.
It is assumed $\rho$ and $\tau$ are bounded above and below by positive
constants.  Table \ref{Units} shows that the energy $E$ has space and time
units $d^2/t^2$.

The wave equation is derived by assuming the energy $E$ is conserved and
the using integration by parts. First
\begin{equation}
\frac{d E}{d t} = \int_{0}^{1}
\left( \rho \, u_t \,u_{tt} + \tau \, u_x \, u_{xt} \right) \, dx \,.
\end{equation}
However
\begin{align*}
\int_0^1 (\tau \, u_x \, u_t)_x \, dx & =
\int_0^1 (\tau \, u_x)_x \, u_t \, dx +
\int_0^1 \tau \, u_x \, u_{xt} \, dx \\
& =  \tau(1) \, u_x(1,t) \, u_t(1,t) - \tau (0) \, u_x(0,t) \, u_t(0,t) \,.
\end{align*}
Consequently
\[
\frac{d E}{d t} = 
\int_{0}^{1}
u_t \left( \rho \, u_{tt} - \left( \tau \, u_{x} \right)_x \right) \, dx
+ \tau(1) \, u_x(1,t) \, u_t(1,t) - \tau(0) \, u_x(0,t) \, u_t(0,t) \,.
\]
But $ u(0,t) = u(1,t) = 0 $ for $t \geq 0$ so $ u_t(0,t) = u_t(1,t) = 0 $
and thus the energy will be conserved provided that
\[
\rho \, u_{tt} = \left( \tau \, u_{x} \right)_x \,,
\]
which is the 1D wave equation with variable material properties.

It is also common to derive the wave equation from Newton's law,
which gives the mass times acceleration is equal to the force which
is the form of this equation.  The form of the equation used here is
\begin{equation}\label{1D VMP Wave Equation}
u_{tt} = \frac{1}{\rho} \, \left( \tau \, u_{x} \right)_x \,.
\end{equation}
In \cite{Dawkins2020} the same form is given for the wave equation but
there $\tau$ is a function both $x$ and $t$.

This second order wave equation can be written as a system by introducing
a velocity $v = v(x,t)$ with units $d/t$ where
\begin{equation} \label{1D VMP System}
u_t = \frac{1}{\rho} \,  v_x \,, \quad
v_t = \tau \, u_{x} \,.
\end{equation}
Also $v$ satisfies a second order differential equation
\[
v_{tt} = \tau \left(\frac{1}{\rho} \,  v_x \right)_x \,.
\]
The second order equation can be expanded to
\begin{equation}\label{Second Order Equations}
u_{tt} = \frac{\tau}{\rho} u_{xx} + \frac{\tau_x}{\rho} u_x \,,\quad
v_{tt} = \frac{\tau}{\rho} v_{xx} - \frac{\tau \rho_x}{\rho^2} v_x \,.
\end{equation}
This brings up a concern that if either $\rho_x/\rho$ or $\tau_x/\tau$ are
large then the first order terms involving $u_x$ or $v_x$ will dominate the
second order terms and then the solution second order equations and
consequently the first order equations could have non wave like behavior?

To apply the ideas in Section \ref{ODEs} to the 1D wave equation two inner
products of functions defined on $[0,1]$ are needed. 
First, if $u1 = u1(x)$ and $u2 = u2(x)$ have dimensions $d$ then let
\begin{equation} \label{rho inner product}
\left< u1 , u2 \right>_\rho = \int_{0}^{1}
u1(x) \, u2(x) \rho(x) \, dx \,,
\end{equation}
Note that $\rho \, dx$ has units 1 so the inner product has units $d^2$.
Next, if $v1 = v1(x)$ and $v2 = v2(x)$ have units $d/t$ then let
\begin{equation} \label{tau inner product}
\left< v1 , v2 \right>_\tau = \int_{0}^{1}
v1(x) \, v2(x) \tau^{-1}(x) \, dx \,.
\end{equation}
Again, note that $\tau^{-1}(x) \, dx$ has units $t^2$ and then because $v1$
and $v2$ have units $d/t$ the inner product has units $d^2$.

To apply the ideas in Section \ref{ODEs} define the operator $A$ as
\[
A u = \tau \, u_x \,.
\]
If $u(0) = u(1) = 0$ then integration by parts gives
\begin{align}
\left< A u , v \right>_\tau
& = \left< \tau \, u_x \, v \right>_\tau \nonumber \\
& = \int_{0}^{1} \tau \, u_x \, v \, \tau^{-1} \, dx \,, \nonumber\\
& = \int_{0}^{1} u_x \, v \, dx \,, \nonumber\\
& =-\int_{0}^{1} u \, v_x \, dx \,, \nonumber\\
& =-\int_{0}^{1} u \, \frac{1}{\rho} \, v_x \, \rho \, dx \,, \nonumber\\
& = \left< u \, A^* v \right>_\rho \,. 
\end{align}
Consequently the adjoint of $A$ is given by
\[
A^*v = - \frac{1}{\rho} v_x \,.
\]
Now the system \eqref{1D VMP System} can be written in the form
given in \eqref{Wave System}:
\begin{equation}
u_t = - A^* \, v \,;\quad v_t = A \, u  \,.
\end{equation}
Additionally $u$ and $v$ will satisfy the second order equations
\begin{equation} \label{2ndOWE}
u_{tt} = - A^* \, A u \,,\quad v_{tt} = - A \, A^* v \,.
\end{equation}

Next
\[
\left< A \, A^* v1 , v2 \right>_\tau = 
 \left< A^* v1 , A^* v2 \right>_\rho = 
  \left< v1 , A \, A^* v2 \right>_\tau  \,,
\]
and
\[
\left< A^* \, A u1 , u2 \right>_\rho = 
 \left< A u1 , A u2 \right>_\tau = 
  \left< u1 , A^* \, A u2 \right>_\rho  \,.
\]
Consequently both $A \, A^*$ and $A^* \, A$ are self adjoint and
positive so that \eqref{2ndOWE} are wave equations.

As in Section \ref{ODEs} a conserved quantity is given by
\[
C(t) = \frac{1}{2} \, \norm{u}_{\rho}^2 +
       \frac{1}{2} \, \norm{v} _{\tau}^2 \,
\]
because if $u(0) = u(1) = 0$ then
\begin{align*}
C_t & = \left< u , u_t \right>_\rho + \left< v , v_t \right>_\tau  \\
    & = \int_{0}^{1}
        u(x) \, \frac{1}{\rho(x)} \, v_x (x) \, \rho(x) \, dx +
       \int_{0}^{1}
        v(x) \, \tau(x) \, u_x (x)  \, \frac{1}{\tau(x)} dx \\
    & = \int_{0}^{1} u(x) \, v_x (x) + v(x) \, u_x (x)  dx \\
    & = \int_{0}^{1} \left( u(x) \, v (x) \right)_x   dx \\
    & = 0 \,.
\end{align*}
Now note that $u_t$ and $v_t$ are also a solution of the system
\eqref{1D VMP System} implying that 
\[
E = \frac{1}{2} \, \norm{u_t}_{\rho}^2 + \frac{1}{2} \, \norm{v_t}_{\tau}^2
\]
is constant. Also the units of $E$ are $d^2/t^2$. A little algebra shows
that $E$ is the energy of the second order 1D wave equation.

\subsection{The Constant Material Properties Wave Equation}

The following standard formulation of the constant material properties (CMP)
1D wave equation provides insight into the variable material properties case.
The wave equation in constant materials is different from just assuming that
$\rho$ and $\tau$ are constant in the above discussion because of the
introduction of the wave speed $c$.

If $\rho$ and $\tau$ are constant then \eqref{1D VMP Wave Equation} becomes
\begin{equation}
u_{tt} = \frac{\tau}{\rho} \, u_{xx} \,.
\end{equation}
Note that $\tau / \rho$ has dimension $d^2/t^2$.
For constant material properties the 1D wave equation is typically
given as
\[
u_{tt} = c^2 \, u_{xx} \,,
\]
where $c>0$ is the constant wave speed with units $d/t$. So
\[
c^2 = \frac{\tau}{\rho} \,.
\]
This equation can also be written as a system
\begin{equation} \label{1D CMP System}
u_t = c \, v_x \,,\quad v_t = c \, u_x \,,
\end{equation}
but now $v$ must have units $d$ not $d/t$ as will be the case for variable
material properties, so $v$ is not a velocity. So this system is not the same
as the first order system for variable material properties that will be
described below. Note that interchanging $u$ and $v$ doesn't change the system.
As before $v$ also satisfies a second order wave equation
\[
v_{tt} = c^2 \, v_{xx} \,.
\]

It is easy to check that for $u(0,t) = u(1,t) = 0$
\begin{equation} \label{CMP conserved quantity}
C = \frac{1}{2} \int_{0}^{1} \left( u^2 + v^2 \right) dx \,,
\end{equation}
is a conserved quantity.  Again note that if $u,v$ are solutions of the
system then so are $u_t$, $v_t$ and then a conserved quantity is given by 
\[
E = \frac{1}{2} \int_{0}^{1} \left( u_t^2 + v_t^2 \right) \rho dx = 
\frac{1}{2} \int_{0}^{1} \left( u_t^2 + c^2 u_x^2 \right) \rho dx \,.
\]
The constant $\rho$ is included so that $E$ has dimensions $d^2/t^2$ which
is correct for energy.  So if $C$ is conserved then so is $E$.

The constant material properties simulation code {\tt Wave1DCMP.m}
will be tested for $t \geq 0$ and the interval $0 \leq x \leq 1$ using
the solution
\begin{align*}
  u(x,t) & = \cos(m \, \pi \, c \, t) \, \sin(m \, \pi \, x), \\
  v(x,t) & = \sin(m \, \pi \, c \, t) \, \cos(m \, \pi \, x).
\end{align*}
where $m$ is a positive integers and $c$ is a positive constant.
Note that the arguments of the $sin$ and $cos$ should be dimensionless
but they are not. This can be fixed by multiplying the arguments by the
density $\rho$ of the string and then assuming $\rho = 1$.
Note that $u(0,t) = u(1,t) = 0$ so this solution satisfies homogeneous
Dirichlet boundary conditions while $v$ satisfies homogeneous Neumann
boundary conditions.  Also, interchanging $u$ and $v$ gives a solution of the
second order equation that satisfies Neumann boundary conditions.  For the
test solution the quantity \eqref{CMP conserved quantity} is $C(t) = 1/2$.
An important topic not included is general boundary conditions for the wave
equation.

\subsection{Staggered Discretizations of the 1D Wave Equations}

\begin{figure}
\begin{center}
   \includegraphics[width=5.00in,trim = 0 320 0 50, clip ]{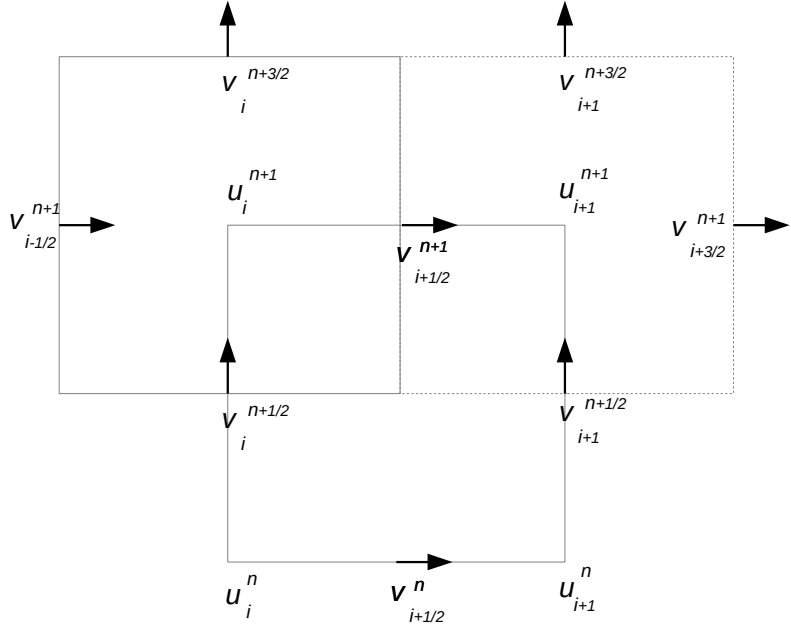}
\caption{Space-Time Staggered Grid}
\label{TimeSpaceGrid}
\end{center}
\end{figure}

As illustrated in Figure \ref{TimeSpaceGrid} two grids that are staggered in
space and time grids will be used to discretize the first order constant
materials properties (CMP) system \eqref{1D CMP System} and the variable
material properties (VMP) system \eqref{1D VMP System}.  The simulation region
is $a \leq x \leq b$ and $ 0 \leq t \leq T$ where $L = b-a > 0$ and $T > 0$.
One grid is the primal grid and the other is the dual grid.  For the primal
grid $Nx$ is the number of discretization points in the spatial region.
Consequently the number of cells in the spatial grid is $Nx-1$ so set
$\Delta x = L/(Nx-1)$.  If $Nt$ is the number of time steps then
$\Delta t = T/Nt$.  Then the primal grids points are given by
\begin{equation} \label{Math Primal Grid}
\left(t^n,x_i\right) = \left(n\,\dt, a+i \, \dx \right) \,,\quad 
0 \leq n \leq Nt ,\, 0 \leq i \leq Nx-1 \,,
\end{equation}
The dual grid points are at the centers of the cell in the primal grid,
\begin{equation} \label{Math Dual Grid}
\left(t^{n+1/2},x_{i+1/2}\right) = 
	\left((n+1/2)\,\dt, a+(i+1/2) \, \dx \right) \,,\quad
0 \leq n \leq Nt ,\, 0 \leq i \leq Nx-2 \,.
\end{equation}
The discretization of $u(x,t)$ on the primal grid the discretization of
$v(x,t)$ on the dual grid are: 
\begin{equation} \label{Math Discretization u and v}
u^n_i = u(x_i,t^n) \,;\quad
v^{n+1/2}_{i+1/2} = v(x_{i+1/2}, t^{n+1/2}) \,.
\end{equation}

\subsection{Constant Material Proprieties (CMP)}

For CMP the inner products and norms of the discrete functions are given by:
\[
\left< u1, u2 \right>
= \sum_{i=0}^{Nx-1} u1_i \, u2_i \, \Delta x :\quad
\left< v1, v2 \right>
= \sum_{i=0}^{Nx-2} v1_{i+1/2} \, v2_{i+1/2} \, \Delta x :
\]
\[
\norm{u}^2 = \left< u \, u \right> :\quad
\norm{v}^2 = \left< v \, v \right> \,.
\]
To have the units of the norms be $d^2$, because $u$ and $v$ have units $d$,
the norms can be multiplied by the constant density $\rho$ so that
$\rho \, dx$ is dimensionless, and then assume that $\rho = 1$.

The CMP system \eqref{1D CMP System} is discretized as
\begin{align}\label{1D UPM Discretized System}
\frac{u^{n+1}_{i} - u^{n}_{i}}{\dt}
	& = c \, \frac{v^{n+1/2}_{i+1/2} - v^{n+1/2}_{i-1/2}}{\dx} \,,
	\nonumber \\
\frac{v^{n+1/2}_{i+1/2} - v^{n-1/2}_{i+1/2}}{\dt} 
	& = c \, \frac{u^n_{i+1}-u^n_i}{\dx} \,.
\end{align}
The initial conditions are $u^0_n$ and $v^{1/2}_{n+1/2}$.

Introduce the finite difference operators
\[
\delta(u)_{i+1/2} = u_{i+1}-u_i \,,\quad \delta(v)_{i} = v_{i+1/2}-v_{i-1/2} \,,
\]
so that the system can be written
\begin{equation}
\frac{u^{n+1}_{i} - u^{n}_{i}}{\dt}
	 = c \, \frac{\delta(v^{n+1/2})_i}{\dx} \,,\quad
\frac{v^{n+1/2}_{i+1/2} - v^{n-1/2}_{i+1/2}}{\dt} 
	 = c \, \frac{\delta(u^n)_{i+1/2}}{\dx} \,.
\label{CMP Discretized System}
\end{equation}
To use the discussion in Section \ref{ODEs} this system must be put into the
form \eqref{Staggered Time Equations}:
\begin{equation}  \label{Discrete half system}
\frac{u^{n+1} - u^n}{\Delta t} = -A^* v^{n+1/2}  \,,\quad
\frac{v^{n+1/2} - v^{n-1/2}}{\Delta t} = A u^{n} \,.
\end{equation}
So
\[ 
A  = c \frac{\delta}{\Delta x}  \,,\quad
A^* = - c \frac{\delta}{\Delta x} .
\]
It is easy to check that $A^*$ is the adjoint operator of $A$.

Now \eqref{Conserved A full} and \eqref{Conserved A* half} give the conserved
quantities
\begin{align}
C^n
& = \norm{u^n}^2 +
\norm{\frac{v^{n+1/2}+v^{n-1/2}}{2}}^2 
-\left(\frac{\dt}{2}\right)^2 \norm{A \, u^n }^2 \,,
\nonumber \\
& = \norm{u^n}^2 +
\norm{\frac{v^{n+1/2}+v^{n-1/2}}{2}}^2 
-\left(\frac{c \, \dt}{2 \, \dx}\right)^2 \norm{\delta  u^n }^2\,,
\end{align}
and 
\begin{align}
C^{n+1/2}
& =
  \norm{\frac{u^{n+1} + u^{n}}{2}}^2
+ \norm{ v^{n+1/2}}^2
-\left(\frac{\dt}{2}\right)^2 \norm{A^* \, v^{n+1/2} }^2 \,,
\nonumber \\
& =
  \norm{\frac{u^{n+1} + u^{n}}{2}}^2
+ \norm{ v^{n+1/2}}^2
-\left(\frac{c \, \dt}{2 \, \dx}\right)^2 \norm{ \delta \, v^{n+1/2} }^2 \,.
\end{align}

The conserved quantities will be positive provided 
\[
\frac{c \, \dt}{2 \, \dx} \norm{\delta} < 1 \,.
\]
Because $\norm{\delta} = 2$ the conserved quantity will be positive if
\[
c \, \frac{\dt}{ \dx} < 1 \,.
\]
This is the Courant-Friedrichs-Lewy (CFL) condition for stability of the
second order discrete equations that are derived below from the first order
discrete system.

\subsection{Variable Material Properties}

For variable material properties $\rho$ is discretized on the primal
grid while $\tau$ is discretized on the dual grid:
\[
\rho_i = \rho(x_i)\,;\quad \tau_{i+1/2} = \tau(x_{i+1/2}) \,.
\]
Then the system \eqref{1D VMP System} is discretized as
\begin{align}\label{1D VMP Discretized System}
\frac{u^{n+1}_{i} - u^{n}_{i}}{\dt}
	& = \frac{1}{\rho_i} \,
\frac{v^{n+1/2}_{i+1/2} - v^{n+1/2}_{i-1/2}}{\dx} \,,
	\nonumber \\
\frac{v^{n+1/2}_{i+1/2} - v^{n-1/2}_{i+1/2}}{\dt} 
	& = \tau_{i+1/2} \, \frac{u^n_{i+1}-u^n_i}{\dx} \,.
\end{align}
As before assume that $u^0$ and $v^{\half}$ are given then the leapfrog time
stepping scheme for $n \geq 0$ is
\[
u^{n+1}_{i} = u^{n}_{i} +
\frac{1}{\rho_i} \, \frac{\dt}{\dx}
\left(v^{n+1/2}_{i+1/2} - v^{n+1/2}_{i-1/2}\right) \,,\quad
v^{n+3/2}_{i+1/2} = v^{n+1/2}_{i+1/2} + \tau_{i+1/2} \, \frac{\dt}{\dx} \,
\left(u^{n+1}_{i+1}-u^{n+1}_i\right) \,.
\]

To use the discussion in Section \ref{ODEs} to study the discretized wave
equation, the system \eqref{1D VMP Discretized System} must be put into the
form \eqref{Staggered Time Equations}:
\begin{equation} 
\frac{u^{n+1} - u^n}{\Delta t} = - A^* \, v^{n+1/2}  \,;\quad
\frac{v^{n+1/2} - v^{n-1/2}}{\Delta t} = A \, u^{n} \,.
\end{equation}
The difference equations \eqref{1D VMP Discretized System} can be written as
\begin{equation} \label{Delta Discretized System}
\frac{u^{n+1}_{i} - u^{n}_{i}}{\dt}
	 = \frac{1}{\rho_i} \, \frac{\delta(v^{n+1/2})_{i}}{\dx} \,,\quad
\frac{v^{n+1/2}_{i+\half} - v^{n-1/2}_{i+1/2}}{\dt} 
	 = \tau_{i+\half} \, \frac{\delta(u^n)_{i+1/2}}{\dx} \,.
\end{equation}
Consequently 
\[ 
A^* = - \frac{\delta}{\rho \, \Delta x} \,,\quad
A =     \frac{\tau \, \delta}{\Delta x} \,.
\]
Below it will be shown that $A^*$ is the adjoint of $A$.

The second order variable material properties wave equation for $u$ is:
\begin{align}
\frac{u^{n+1}_{i} - u^{n}_{i}}{\dt} -
\frac{u^{n}_{i} - u^{n-1}_{i}}{\dt} & =
- A^* \, v^{n+1/2}_i - (-A^*) \, v^{n-1/2}_i
; \nonumber \\ 
\frac{u^{n+1}_{i} - u^{n}_{i}}{\dt} -
\frac{u^{n}_{i} - u^{n-1}_{i}}{\dt} & =
\frac{1}{\rho_i} \, \frac{\delta(v^{n+1/2})_{i}}{\dx}
-\frac{1}{\rho_i} \, \frac{\delta(v^{n-1/2})_{i}}{\dx}
; \nonumber \\
\frac{u^{n+1}_{i} - u^{n}_{i}}{\dt}
- \frac{u^{n}_{i} - u^{n-1}_{i}}{\dt} & =
\frac{1}{\rho_i} \, \frac{v^{n+1/2}_{i+1/2} - v^{n+1/2}_{i-1/2}}{\dx}
- \frac{1}{\rho_i} \, \frac{v^{n-1/2}_{i+1/2} - v^{n-1/2}_{i-1/2}}{\dx}
; \nonumber \\
\frac{u^{n+1}_{i} - 2 \, u^{n}_{i} + u^{n-1}_{i}}{\dt} & =
\frac{v^{n+1/2}_{i+1/2}-v^{n+1/2}_{i-1/2}-v^{n-1/2}_{i+1/2}+v^{n-1/2}_{i-1/2}}
{\dx \, \rho_i} 
; \nonumber \\
\frac{u^{n+1}_{i} - 2 \, u^{n}_{i} + u^{n-1}_{i}}{\dt^2} & =
\frac{1}{\dx  \, \rho_i}
\frac{v^{n+1/2}_{i+1/2} - v^{n-1/2}_{i+1/2} - \left( v^{n+1/2}_{i-1/2}
- v^{n-1/2}_{i-1/2} \right)}{\dt} 
; \nonumber \\
\frac{u^{n+1}_{i} - 2 \, u^{n}_{i} + u^{n-1}_{i}}{\dt^2} & =
\frac{1}{\dx  \, \rho_i}
\left(\tau_{i+1/2} \, \frac{u^n_{i+1}-u^n_i}{\dx}
- \tau_{i-1/2} \, \frac{u^n_{i}-u^n_{i-1}}{\dx} \right) 
; \nonumber \\
\frac{u^{n+1}_{i} - 2 \, u^{n}_{i} + u^{n-1}_{i}}{\dt^2} & =
\frac{1}{\rho_i}
\frac{
  \tau_{i+1/2} \, u^n_{i+1}
- ( \tau_{i+1/2} + \tau_{i-1/2} )  u^n_{i}
+ \tau_{i-1/2} \, u^n_{i-1}
}{\dx^2}
.\nonumber
\end{align}
A similar argument gives the second order equation for $v$ as
\begin{align}
\frac{v^{n+1/2}_{i+1/2} - 2 \, v^{n-1/2}_{i+1/2} - v^{n-3/2}_{i+1/2}}{\dt^2} 
& = \tau_{i+1/2} \,
\frac{ \rho_{i+1}^{-1}  v^{n-1/2}_{i+3/2}  
- (\rho_{i+1}^{-1} + \rho_{i}^{-1} ) v^{n-1/2}_{i+1/2}
+ \rho_{i}^{-1} v^{n-1/2}_{i-1/2}}{\dx^2}  \,.
\end{align}

To study the conservation properties of the discretized system
\eqref{1D VMP Discretized System} weighted inner products and norms analogous
to \eqref{rho inner product} and \eqref{tau inner product} will be needed.
So if $u1_i = u1(x_i)$ and $u2_i = u2(x_i)$ and $\rho_i = \rho(x_i)$ then let
\begin{equation} \label{discrete rho inner product}
\left< u1 , u2 \right>_\rho = \sum_{i=1}^{Nx}
u1_i \, u2_i \rho_i \, \Delta x \,.
\end{equation}
Next, if $v1_{i+1/2} = v1(x_{i+1/2})$ and $v2_{i+1/2} = v2(x_{i+1/2})$
then let
\begin{equation} \label{discrete tau inner product}
\left< v1 , v2 \right>_\tau = \sum_{i=0}^{Nx-1}
v1_{i+1/2} \, v2_{i+1/2} \, \tau^{-1}_{i+1/2} \, \Delta x \,.
\end{equation}
As above, let the inner products and norms with no subscript be the
standard mean square norm and the norm with the subscript $\infty$
be the maximum of the absolute value.

First, these inner products can be used to show that the adjoint of $A$
is $A^*$ that was defined above because
\begin{align*}
\left< A u , v \right>_\tau \
& = \left< \frac{\tau \, \delta u}{\Delta x} , v \right>_\tau \,, \\
& = \sum_0^{Nx-1} \frac{\tau_{i+\half} \, \delta u_{i+\half}}{\Delta x} \, 
	v_{i+\half} \, \tau_{i+\half}^{-1} \Delta x \,, \\
& = \sum_0^{Nx-1} {\delta u}_{i+\half} \, v_{i+\half} \,, \\
& = - \sum u_i \, {\delta v}_i \,, \\
& = - \sum u_i \, \frac{{\delta v}_i}{{\rho_i} \, \Delta x }
      \, \rho_i \, \Delta x , \\
& = \left< u , - \frac{\delta v}{\rho \, \Delta x} \right>_\rho \,, \\
& =  \left< u , A^* \, v \right>_\rho \,.
\end{align*}

The discrete wave equation has conserved quantities that can be
derived from those defined in
\eqref{Conserved A full}
and
which are
\eqref{Conserved A* half}
\begin{align*}
C^n
= &
\norm{u^{n}}_\rho^2 +
\norm{\frac{v^{n+1/2} + v^{n-1/2}}{2}}_\tau^2 -
\left(\frac{\Delta t}{2} \right)^2 \norm{ A \, u^{n} }_\tau^2  \,, \\
= &
\norm{u^{n}}_\rho^2 +
\norm{\frac{v^{n+1/2} + v^{n-1/2}}{2}}_\tau^2 -
\left(\frac{\Delta t}{2} \right)^2
\norm{ \frac{\tau \, \delta u^{n} }{\Delta x} }_\tau^2  \,. \\
\end{align*}
\begin{align*}
C^{n+1/2}
= &
\norm{ v^{n+1/2}}_\tau^2 
+ \norm{\frac{u^{n+1} + u^{n}}{2}}_\rho^2 
- \left(\frac{\Delta t}{2}\right)^2 \norm{ A^* \, v^{n+1/2} }_\rho^2  \,, \\
= &
\norm{ v^{n+1/2}}_\tau^2 
+ \norm{\frac{u^{n+1} + u^{n}}{2}}_\rho^2 
- \left(\frac{\Delta t}{2}\right)^2
\norm{ \frac{\delta v^{n+1/2}}{\rho \, \Delta x} }_\rho^2  \,. \\
\end{align*}

To estimate when the conserved quantity is positive compute
\begin{align*}
\norm{\frac{\tau \, \delta u^{n} }{\Delta x} }_\tau^2 = &
\left< \frac{\tau \, \delta u^{n} }{\Delta x} ,
\frac{\tau \, \delta u^{n} }{\Delta x} \right>_\tau \\
= & \left< \tau \, \frac{\delta u^{n} }{\Delta x} ,
\frac{\delta u^{n} }{\Delta x} \right>  \\
\leq & \frac{\max(\tau)}{\Delta x^2} \left<\delta u^n , \delta u^n \right> \\
\leq & \frac{ \max(\tau) }{\Delta x^2} \norm{\delta u^n}^2 \\
\leq & \max(\tau)\left(\frac{2}{\Delta x}\right)^2  \, \norm{u^n}^2 \\
\leq & \max(\tau)\left(\frac{2}{\Delta x}\right)^2  \,
	\left< u^n , u^n \right> \\
\leq & \max(\tau) \left(\frac{2}{\Delta x}\right)^2  \,
	\left< u^n , \frac{u^n}{\rho} \right>_\rho \\
\leq & \max(\tau) \left(\frac{2}{\Delta x}\right)^2  \,
	{\frac{1}{\min(\rho)}} \left< u^n , u^n \right>_\rho \\
\leq & \frac{\max(\tau)}{\min(\rho)}
	\left(\frac{2}{\Delta x}\right)^2 \norm{u^n}_\rho^2 \,.
\end{align*}
As noted above when $\rho$ and $\tau$ are constant the wave speed $c$ is given
by $c^2 = {\tau}/{\rho}$ so in the variable materials case
\begin{equation}\label{wave speed}
s = \sqrt{\frac{\max(\tau)}{\min(\rho)}}
\end{equation}
gives an estimate for the maximum wave speed.  Consequently
\[
C^n \geq 
\left(1 -
\left(s \, \frac{\Delta t}{\Delta x}\right)^2
\right)
\norm{u^n}_\rho^2
+ \norm{\frac{v^{n+1/2} + v^{n-1/2}}{2}}_\rho^2 \,.
\]

Similarly
\[
C^{n+1/2} \geq \norm{\frac{u^{n+1} + u^{n}}{2}}_\tau^2 +
\left( 1 - \left( s \, \frac{\Delta t}{\Delta x} \right)^2 \right)
\norm{ v^{n+1/2}}_\tau^2 \,.
\]

To assure that both conserved quantities are positive requires that
$ s \Delta t / \Delta x \leq 1$ which corresponds to the standard CFL
stability condition when $\rho$ and $\tau$ are constant. So for this
constraint the discretization is stable.
Note that if $\rho$ goes to zero or $\tau$ becomes large then $\Delta t$
must go to zero to maintain stability. So it is important for
density to be bounded below and the Young's modulus to be bounded above.

\subsection{Implementing the Discrete Systems}

Two codes will be created, {\tt Wave1DCMP.m} for constant material properties
and a more complex {\tt Wave1DVMP.m} for variable material properties. The
notation for the spatial derivatives is based on the notation in three
dimensions where $u$ will be a scalar while $v$ is a vector. So then
in the system \eqref{Discrete half system} $A$ must be the gradient while
$-A^*$ must be the divergence.  So the simulation codes will use the notation 
\begin{align*}
& {\tt CMP} && A = c \, {\tt Grad1}      &&
    A^* = - c \, {\tt Div1} \,, \\
& {\tt VMP} && A = {\tau} \, {\tt Grad1} &&
    A^* = - \frac{1}{\rho} \, {\tt Div1} \,.
\end{align*}

The spatial primal grid is obtained by replacing $i$ by $i-1$ in 
\ref{Math Primal Grid} to get
\[
xp(i) = a + (i-1) \, \Delta x \quad 1 \leq i \leq Nx \,.
\]
The dual spatial grid is obtained by replacing $i$ by $i-1$ in
\ref{Math Dual Grid} to get
\[
xd(i) = a+ (i-1/2) \, \Delta x  \,,\quad 1 \leq i \leq Nx-1 \,.
\]
Because $\Delta x = (b-a)/(Nx-1)$, $xp(1) = a$ and $xp(Nx) = b$.  Also note
that
\begin{align*}
xd(1) & = a +\Delta x /2 \, \\
xd(Nx-1)
& = a + (i-1-1/2) \dx  \\
& = a + (i-1)\, \dx - \dx/2 \\
& = a - (b-a) \dx / 2. \\
& = b - \dx / 2.
\end{align*}
The time grids are not stored in arrays so there in no need for a
special notation.

In \ref{Math Discretization u and v} the function $u(x,t)$ is discretized on
the primal grid while the function $v(x,t)$ discretized on the dual grid.
So in the simulation programs
\[
u^n_i \approx u(xp(i), n \, \dt)\,, \quad v^n_i \approx v(xd(i),(n+1/2) \dt)
 \,,\quad 0 \leq n \leq Nt \,.
\]
Note that for $n = Nt$, $t^{n+1/2} = (N_t+1/2) \, \dt = T + \dt/2$
and that this value of $v$ will be used to estimate $v$ at the final time
$T$ using
\[
v^{Nt}_{i+1/2} = \frac{v^{Nt+1/2}_{i+1/2}+ v^{Nt-1/2}_{i+1/2}}{2} \,.
\]

Both simulation codes will start with the initial data $u^{0} = u(x,0)$ and
$v^{0} = v(x,\Delta x/ 2)$.  Then the update equations for the simulation
code are obtained from \ref{Discrete half system}.  For $u$ replace $n+1/2$
by $n$ in the first equation while for $v$ replace $n$ by $n+1$ in
the second equation and then replace $n+3/2$ by $n+1$ and $n+1/2$ by $n$ to
obtain
\begin{equation} \label{Update Equations}
u^{n+1}_{i} = u^{n}_{i} + \dt \, c \, ({\tt Div1} v^n)_i \,,\quad
v^{n+1}_{i} = v^{n}_{i} + \dt \, c \, ({\tt Grad1} u^{n+1})_{i} \,,\quad
0 \leq i \leq Nt - 1 \,.
\end{equation}
Not that the second equation depends on the results of the first 
equation so the order of evaluation is critical.

The implementation of the code {\tt Wave1DVMP.m} is essentially the same as that
of {\tt Wave1DCMP.m}. Also see the simulation codes for the implementation of
the conserved quantities. Consequently a preliminary test for variable
materials simulation code is given by assuming the material properties $\rho$
and $\tau$ are constant and then use the analytic solution
\begin{align*}
  u(x,t) & = \cos(m \, \pi \, \frac{1}{\rho} \, t) \, \sin(m \, \pi \, x), \\
  v(x,t) & = \sin(m \, \pi \, \tau \, t) \, \cos(m \, \pi \, x).
\end{align*}
where $m$ is a positive integer.  Note that if $c = 1/\rho = \tau$ then this
is the same as the solution for the constant material code. In fact both
codes give the same result in this case.

\subsection{Simulation Results}

Simulations show that for {\tt Wave1DCMP.m} and {\tt Wave1DVMP.m} the
conserved quantities are constant within a small multiple of {\tt eps} and
that the constant increases slowly with increasing grid size.
The constant material code {\tt Wave1DCMP.m} produces at least second order
accurate solutions but computes the solution to forth order accuracy when the
final time is a multiple of the half period of the solution.
The situation for {\tt Wave1DVMP.m} is more complex as it does not alway
achieve second order accuracy but is at least first order accurate. A large
number of simulations are presented in an attempt to understand this.
The results are noisy especially when the errors are small. 

Because the simulation codes use centered differences it is to be expected
that the simulations should have order of convergence 2, that is, for
decreasing $\dx$ and $\dt$ the error in the solution should behave like
\[
Er = C_1 \dx^2 + C_2 \dt^2
\]
where $C_1$ and $C_2$ are constants.  This was tested by setting
$\dt = C_3 \, \dx$ for some constant $C_3$ and then for a sequence of
decreasing spatial step sizes $\dx$ the errors should behave like
$Er = C  \dx^2$ for some constant $C$. To estimate the order of
convergence, use two simulations with different values of $\dx$ and assume
that $Er_1 = C \, \dx_1^p$ and $Er_2 = C \, \dx_2^p$  and then $p$ can be
estimated by 
\[
p = \frac{\ln (Er_1) -\ln (Er_2)}{ \ln (\dx_1)-\ln( \dx_2)} \,.
\]

Because it is easy to find analytic solutions for the constant coefficient
wave equation the code {\tt Wave1DCMP.m} is tested by computing the errors
using an analytic solution.  For variable material properties, typically
there are no analytic solutions.  So in {\tt Wave1DVMP.m} the errors are
estimated by choosing the number of spatial grids point $Nx = 2^k+1$ for some
$k>0$. Increasing $k$ to $k+1$ will cause the refined grid points to be at
the edges and centers of the course grid cells. The difference between the
course grid solution and the fine grid solution at the course grid points
gives an estimate for the error in $u$ on the course grid points.
In the error analysis it more informative to choose a decreasing sequence of
$\dx_k$ and then estimate the error in the solution $Er_k(x_i)$ and then
plot $Er_k(x_i)/\dx_k^2$ as a function of $x_i$. It is important that this
illuminates where in space the errors are occurring.

For the time discretizations it is important for stability to keep the ratio
of $dx/dt$ constant. First if $Nx = 2^k + 1$ then $\dx = L/2^k$.  Next for
some integer $f$ choose $Nt = 2^{k+f}$ so that $\dt = T/2^{k+f}$ and then
$\dx/\dt = 2^{-f} \, L/T$ is constant for all $k$. The value of $f$ is chosen
so that the simulations are stable for all $k$.

\begin{figure}
\begin{center}
\begin{tabular}{cc}
A & B \\
\includegraphics[width = 3.0in, height = 2.5in, trim = 0mm 0mm 0mm 0mm, clip]{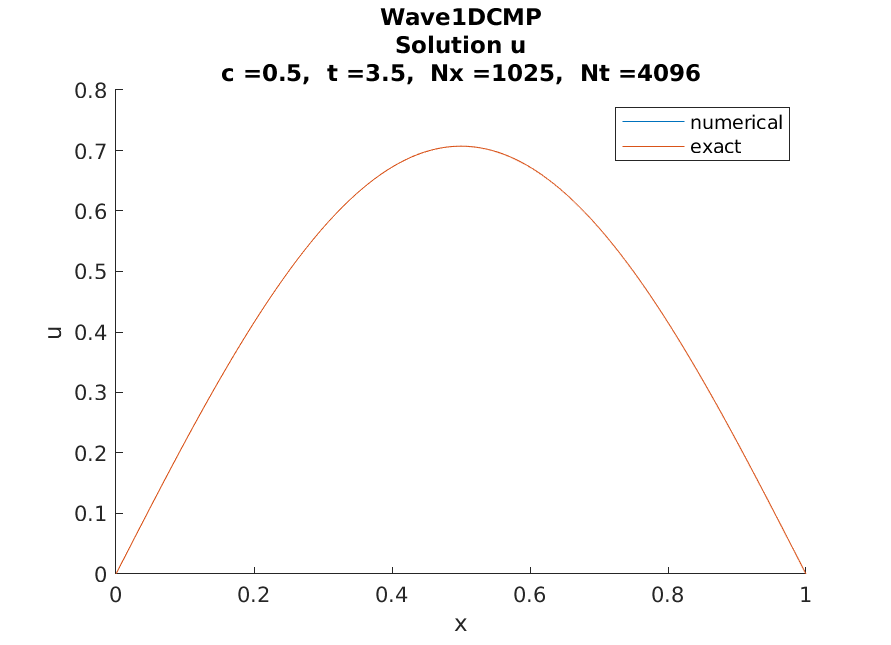} &
\includegraphics[width = 3.0in, height = 2.5in, trim = 0mm 0mm 0mm 0mm, clip]{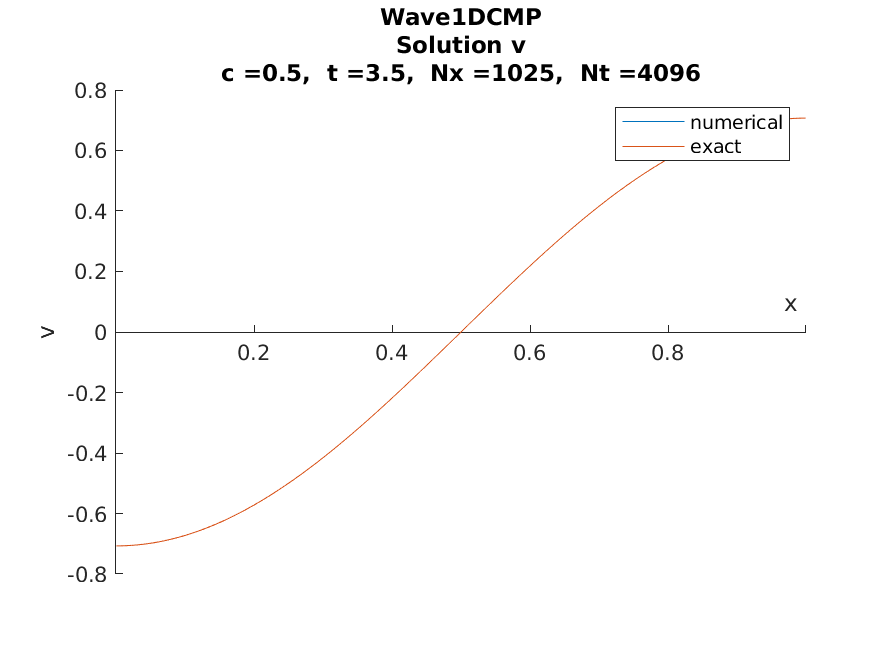} \\
C & D \\
\includegraphics[width = 3.0in, height = 2.5in, trim = 0mm 0mm 0mm 0mm, clip]{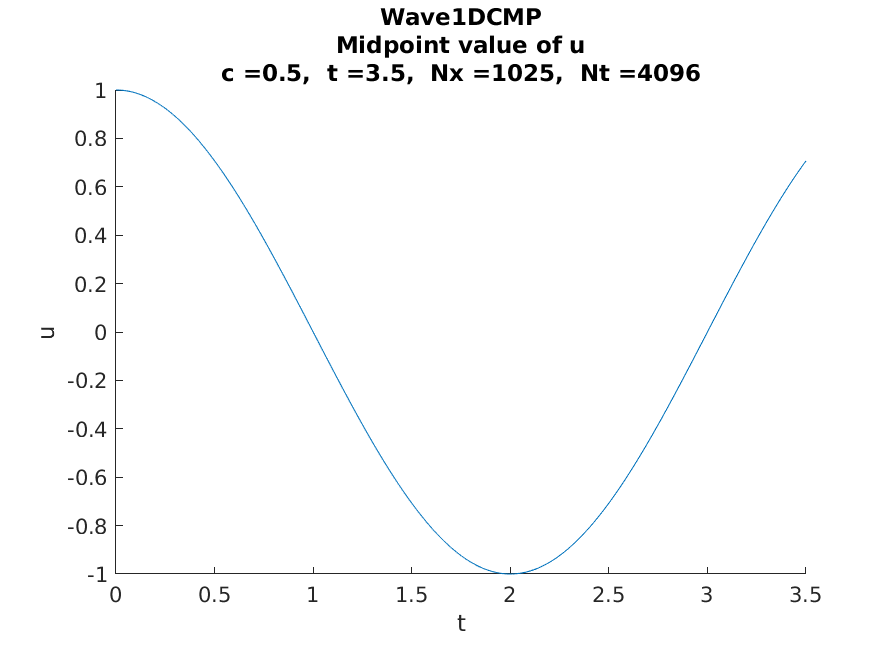} &
\includegraphics[width = 3.0in, height = 2.5in, trim = 0mm 0mm 0mm 0mm, clip]{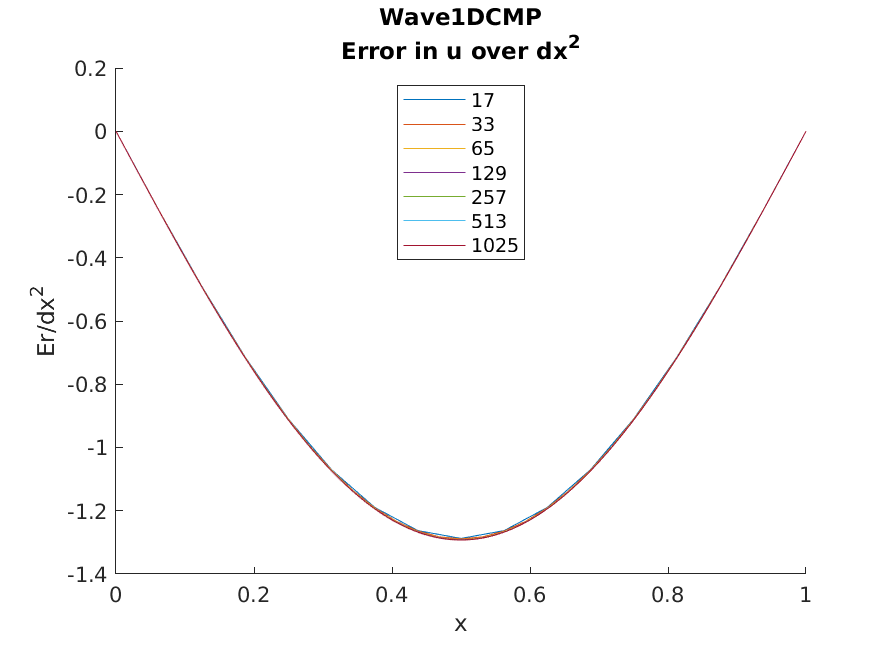} \\
\end{tabular}
\caption{Simulations using {\tt Wave1DCMP.m}.}
\label{Wave1DCM Solution}
\end{center}
\end{figure}

Figure \ref{Wave1DCM Solution} displays the results of seven simulations
run by {\tt Wave1DCMP.m}. The number of spatial grid points for each
simulation is given in figure \ref{Wave1DCM Solution} part D.  The initial
data were $u(x,0) = sin(\pi \, x)$ and $v(x,0) = 0$.
Parts A and B of figure \ref{Wave1DCM Solution} show the numerical and exact
solutions $u$ and $v$ at the final time for for the run using the most grid
points.  These solution are identical up to graphics resolution.  Part C of
figure \ref{Wave1DCM Solution} shows the value of the solution at the midpoint
of the interval for the run using the largest number of grid points.  So the
time interval used was slightly smaller than the period of the solution.
Part D of figure \ref{Wave1DCM Solution} plots the error in the solution
divided by $\dx^2$ for each of the seven simulations. Again all of these plots
are identical up to graphics resolution confirming the second order
convergence of the solution. 

\begin{figure}
\begin{center}
\begin{tabular}{cc}
$\rho(x) = 1+x/2 \,, \tau(x) = 1$ & $\rho(x) = 1 \,, \tau(x) = 1+x/2$ \\
\includegraphics[width = 3.0in, height = 2.5in, trim = 0mm 0mm 0mm 0mm, clip]{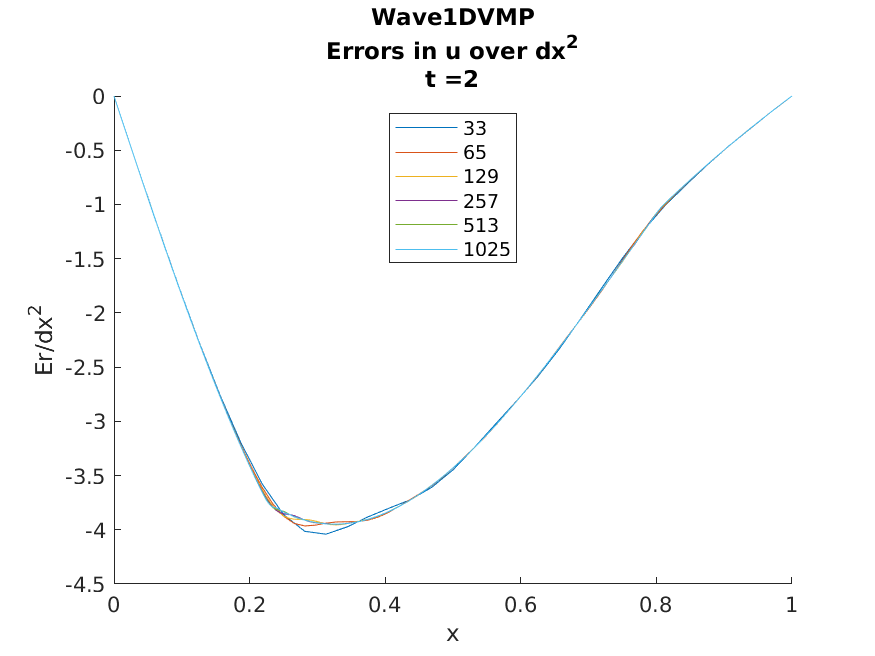} &
\includegraphics[width = 3.0in, height = 2.5in, trim = 0mm 0mm 0mm 0mm, clip]{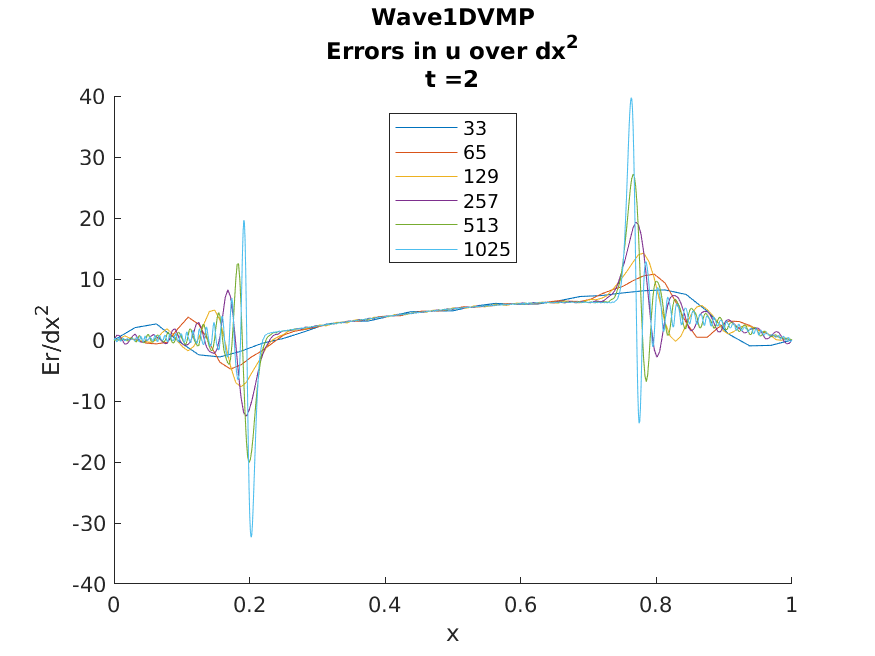} \\
$\rho(x) = 1-x/2 \,, \tau(x) = 1$ & $\rho(x) = 1 \,, \tau(x) = 1-x/2$ \\
\includegraphics[width = 3.0in, height = 2.5in, trim = 0mm 0mm 0mm 0mm, clip]{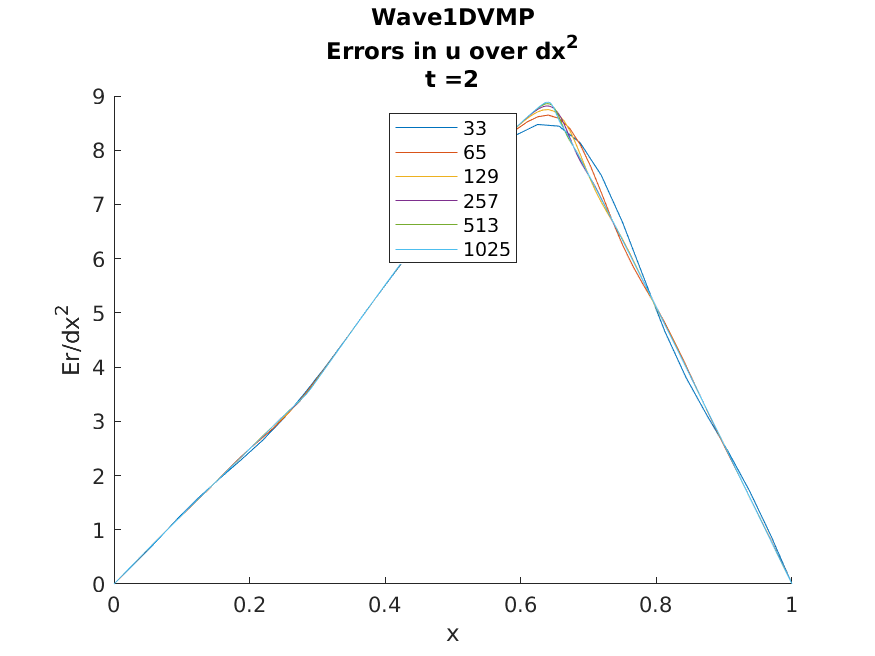} &
\includegraphics[width = 3.0in, height = 2.5in, trim = 0mm 0mm 0mm 0mm, clip]{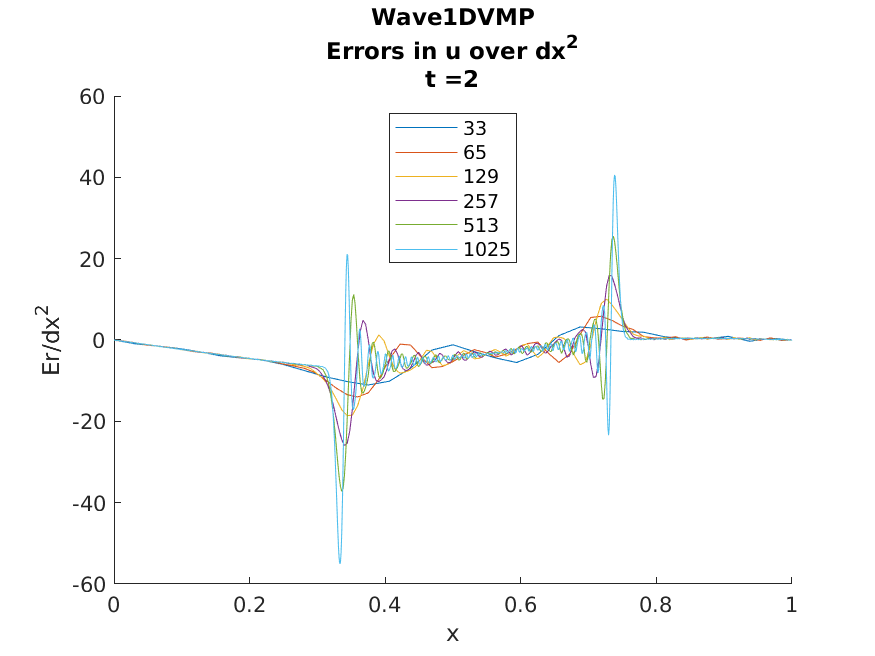} \\
\end{tabular}
\caption{Simulations using {\tt Wave1DVMP.m} with linear material properties.}
\label{Wave1DVM Linear Materials}
\end{center}
\end{figure}

\begin{figure}
\begin{center}
\begin{tabular}{cc}
$\rho(x) = 1 \,, \tau(x) = 1-x/2$ & $\rho(x) = 1 \,, \tau(x) = 1-x/2$ \\
\includegraphics[width = 3.0in, height = 2.5in, trim = 0mm 0mm 0mm 0mm, clip]{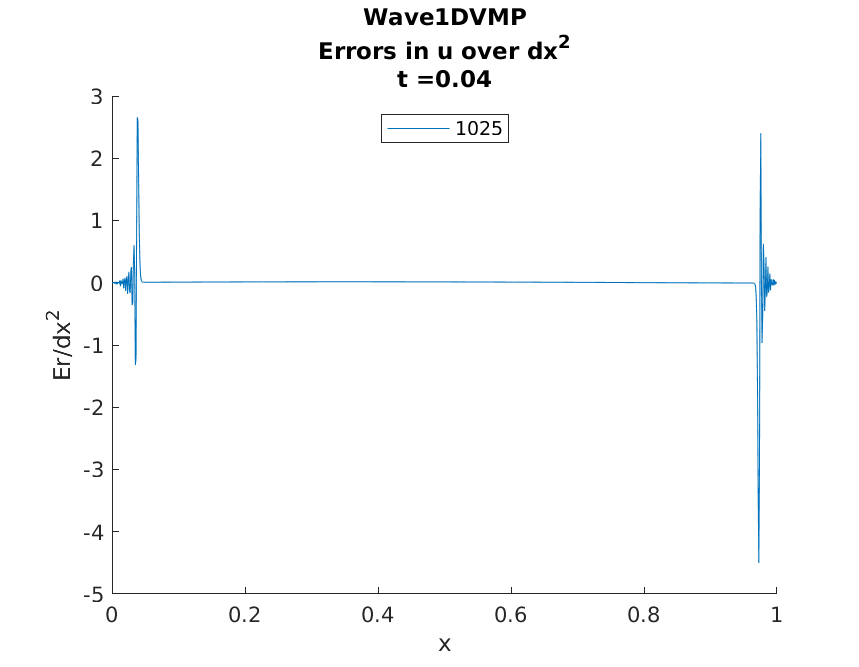} &
\includegraphics[width = 3.0in, height = 2.5in, trim = 0mm 0mm 0mm 0mm, clip]{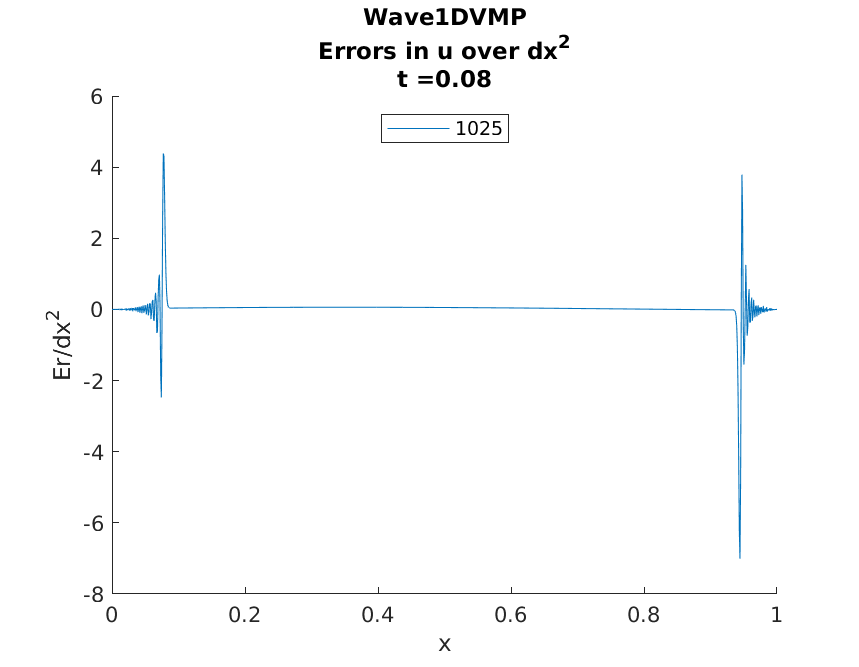} \\
$\rho(x) = 1 \,, \tau(x) = 1-x/2$ & $\rho(x) = 1 \,, \tau(x) = 1-x/2$ \\
\includegraphics[width = 3.0in, height = 2.5in, trim = 0mm 0mm 0mm 0mm, clip]{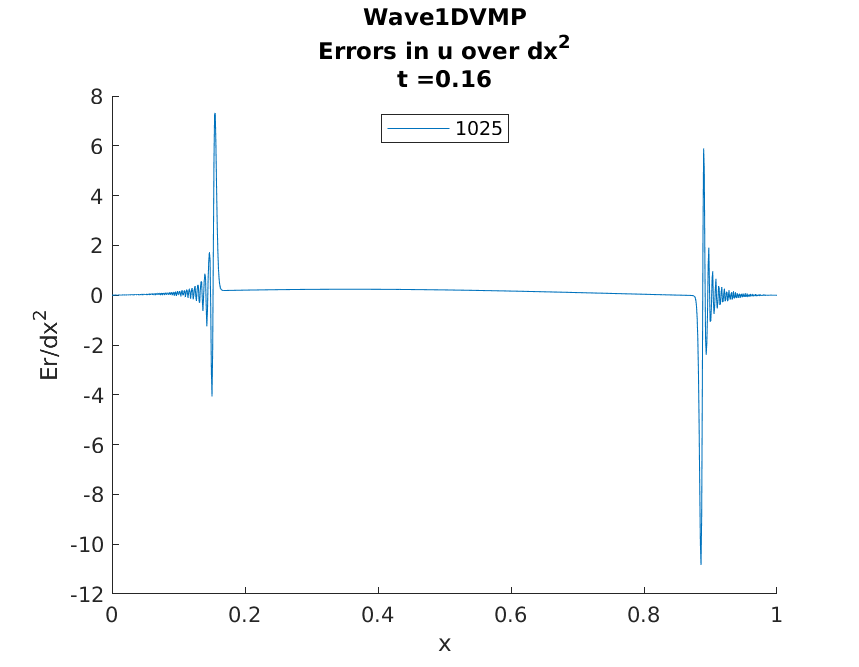} &
\includegraphics[width = 3.0in, height = 2.5in, trim = 0mm 0mm 0mm 0mm, clip]{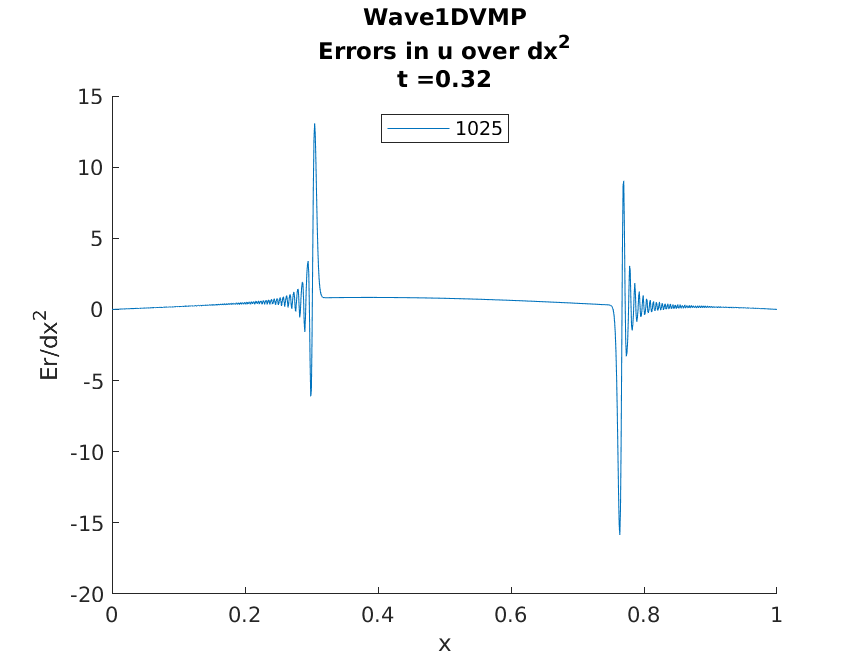} \\
\end{tabular}
\caption{Errors Near The Boundary for {\tt Wave1DVMP.m}. }
\label{Wave1DCM Boundary Errors}
\end{center}
\end{figure}

\begin{figure}
\begin{center}
\begin{tabular}{cc}
$ p = 1 ,\, q = 1 $ & $ p = 1 ,\, q = 2 $ \\
\includegraphics[width = 3.0in, height = 2.5in, trim = 0mm 0mm 0mm 0mm, clip]{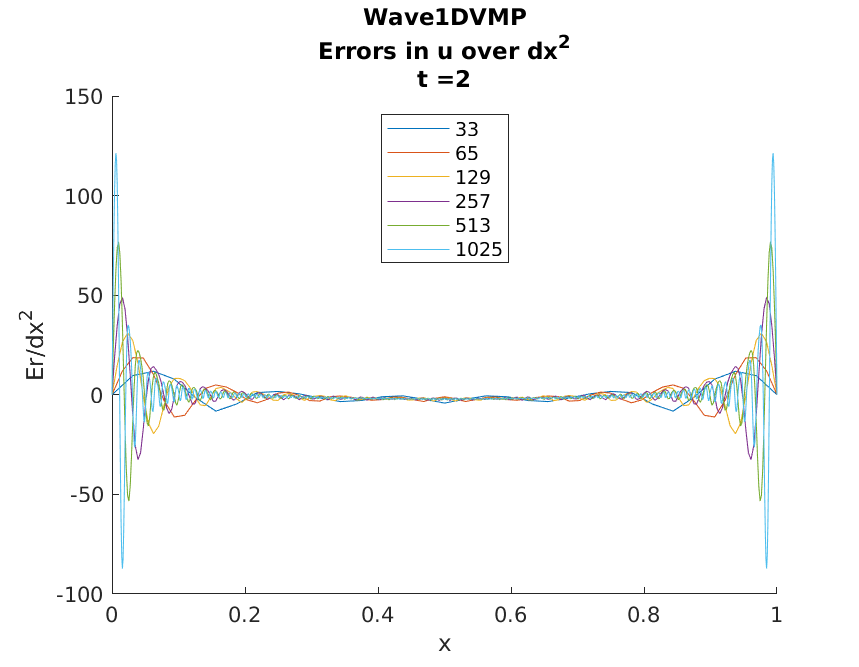} &
\includegraphics[width = 3.0in, height = 2.5in, trim = 0mm 0mm 0mm 0mm, clip]{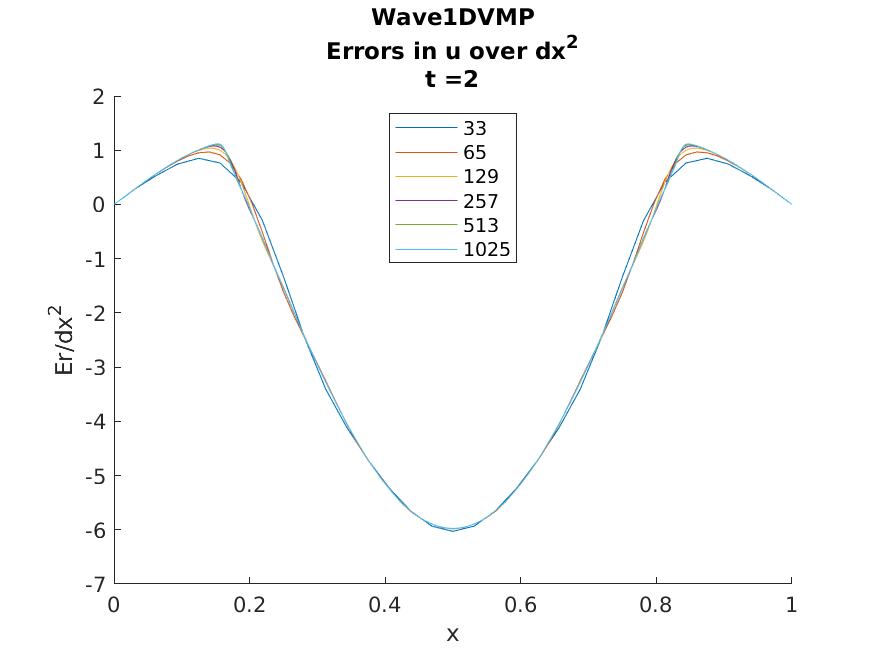} \\
$ p = 2 ,\, q = 1 $ & $ p = 2 ,\, q = 2 $ \\
\includegraphics[width = 3.0in, height = 2.5in, trim = 0mm 0mm 0mm 0mm, clip]{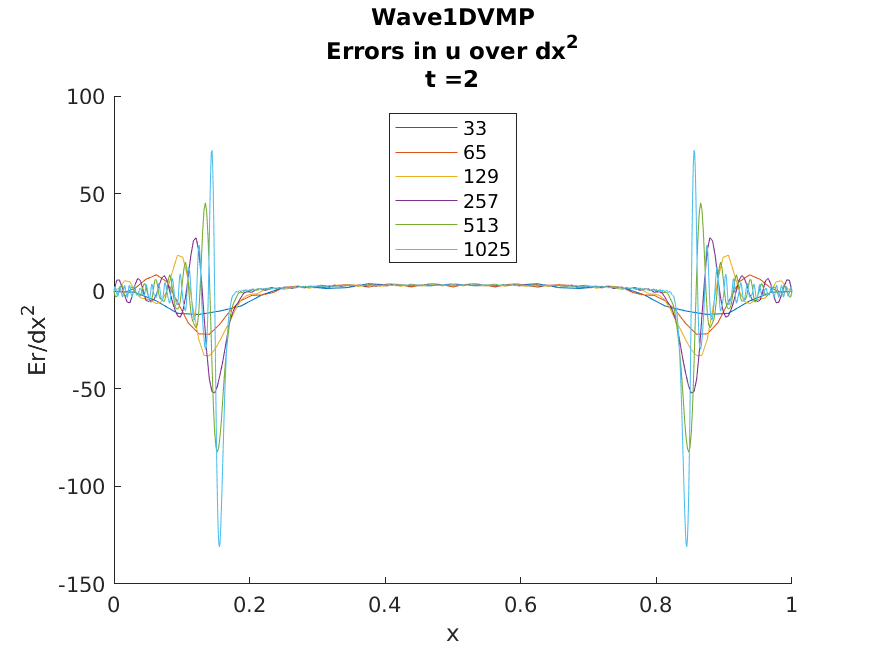} &
\includegraphics[width = 3.0in, height = 2.5in, trim = 0mm 0mm 0mm 0mm, clip]{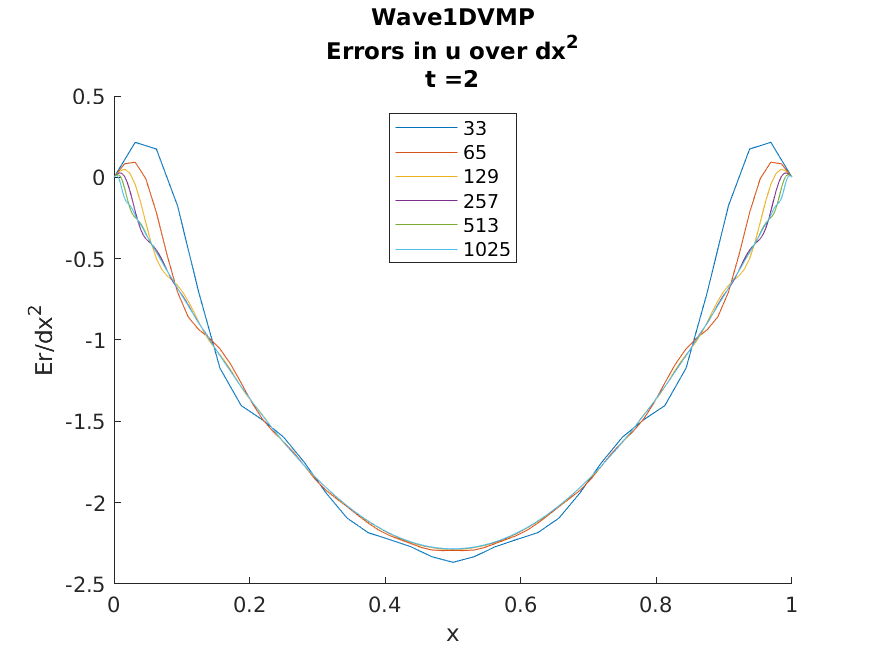} \\
\end{tabular}
\caption{Boundary Errors For The Bump in {\tt Wave1DVMP.m}. }
\label{Wave1DCM Bump Errors}
\end{center}
\end{figure}

The situation for {\tt Wave1DVMP.m} is a bit more complex.  The first test was
to define the material properties functions in terms of the wave speed, that
is, $\rho(x) = 1/c$ and $\tau(x) = c$ and then show that {\tt Wave1DVMP.m}
gave the same results as {\tt Wave1DCMP.m}. The results from the two codes are
the same.
In the case that at least one of $\rho$ and $\tau$ are not constant the
solutions can have small oscillatory errors that cause the convergence rate
to be less that 2 but still larger than 1.  To illustrate this, the errors
$Er(x)$ divided by $\Delta x^2$ are plotted for several examples in figures
\ref{Wave1DVM Linear Materials},
\ref{Wave1DCM Boundary Errors},
\ref{Wave1DCM Bump Errors},
\ref{PWLerrors} and
\ref{Wave1DCM Jump Errors} \,.
If the numerical solutions are converging at order 2 then the plots of
$Er(x)/\Delta x^2$ will be essentially independent of $\Delta x$.
When the order of convergence is less than 2 these plot will show increasing
$Er(x)/\Delta x^2$.  The legend gives the number of points in the grids and
then $\Delta x = 1/(Nx-1)$.

A simple case to start with is to choose one of $\rho$ and $\tau$ to be a
linear function.  As shown in Figure \ref{Wave1DVM Linear Materials} there
are some surprises.  First the legends give the number of spatial points
$N_x$ so for $Nx = 1025$ the errors are order $10^{-6}$. As the length
of the interval is 1, $\Delta x^2 \approx 10^{-6}$ so the errors are small.
Note that constant $\rho$ simulations are second order accurate so it is
the variable $\tau$ that causes the small oscillatory error.
Next Figure \ref{Wave1DCM Boundary Errors} shows that for one choice of
$\rho$ and $\tau$ the larger errors start at the boundary and move to the
interior. The final times double from frame to frame so roughly the oscillatory
errors move away from the boundary linearly and grow linearly. 

In Figure \ref{Wave1DCM Bump Errors} 
\begin{align*}
\rho(x) = 1+(2 \, x \, (1-x))^p \,,\quad 
\tau(x) = 1+(2 \, x \, (1-x))^q \,.
\end{align*}
The panels contain plots of the solution errors fore several values of
$p$ and $q$.  The number of grid points are given in the legend.
The final time $t =2$ is close to a ''period'' of the solution.
Note that for $p > 1$ and $q > 1$ the derivatives of $\rho$ and $\tau$
are zero at the boundary while if $p=1$ and $q=1$ the derivatives are not
zero. For $q = 2$ and $p = 1$ or $p = 2$ there are no oscillatory errors
while if $q =1$ and $p = 1$ or $p = 2$ there are oscillatory errors.
It seems that the derivative of $\tau$ not being zero at boundary
causes the oscillatory errors.

\begin{figure}
\begin{center}
   \includegraphics[width=3.00in,trim = 0 0 0 0, clip ]{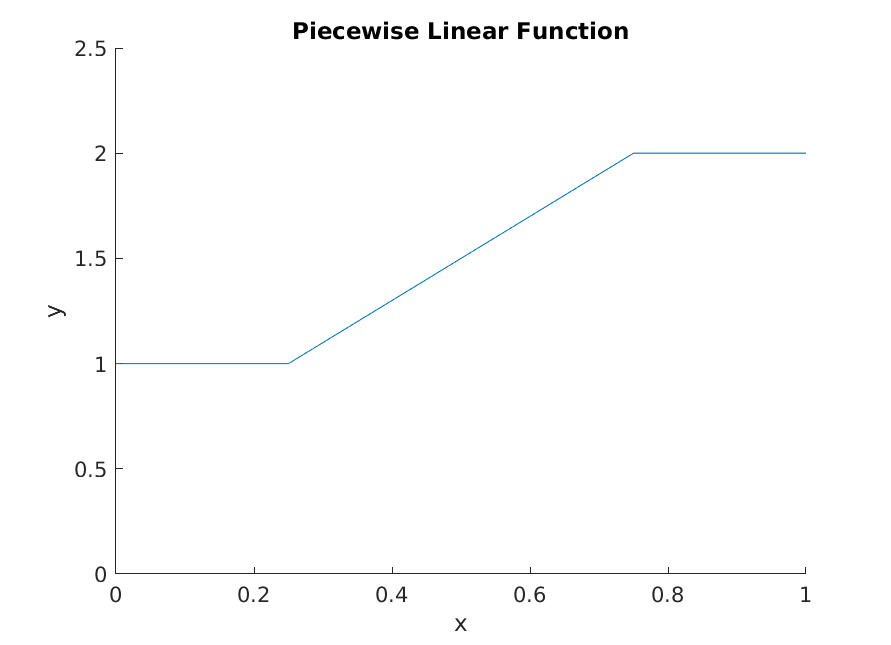}
\caption{Piecewise Linear Material Property}
\label{PWLM}
\end{center}
\end{figure}

To test the hypothesis that nonzero derivatives at the boundary cause the
oscillatory errors in the next examples either $\rho$ or $\tau$ is a piecewise
linear functions as shown in Figure \ref{PWLM}. If $H$ is the Heaviside
function, that is, $H(x) = 1$ if $x \geq 0$ and $H(x) = 0$ for $x<0$ then
$\rho$ or $\tau$ are given by
\begin{equation}
   c \, (1-H(x-a))
 + d \, H(x-b)
 + \frac{(a \, d-b \, c) + (c-d) \, x}{a-b} \, (H(x-a)-H(x-b))\,.
\end{equation}
In Figure \ref{PWLM} and in the simulations $ a = 1/4$, $ b = 3/4$, $ c = 1$,
$ d = 2$.

\newpage

Figure \ref{PWLerrors} shows that the solutions errors.
Note that $y'(0) = y'(1) = 0$. So the conjecture above that the derivative
at the boundary being zero eliminates the oscillations is not correct.
Additionally the discontinuities in the derivative $d \tau / d x$ do cause
similar oscillatory errors. Again the errors are a small multiple of
$\Delta x^2$ that grows slowly with decreasing $\Delta x$.

\begin{figure}
\begin{center}
\begin{tabular}{cc}
Piecewise linear $\rho$ & Constant $\rho$ \\
Constant $\tau$ & Piecewise linear $\tau$ \\
\includegraphics[width = 3.0in, height = 2.5in, trim = 0mm 0mm 0mm 0mm, clip]{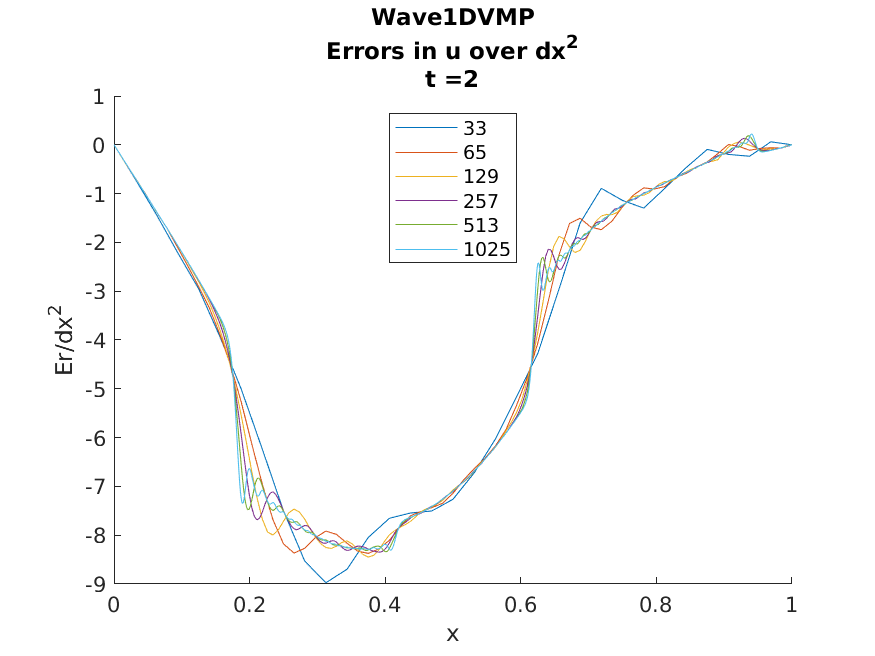} &
\includegraphics[width = 3.0in, height = 2.5in, trim = 0mm 0mm 0mm 0mm, clip]{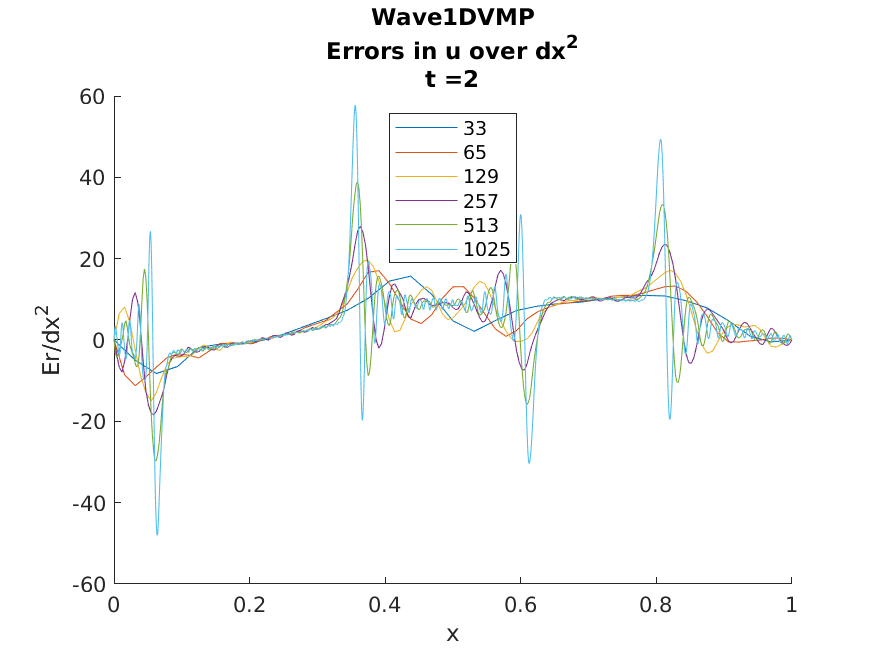} \\
\end{tabular}
\caption{Errors For Piecewise Linear $\rho$ or $\tau$ in {\tt Wave1DVMP.m}}.
\label{PWLerrors}
\end{center}
\end{figure}

The final test is to see what discontinuities in $\rho$ and $\tau$ do to the
errors, that is one of $\rho$ or $\tau$ jumps up or down:
\[
\rho(x) = 1 \pm H(x)/2 \,;\quad \tau(x) = 1 \pm H(x)/2 \,.
\]
Figure \ref{Wave1DCM Jump Errors} show that the jump causes oscillations
similar to those seen in the previous figures near the boundary of the $x$
interval.

\begin{figure}
\begin{center}
\begin{tabular}{cc}
$\rho$ jumps up & $\rho$ jumps down \\
\includegraphics[width = 3.0in, height = 2.5in, trim = 0mm 0mm 0mm 0mm, clip]{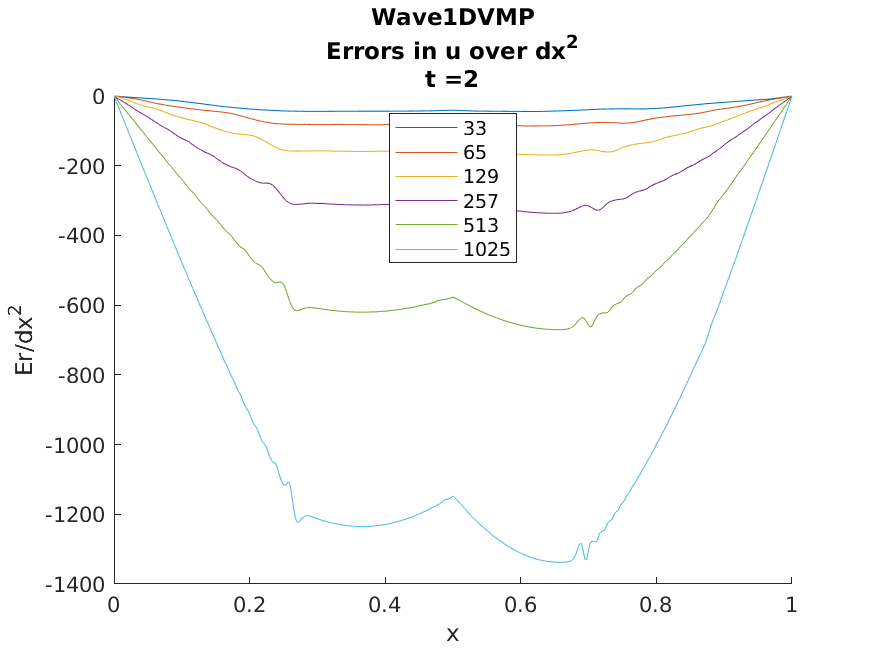} &
\includegraphics[width = 3.0in, height = 2.5in, trim = 0mm 0mm 0mm 0mm, clip]{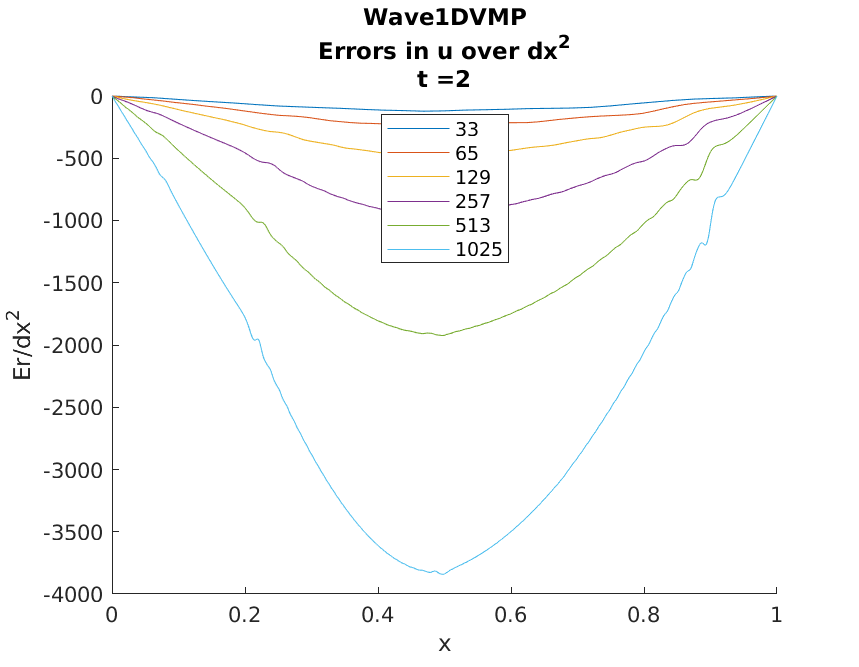} \\
$\tau$ jumps up & $\tau$ jumps down \\
\includegraphics[width = 3.0in, height = 2.5in, trim = 0mm 0mm 0mm 0mm, clip]{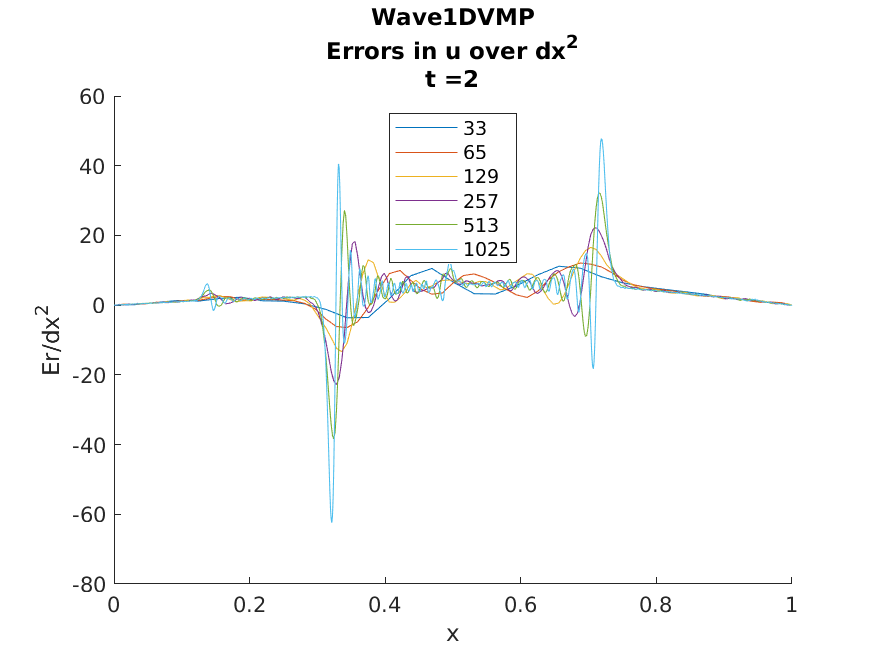} &
\includegraphics[width = 3.0in, height = 2.5in, trim = 0mm 0mm 0mm 0mm, clip]{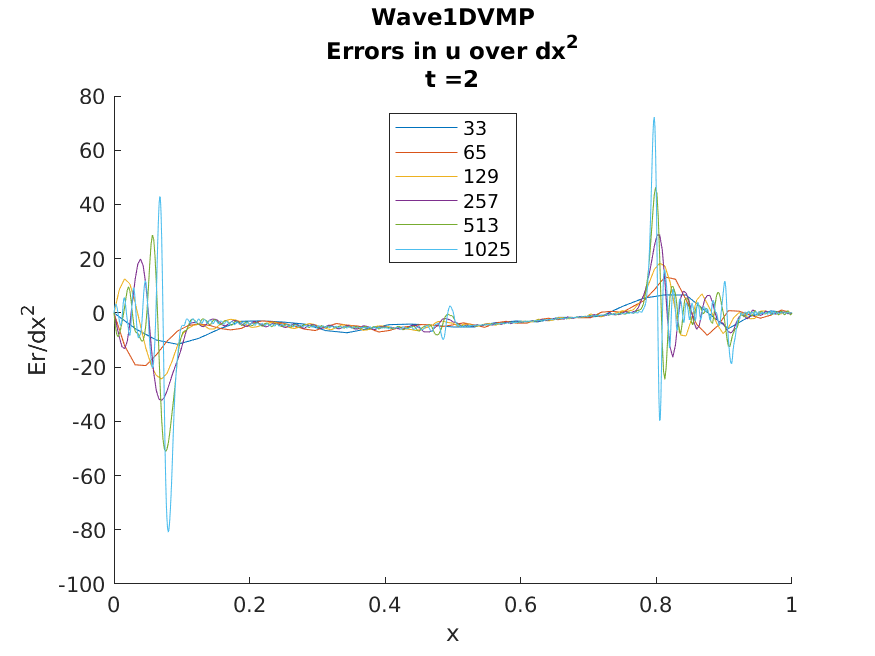} \\
\end{tabular}
\caption{$\rho$ or $\tau$ have a jump at $x=1/2$ in {\tt Wave1DVMP.m}. }
\label{Wave1DCM Jump Errors}
\end{center}
\end{figure}

\newpage \clearpage
\setcounter{equation}{0}
\section{Three Dimensional Variable Coefficient Differential Operators
\label{Second Order D0s}}
This section introduces terminology based on that used in differential geometry
of exact sequences and diagram chasing that are used to create second order
variable coefficient differential operators that are be used to create 
second order wave equations. However, using the ideas from differential
geometry in the mimetic finite difference \cite{RobidouxSteinberg2011} setting
poses a new problem as there are two interlaced discretization grids.
Consequently two exact sequences as shown in Figure \ref{Exact-Sequences} are
need to study the discretizations. In the continuum the two exact sequences
are not needed but are useful when the continuum is discretized.  The use of
diagram chasing is critical to obtain differential operators that can be
discretized using the mimetic finite difference methods as described in
\cite{RobidouxSteinberg2011}.

After the operators are defined, inner products are introduced and the adjoints
of all of the operators obtained from the double exact sequence are computed.
This is then used to show that the second order operators are self-adjoint and
either positive or negative and thus suitable for defining wave equations.
In this continuum setting there are two equivalent derivations of each
operator. In the discrete setting these operators will be different.

\subsection{Exact Sequences} \label{Exact Sequences}
\begin{figure}
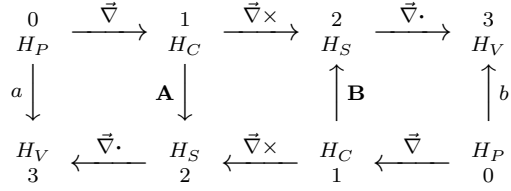

\begin{equation*}
  \begin{CD}
		   {\mystack{0}{H_P}}
        @>\grad >> {\mystack{1}{H_C}}
        @>\curl >> {\mystack{2}{H_S}}
        @>\divg >> {\mystack{3}{H_V}}
\\
@V{a}VV @V{{\bf A}}VV @AA{{\bf B}}A @AA{b}A @. \\
                   {\mystack{H_V}{3}}
        @<\divg << {\mystack{H_S}{2}}
        @<\curl << {\mystack{H_C}{1}}
        @<\grad << {\mystack{H_P}{0}}
\\
  \end{CD}
\end{equation*}
\caption{Continuum Double Exact Sequence Diagram}
\label{Exact-Sequences}
\end{figure}

In the double exact sequence diagram \ref{Exact-Sequences} the bottom row is
the same as the top row written in opposite order. For reasons which will
become clear when the operators are discretized,
$P$ stands for points,
$C$ stands for curves,
$S$ stands for surfaces,
$V$ stands for volumes.
In this diagram $H_P$ and $H_V$ are linear spaces of smooth scalar functions
depending on the spatial variables $(x,y,z)$ and $H_C$ and $H_S$ are linear
spaces of smooth vector valued functions that also depend on $(x,y,z)$. In this
discussion all of the functions converge rapidly to zero as $x^2+y^2+z^2$
becomes large.

\begin{table}[ht]
\begin{center}
\begin{tabular}{|l|l|}
\hline
$\text{if }f \in H_P$        $\text{ then } \grad f \in H_C$ &
$\text{if }f \in H_P$        $\text{ then } a \, f \in H_V$ \\
$\text{if }\vec{v} \in H_C$  $\text{ then } \curl \vec{v} \in H_S$ &
$\text{if }\vec{v} \in H_C$  $\text{ then } {\bf A} \, \vec{v} \in H_S$ \\
$\text{if }\vec{w} \in H_S$  $\text{ then } \divg \vec{w} \in H_V$ &
$\text{if }\vec{w} \in H_C$  $\text{ then } {\bf B} \, \vec{w} \in H_S$ \\
& $\text{if }g \in H_P$        $\text{ then } b \, g \in H_V$  \\
\hline
\end{tabular}
\caption{First order operators on the left, material property operators on
the right.}
\label{Basic Operators}
\end{center}
\end{table}

The first order differential operators in the exact sequence
\ref{Exact-Sequences} are the gradient $\grad$, curl or rotation $\curl$, and
divergence $\divg$.  The diagram is called {\em exact} because the curl of the
gradient is zero $\curl \, \grad = 0$ and the divergence of the curl is zero
$\divg \, \curl = 0$.  The scalar functions $a$ and $b$ and also the matrix
valued functions ${\bf A}$ and ${\bf B}$ are used to describe material
properties and act by multiplication as described in Table
\ref{Basic Operators}.
The functions $a$ and $b$ are bounded above and below by positive constants.
The matrix functions are symmetric positive definite and the eigenvalues of
the matrices are bounded above and below by positive constants.  To begin,
these functions will be assumed smooth but later this will not be assumed so
that they can be used to describe discontinuous material properties.

So in Figure \ref{Basic Operators} the horizontal arrows represent the action
of the differential operators while the vertical arrows represent multiplication
by scalar functions $a$ and $b$ and by $3 \times 3$ matrices ${\bf A}$ and
${\bf B}$ that are known as {\em star} operators in differential geometry.
The differential operators and material property functions give mappings
between the spaces in the double exact sequence as described in Table
\ref{Basic Operators}.  Note that the differential operators are not
invertible, but that the conditions on $a$, $b$, ${\bf A}$ and ${\bf B}$ imply
that they are invertible. This will play a critical role in diagram chasing.

The integers in the double exact sequence give the spatial dimension of the
function in the spaces, that is if $f \in H_P$ the $f$ has no spatial
dimension while if $g \in H_V$ then $g$ has spatial dimension $1/d^3$.
Also if $\vec{v} \in H_C$ then $\vec{v}$ has spatial dimension $1/d$
and if $\vec{w} \in H_S$ then $\vec{w}$ has spatial dimension $1/d^2$.
The differential operators all have dimension $1/d$. Moreover $a$ and $b$
have spatial dimensions $1/d^3$ while $\bf A$ and $\bf B$ have dimensions
$1/d$. The directions of the vertical arrows were chosen so that $a$, $b$,
$\bf A$ and $\bf B$ have dimensions $1/d^k$ for $k>0$.
Importantly, the operators and spaces in Table \ref{Basic Operators} are
dimensionally consistent! In some applications $a$ and $b$ are densities.

\subsection{Diagram Chasing and Second Order Operators}

\begin{table}[ht]
\begin{center}
\begin{tabular}{|r|r||r|r|}
\hline
\multicolumn{2}{|c|}{Upper Row} &
\multicolumn{2}{|c|}{Bottom Row} \\
\hline
\multicolumn{4}{|c|}{First Box} \\
\hline
  $f \in H_P$                       
& $\vec{v} \in H_C$
& $\vec{w} \in H_S$
& $g \in H_V$ \\
  $\grad f \in H_C$                 
& ${\bf A} \vec{v} \in H_S$
& $\divg \vec{w} \in H_V$
& $a^{-1} g \in H_P$ \\
  ${\bf A}  \grad f \in H_S$              
& $\divg {\bf A} \vec{v} \in H_V$
& $a^{-1} \divg \vec{w} \in H_P$
& $\grad a^{-1} g \in H_C$ \\
  $\divg {\bf A}  \grad f \in H_V$        
& $a^{-1} \divg {\bf A} \vec{v} \in H_P$
& $\grad a^{-1} \divg \vec{w} \in H_C$
& ${\bf A} \grad a^{-1} g \in H_S$ \\
  $a^{-1} \divg {\bf A}  \grad f \in H_P$ 
& $\grad a^{-1} \divg {\bf A} \vec{v} \in H_C$
& ${\bf A} \grad a^{-1} \divg \vec{w} \in H_S$
& $\divg {\bf A} \grad a^{-1} g \in H_V$ \\
\hline 
\multicolumn{4}{|c|}{Second Box} \\
\hline
  $\vec{v} \in H_C$                       
& $\vec{w} \in H_S$
& $\vec{v} \in H_C$
& $\vec{w} \in H_S$ \\
  $\curl \vec{v} \in H_S$                       
& ${\bf B}^{-1} \vec{w} \in H_C$
& $\curl \vec{v} \in H_S$
& ${\bf A}^{-1} \vec{w} \in H_C$ \\
  ${\bf B}^{-1} \curl \vec{v} \in H_C$                       
& $\curl {\bf B}^{-1} \vec{w} \in H_S$
& ${\bf A}^{-1} \curl \vec{v} \in H_C$
& $\curl {\bf A}^{-1} \vec{w} \in H_S$ \\
  $\curl {\bf B}^{-1} \curl \vec{v} \in H_S$                       
& ${\bf A}^{-1} \curl {\bf B}^{-1} \vec{w} \in H_C$
& $\curl {\bf A}^{-1} \curl \vec{v} \in H_S$
& ${\bf B}^{-1} \curl {\bf A}^{-1} \vec{w} \in H_C$ \\
  ${\bf A}^{-1} \curl {\bf B}^{-1} \curl \vec{v} \in H_C$
& $\curl {\bf A}^{-1} \curl {\bf B}^{-1} \vec{w} \in H_S$
& ${\bf B}^{-1} \curl {\bf A}^{-1} \curl \vec{v} \in H_C$
& $\curl {\bf B}^{-1} \curl {\bf A}^{-1} \vec{w} \in H_S$ \\
\hline 
\multicolumn{4}{|c|}{Third Box} \\
\hline
  $\vec{w} \in H_S$
& $g \in H_V$
& $f \in H_P$
& $\vec{v} \in H_C$ \\
  $\divg \vec{w} \in H_V$
& $b^{-1} g \in H_P$
& $\grad f \in H_C$
& ${\bf B}^{-1} \vec{v} \in H_S$ \\
  $b^{-1} \divg \vec{w} \in H_P$
& $\grad {b^{-1}} g \in H_C$ 
& ${\bf B} \grad f \in H_S$
& $\divg {\bf B}^{-1} \vec{v} \in H_V$ \\
  $\grad b^{-1} \divg \vec{w} \in H_C$
& ${\bf B} \grad b^{-1} g \in H_S$
& $\divg {\bf B} \grad f \in H_V$
& $b^{-1} \divg {\bf B}^{-1} \vec{v} \in H_P$ \\
  ${\bf B} \grad b^{-1} \divg \vec{w} \in H_S$
& $\divg {\bf B} \grad b^{-1} g \in H_V$
& $b^{-1} \divg {\bf B} \grad f \in H_P$
& $\grad b^{-1} \divg {\bf B}^{-1} \vec{v} \in H_C$ \\
\hline 
\end{tabular}
\caption{Fundamental second order differential operators. The
left two columns start with spaces in the top row in Figure
\ref{Exact-Sequences}
while the right two columns start with spaces in the bottom row.}
\label{General Fundamental Operators}
\end{center}
\end{table}

\begin{table}[ht]
\begin{center}
\begin{tabular}{|r|r|r|r|}
\hline
\multicolumn{4}{|c|}{First Box} \\
\hline
  $f \in H_P$                       
& $\vec{v} \in H_C$
& $\vec{w} \in H_S$
& $g \in H_V$ \\
  $c^2 \divg  \grad f \in H_P$ 
& $c^2 \grad \divg \vec{v} \in H_C$
& $c^2 \grad \divg \vec{w} \in H_S$
& $c^2 \divg \grad  g \in H_V$ \\
\hline 
\multicolumn{4}{|c|}{Second Box} \\
\hline
  $c^2 \vec{v} \in H_C$                       
& $c^2 \vec{w} \in H_S$
& $c^2 \vec{v} \in H_C$
& $c^2 \vec{w} \in H_S$ \\
  $c^2 \curl \curl \vec{v} \in H_C$                       
& $c^2 \curl \curl \vec{w} \in H_S$
& $c^2 \curl \curl \vec{v} \in H_C$
& $c^2 \curl \curl \vec{w} \in H_S$ \\
\hline 
\multicolumn{4}{|c|}{Third Box} \\
\hline
  $\vec{w} \in H_S$
& $g \in H_V$
& $f \in H_P$
& $\vec{v} \in H_C$ \\
  $c^2 \grad \divg \vec{w} \in H_S$
& $c^2 \divg \grad g \in H_V$
& $c^2 \divg \grad f \in H_P$
& $c^2 \grad \divg \vec{v} \in H_C$ \\
\hline 
\end{tabular}
\caption{The general operators can be simplified using
the assumptions in \ref{Simplifying Assumptions}.
\label{Simplified General Operators}}
\end{center}
\end{table}

Table \ref{General Fundamental Operators} gives all of the possible second
order operators given by diagram chasing. As an example of diagram chasing,
consider
$f \in H_P$ so that
$\grad f \in H_C$ and then
${\bf A} \grad f \in H_S$ so that
$\divg {\bf A} \grad f \in H_V $ and finally
$a^{-1} \divg {\bf A} \grad f \in H_S $.
Consequently $a^{-1} \divg {\bf A} \grad$ maps $H_S$ into $H_S$.
The gives the upper left entry in Table \ref{General Fundamental Operators}.
The remaining operators are created similarly.

For diagram chasing it is important that the mappings $a$, $b$, ${\bf A}$ and
${\bf B}$ are invertible while $\grad$, $\curl$ and $\divg$ are not invertible.
Consequently to create a second order operator, only going clockwise around a
square in Figure \ref{Exact-Sequences} is allowed. However it is possible to
start and any corner, so this gives twelve operators, four corners times three
squares.  In the continuum some of these operators are essentially the same,
for example $a^{-1} \divg {\bf A}  \grad$ and $b^{-1} \divg B \grad$. What is
important is that when these operators are discretized they will not be the
same.  The assumptions that ${\bf B}= {\bf A}$ and $b = a$ reduces the number
of operators to six. When the material properties are constant it is useful
to introduce that wave speed and then then use the simplifying assumptions
\begin{equation} \label{Simplifying Assumptions}
a = 1/c^3, b = 1/c^3, {\bf A} = {\bf I}/c, {\bf B} = {\bf I}/c,
\end{equation}
where $\bf I$ is the identity matrix then the operators simplify to those
in Table \ref{Simplified General Operators}.
So under the simplifying assumptions there are only three distinct second order operators:
\[
\lap f = \divg \grad f \,; \quad \quad
\curl \curl \vec{v} \,; \quad \quad
\grad \divg  \vec{w} \,.
\]
where $\lap$ is the scalar Laplacian.

\subsection{Additional Second Order Operators}
Note that in Table \ref{General Fundamental Operators}
there are two operators defined on $H_P$, four operators defined on $H_C$,
four operators defined on $H_S$ and two operators defined on $H_V$. If
linear operators are defined on the same space then linear combinations
of these operators are again linear operators. The two operators defined
on $H_P$ and the two defined on $H_V$ are essentially the same so linear
combinations are not interesting. For any operator in the {\em Top Row}
boxes, there is a similar operator in the {\em Bottom Row} that can
be obtained by interchanging $a$ with $b$ and ${\bf A}$ with ${\bf B}$.
For $\vec{v} \in H_C$ and for $\vec{w} \in H_S$ define
\begin{align}
{\bf VL}_1 (\vec{v}) & = 
\grad a^{-1} \, \divg \, {\bf A} \, \vec{v} -
{\bf A}^{-1} \, \curl \, {\bf B}^{-1} \, \curl \, \vec{v} \, \nonumber \\
{\bf VL}_2 (\vec{w}) & = 
{\bf B} \, \grad b^{-1} \, \divg \, \vec{w} -
\curl \, {\bf A}^{-1} \, \curl \, {\bf B}^{-1} \, \vec{w} \, \nonumber \\
{\bf VL}_3 (\vec{v}) & =
\grad \, b^{-1} \, \divg \, {\bf B} \, \vec{v} - 
{\bf B}^{-1} \, \curl \, {\bf A}^{-1} \, \curl \, \vec{v} \label{Linear Combination Operators} \\
{\bf VL}_4 (\vec{w}) & =
{\bf A} \, \grad \, a^{-1} \, \divg \, \vec{w}-
\curl \, {\bf B}^{-1} \, \curl \, {\bf A}^{-1} \, \vec{w} \nonumber \\
\nonumber
\end{align}

Under the simplifying assumptions \ref{Simplifying Assumptions} and with
$a = b = 1$ and ${\bf A} = {\bf B} = {\bf I}$ these operators become
\[
\lapv \vec{v} = \grad \divg \vec{v} - \curl \curl \vec{v} \,,\quad
\lapv \vec{w} = \grad \divg \vec{w} - \curl \curl \vec{w} \,,
\]
where $\lapv$ is the vector Laplacian
\[
\lapv \vec{v} = \left(\lap v_1,\lap v_2,\lap v_3\right) \,,\quad
\lapv \vec{w} = \left(\lap w_1,\lap w_2,\lap w_3\right) \,.
\]
Operators like these appear in the elastic wave equation which will
be studied in Section \eqref{Elastic Wave Equation}.

\subsection{Inner Products}
Applying the mimetic ideas to physical problems requires the use of inner
products on the spaces $H_P$, $H_C$, $H_S$ and $H_V$. It is important that
the inner products do not have a spatial dimension.  Two bilinear forms
will help simplify the notation. Let
$f = f(x,y,z) \in H_P$,\,
$g = g(x,y,z) \in H_V$,\,
$\vec{v} = \vec{v(x,y,z)} \in H_C$
and
$\vec{v} = \vec{w(x,y,z)} \in H_S$ and then define 
\begin{align}
\lbilin f , g \rbilin & = 
\int_{\real^3} f(x,y,z) \, g(x,y,z) \, dx \, dy \, dz \,, \nonumber \\
\lbilin \vec{v} , \vec{w} \rbilin & = 
\int_{\real^3} \vec{v}(x,y,z) \bdot \vec{w}(x,y,z) \, dx \, dy \, dz \,.
\label{bilinear form}
\end{align}
These bilinear forms are dimensionless because $dx$, $dy$ and $dy$ have
dimension $d$ while $f$ has dimension $0$, $g$ has dimension $1/d^3$,
$\vec{v}$ has dimension $1/d$ and $\vec{w}$ has $1/d^2$.

The inner product on the function spaces must use a weight function to
be dimensionless: 
\begin{itemize}
\item
for $f_1 \,, f_2 \in H_P$ set
$\left< f_1 , f_2 \right>_P = \lbilin a \, f_1 , f_2 \rbilin$
\item
for $\vec{v}_1 \,, \vec{v}_2 \in H_C$ set
$ \left< \vec{v}_1 , \vec{v}_2 \right>_C = 
\lbilin {\bf A} \vec{v}_1 , \vec{v}_2 \rbilin$
\item
for $\vec{w}_1 \,, \vec{w}_2 \in H_S$ set
$\left< \vec{w}_1 , \vec{w}_2 \right>_S = 
\lbilin {\bf A}^{-1} \, \vec{w}_1 , \vec{w}_2 \rbilin$ 
\item
for $g_1 \,, g_2 \in H_V$ set
$\left< g_1 , g_2 \right>_V = 
\lbilin a^{-1} \, g_1 , g_2 \rbilin$
\end{itemize}
As usual
$\norm{f}_P^2 = \left< f , f \right>_P$,
$\norm{\vec{v}}_C^2 = \left< \vec{v} , \vec{v} \right>_C$,
$\norm{\vec{w}}_S^2 = \left< \vec{w} , \vec{w} \right>_S$,
$\norm{g}_V^2 = \left< g , g \right>_V$.
Additional inner products can be made by replacing $a$ by $b$ and
${\bf A}$ by ${\bf B}$. To be inner products it is important that
$a$ and $b$ are positive and that ${\bf A}$ and ${\bf B}$ are
symmetric and positive definite matrices.

\subsection{Adjoint Operators \label{Continuum Adjoint Operators}}
Adjoints are commonly defined for operators mapping a space into itself
but most of the operators used here are mapping between two different
spaces, so the adjoint is defined as in Section \ref{Continuous Time}.
The discussion in that section shows that if $X$, $Y$, and $Z$ are 
linear spaces and $O_1$ and $O_2$ are linear operator such that
\[
X \overset{O_1}{\rightarrow} Y \overset{O_2}{\rightarrow} Z \\
\]
then
\[
Z \overset{O_2^*}{\rightarrow} Y \overset{O_1^*}{\rightarrow} X \,
\]
and 
\begin{equation}
(O_1 O_2)^* = O_2^* \, O_1^*  \,.
\label{Product Adjoints}
\end{equation}
Because diagram chasing gives operators as compositions of other
operators, this will be used many times.

The adjoints of the operators in Table \ref{Basic Operators} are
\begin{align}
\grad^* & = - a^{-1} \, \divg \, {\bf A} &
\grad^* & = - b^{-1} \, \divg \, {\bf B}
\,,\nonumber\\
\curl^* & = + {\bf A}^{-1} \, \curl \, {\bf B}^{-1} &
\curl^* & = + {\bf B}^{-1} \, \curl \, {\bf A}^{-1}
\,,\nonumber\\
\divg^* & = - {\bf B} \grad  b^{-1} &
\divg^* & = - {\bf A} \grad  a^{-1} 
\,,\nonumber\\
{\bf A}^* & = {\bf A}^{-1} &
{\bf B}^* & = {\bf B}^{-1}
\,, \label{Adjoints} \\
a^* & = a^{-1} &
b^* & = b^{-1}
\nonumber\,,
\end{align}
where the column on the left contains differential operators from the top
row in Figure \ref{Exact-Sequences} and the column on the right contains
differential operators from the bottom row in Figure \ref{Exact-Sequences}.

The proofs of the formulas in \ref{Adjoints} are straight forward.
For the gradient let $f \in H_P$ and $\vec{v} \in H_C$ so that
\begin{align*}
\left< \grad f , \vec{v} \right>_C 
 & = \lbilin {\bf A} \, \grad f , \vec{v} \rbilin \\
 & = \lbilin \grad f , {\bf A} \, \vec{v} \rbilin \\
 & = - \lbilin f , \divg \, {\bf A} \vec{v} \rbilin \\
 & = - \lbilin a \, f , a^{-1} \, \divg \, {\bf A} \vec{v} \rbilin \\
 & =   \left< f , - a^{-1} \, \divg \, {\bf A} \vec{v} \right>_P
\end{align*}
For the curl let $\vec{v} \in  H_C$ and $\vec{w} \in H_S$ so that
\begin{align*}
\left< \curl \vec{v} , \vec{w} \right>_S 
 & = \lbilin {\bf A}^{-1} \, \curl \vec{v} , \vec{w} \rbilin \\
 & = \lbilin \curl \vec{v} , {\bf A}^{-1} \, \vec{w} \rbilin \\
 & = \lbilin \vec{v} , \curl \, {\bf A}^{-1} \vec{w} \rbilin \\
 & = \lbilin {\bf B} \, \vec{v} ,
	{\bf B}^{-1} \, \curl \, {\bf A}^{-1} \vec{v} \rbilin \\
 & = \left< \vec{v} , {\bf B}^{-1} \, \curl \, {\bf A}^{-1} \vec{w} \right>_C
\end{align*}
For the divergence let $\vec{w} \in H_S$ and $ g \in H_V$ so that
\begin{align*}
\left< \divg \vec{w} , g \right>_V
 & = \lbilin a^{-1} \divg \vec{w} , g \rbilin  \\
 & = \lbilin \divg \vec{w} , a^{-1} g \rbilin  \\
 & = - \lbilin \vec{w} , \grad ( a^{-1} g ) \rbilin  \\
 & = - \lbilin {\bf A}^{-1}  \vec{w} , {\bf A} \grad ( a^{-1} g ) \rbilin  \\
 & =   \left< \vec{w} , - {\bf A} \grad ( a^{-1} g ) \right>_S  \\
\end{align*}
For the operator ${\bf A}$ let $\vec{v} \in H_C$ and $\vec{w} \in H_S$ so that
\begin{align*}
\left< {\bf A} \vec{v} , \vec{w} \right>_S  
 & = \lbilin {\bf A}^{-1} {\bf A} \vec{v} , \vec{w} \rbilin \\
 & = \lbilin {\bf A} \vec{v} , {\bf A}^{-1} \vec{v} \rbilin \\
 & = \left< \vec{v} , {\bf A}^{-1} \vec{v} \right>_C
\end{align*}
For the operator $a$ let $ f \in H_P$ and $g \in H_V$ so that
\begin{align*}
\left< a \, f , g \right>_V 
 & = \lbilin a^{-1} a \, f , g \rbilin  \\
 & = \lbilin a \, f , a^{-1} g \rbilin  \\
 & = \left< f , a^{-1} g \right>_P  \\
\end{align*}
Similar arguments give the adjoint operators for operators containing
$b$ and ${\bf B}$.
To keep the notation easy to read it has not been specified whether to use
$a$ or $b$ and whether to use ${\bf A}$ or ${\bf B}$ in the inner products
when computing adjoints. This is clear from the context.

It is now straight forward to compute the adjoints of the second order
operators in Table \ref{General Fundamental Operators}:
\begin{align}
\left( a^{-1} \divg {\bf A} \grad \right)^* 
& = \grad^* \, {\bf A}^* \, \divg^* \, {a^{-1}}^* \nonumber\\
& = a^{-1}\, \divg\, {\bf A}\, {\bf A}^{-1} \,
	{\bf A} \, \grad \, a^{-1} \, a \label{Example 1}\\
& = a^{-1} \, \divg \, {\bf A} \, \grad \,. \nonumber
\end{align}
\begin{align}
\left({\bf A}^{-1} \, \curl \, {\bf B}^{-1} \, \curl \, \right)^*
& = \curl^* \, {{\bf B}^{-1}}^* \, \curl^* \, {{\bf A}^{-1}}^* \nonumber\\
& = {\bf A}^{-1} \, \curl \, {\bf B}^{-1} \, {\bf B} \,
	{\bf B}^{-1} \, \curl \, {\bf A}^{-1} \, {\bf A} \label{Example 2} \\
& = {\bf A}^{-1} \, \curl \, {\bf B}^{-1} \, \curl \, . \nonumber
\end{align}
\begin{align}
\left( {\bf B} \, \grad b^{-1} \, \divg \right)^*
& = \divg^* \, {b^{-1}}^* \, \grad^* \, {\bf B}^* \nonumber \\
& = {\bf B} \, \grad \, b^{-1} \, b \, b^{-1} \,
	\divg \, {\bf B} \, {\bf B}^{-1} \label{Example 3} \\
& = {\bf B} \, \grad \, b^{-1} \, \divg  \nonumber
\end{align}
So these three operators are self-adjoint and similar arguments show that all
operators in Table \ref{General Fundamental Operators} are self-adjoint.

\subsection{Positive and Negative Second Order Operators}
Arguments like those in the previous sections can be use to show that
the second order operators in Table \ref{General Fundamental Operators} are
either positive or negative, those that contain two curl operators are
positive while those that contain a gradient and divergence are negative.

Let $f \in H_P$ and then consider
\begin{align*}
\left< a^{-1} \, \divg \, {\bf A} \, \grad \, f, f \right>_P
& = \left< \divg \, {\bf A} \, \grad \, f, a \, f \right>_V
\\
& = - \left< {\bf A} \, \grad \, f, {\bf A} \, \grad \, a^{-1} \, a \,f \right>_S \\
& = - \left< {\bf A} \, \grad \, f, {\bf A} \, \grad \, f \right>_S \\
& \leq 0
\end{align*}
Let $\vec{v} \in H_C$ and then consider
\begin{align*}
\left< {\bf A}^{-1} \, \curl \, {\bf B}^{-1} \, \curl \, \vec{v} , \vec{v} \right>_C
& = \left<  \curl \, {\bf B}^{-1} \, \curl \, \vec{v} , {\bf A} \, \vec{v} \right>_S \\
& = \left< {\bf B}^{-1} \, \curl \, \vec{v} , {\bf B}^{-1} \, \curl \, {\bf A}^{-1} {\bf A} \vec{v} \right>_C \\
& = \left< {\bf B}^{-1} \, \curl \, \vec{v} , {\bf B}^{-1} \, \curl \, \vec{v} \right>_C \\
& \geq 0
\end{align*}
Let $\vec{w} \in H_S$ and then consider
\begin{align*}
\left< {\bf B} \, \grad b^{-1} \, \divg \, \vec{w} , \vec{w} \right>_S
& = \left< \grad b^{-1} \, \divg \, \vec{w} , {\bf B}^{-1} \, \vec{w} \right>_S \\
& = 
\left< b^{-1} \, \divg \, \vec{w} , \grad^* \, {\bf B}^{-1} \, \vec{w} \right>_P \\
& = - \left< b^{-1} \, \divg \, \vec{w} , b^{-1} \, \divg \, {\bf B} \, {\bf B}^{-1} \, \vec{w} \right>_P \\
& = - \left< b^{-1} \, \divg \, \vec{w} , b^{-1} \, \divg \, \vec{w} \right>_P \\
& \leq 0
\end{align*}
These results capture the important features of the 
second order differential operators.

\newpage \clearpage
\setcounter{equation}{0}
\section{3D Wave Equations With Variable Material Properties
\label{3D Wave Equations}}

Wave equation can be derived from Newton's law and thus have the form
\begin{equation}
\rho \, \frac{d^2 W}{d t^2} = \mathscr{A} \, W \,,
\end{equation}
where $W = W(t,x,y,z)$ is a scalar or vector function, $\rho$ is the density of
the material that the wave is traveling in and $\mathscr{A}$ is a second order
spatial differential operator. Here such equations will written in the form
\begin{equation}
\frac{d^2 W}{d t^2} = \mathscr{B} \, W \,,\,
\mathscr{B} = \frac{1}{\rho} \mathscr{A} \,.
\end{equation}
where $\mathscr{B}$ is one of the second order operators in Table
\ref{General Fundamental Operators}.
The general case provides 12 possible wave equations. Another four equations
are obtained using the operators in \eqref{Linear Combination Operators}.
Note that many of these equations are equivalent but will not be equivalent
in the discrete setting because of the use of two grids in the discretization.
A critical point for wave equations is that the operator $\mathscr{B}$ must be
negative definite so that all eigenvalues of this operator are real and
strictly less than zero. This will guarantee that the solutions are oscillatory.
As will be seen below, if the material properties are constant and trivial
and $c$ is the wave speed then $\mathscr{B}$ will be a $c^2$ times one of
$\divg \grad$, $\grad \divg$ or $\curl \curl$.

Next the second order equations are written as first order systems that are
created using operators generated by going half way around the squares in the
double exact sequence diagram \ref{Exact-Sequences}.  Note that the dependent
variables in the  first order systems are not dimensionless because they come
from spaces containing function with different spatial units. As before,
the first order systems are used to create conservation laws.

Finally, the connection between the notation in the literature and the notation
in this paper for wave, Maxwell and elastic wave equation is discussed.  The
methods used here cannot produce the general elastic wave equation.

\subsection{Second order Wave Equations}

The well known scalar wave equation is given by choosing $W = f \in H_P$ and
then choosing $\mathscr{B}$ to be the operator in the upper left cell of Table
\ref{General Fundamental Operators} to get
\begin{equation}\label{Second Order Scalar}
\frac{d^2 f}{d t^2} = a^{-1} \divg {\bf A} \grad  f \,.
\end{equation}
To understand how this equation relates to the standard wave equation assume
that $a$ is a constant and that $\bf A$ is a consant $A$ times the identity
matrix. In this case $A/a \units d^2/t^2$ that is $A/a = c^2$ where $c$ is the
wave speed. So under these simplifying assumptions the above equation becomes
the standard scalar wave equation
\[
\frac{d^2 f}{d t^2} = c^2 \, \divg \, \grad  f \,.
\]
The initial conditions for the second order equation are $f(0)$ and $f'(0)$.

A vector wave equation can be generated by choosing $W = \vec{w} \in H_S$ and
then choosing $\mathscr{B}$ to be the operator in the lower left cell of
Table \ref{General Fundamental Operators} to get
\begin{equation}
\frac{d^2 \vec{w}}{d t^2} = {\bf B} \, \grad b^{-1} \, \divg \, \vec{w} \,.
\end{equation}
To understand how this equation relates to a standard wave equation assume
that $b$ is a constant and that $\bf B$ is a consant $B$ times the identity
matrix. In this case $B/b \units d^2/t^2$ that is $B/b = c^2$ where $c$ is the
wave speed. So under these simplifying assumptions the above equation becomes
the standard vector wave equation
\[
\frac{d^2 \vec{w}}{d t^2} = c^2 \, \grad \, \divg \, \vec{w} \,.
\]
For this equation the initial conditions are $\vec{w}(0)$ and $\vec{w}'(0)$.

Maxwell's second order equation can be generated by choosing
$W = \vec{v} \in H_C$ and $\mathscr{B}$ the operator in the center left cell of
Table \ref{General Fundamental Operators} to get
\begin{equation}\label{Maxwell's Equation}
\frac{d^2 \vec{v}}{d t^2} =
- {\bf A}^{-1} \, \curl \, {\bf B}^{-1} \, \curl \, \vec{v} \,.
\end{equation}
Under the simplifying assumptions this becomes
\[
\frac{d^2 \vec{v}}{d t^2} = - c^2 \, \curl \, \curl \, \vec{v} \,,\quad
\]
which is Maxwell's second order equation in uniform materials. In total twelve
equations can be created all of which under the simplifying assumptions
reduce to one of the three wave equations just discussed. However, in the
discretization the twelve general equations offer flexibility in modeling
physical problems with variable material properties.

Additional second order wave equations can be made from the two term
second order operators ${\bf VL}_1$, ${\bf VL}_2$, ${\bf VL}_3$, and
${\bf VL}_4$ in \eqref{Linear Combination Operators}. For example
\begin{equation}
\frac{d^2 \vec{v}}{d t^2} =
\grad a^{-1} \divg {\bf A} \vec{v} -
{\bf A}^{-1} \curl {\bf B}^{-1} \curl \vec{v}  \,.
\label{Curl Curl}
\end{equation}
Under the simplifying assumptions this equation becomes
\[
\frac{d^2 \vec{v}}{d t^2} =
c^2 \, \grad \, \divg \, \vec{v} -  c^2 \, \curl \, \curl \vec{v}  \,,
\]
which is a special case of the elastic wave equation
\eqref{dimensionless elastic equation}.  All four equations created this
way reduce to the elastic wave equation under the simplifying assumptions.
Note that the space time units of the dependent variable do not
play a role in the second order equations. The study of first order system
will clarify what the units of $a$, $b$, ${\bf A}$ and ${\bf B}$ must be.

\subsection{First Order Systems and Conserved Quantities
\label{First Order Systems and Conserved Quantities}}

There is a natural way to use diagram chasing to write the second order wave
equations as a system of first order equations and then use this to define 
conserved quantities.  Note that the first order equations contain functions
from different spaces with different spatial units.  The main idea is to
choose two function in diagonally opposite corners of one of the squares
in Figure \ref{Exact-Sequences} and do a diagram chase.  Consequently there
are lots of first order systems!

For example for equation \eqref{Second Order Scalar}, because $f \in H_P$,
choose $\vec{w} \in H_S$ and then set
\begin{equation}
\frac{d \vec{w}}{d t}  = {\bf A} \grad  f 
\,,\quad
\frac{d f}{d t} = a^{-1} \divg \, \vec{w}
\label{System 1}
\end{equation}
The initial conditions for this system are $f(0)$ and $w(0)$
This also gives the vector wave equation
\begin{equation}\label{Second Order Vector}
\frac{d^2 \vec{w}}{d t^2} = {\bf A} \, \grad a^{-1} \divg \, \vec{w}
\end{equation}
for which the initial conditions are $\vec{w}(0)$ and $\vec{w}'(0)$.

The exact sequence \ref{Exact-Sequences} requires that the spatial units the
dependent variable in system \eqref{System 1} are $f \units 1$ and
$\vec{w} \units 1/d^2$. At least is some applications $a$ is a density so
assume $a \units 1/d^3$. If $u_f$ are the units of $f$ and $u_w$ are the
units of $w$ then the two first order equation give
\begin{equation} \label{Units u w}
\frac{u_w}{t} = u_A \frac{1}{d} \, u_f \,,\quad
\frac{u_f}{t} = d^3 \, \frac{1}{d} \, u_w \,,
\end{equation}
which implies
\[
u_A = \frac{1}{d \, t^2}.
\]
This then implies that
\[
\frac{A}{a} \units \frac{d^2}{t^2}
\]
which is required by the constant coefficient case where $A/a = c^2$
where $c$ is the wave speed. Now \ref{Units u w} become 
\[
u_w = \frac{1}{t\, d^2} \, u_f \,,\quad
u_f = t \, d^2 \, u_w \,.
\]

Now if $f \units 1$, that is $u_f = 1$, the $w \units 1/t \, d^2$.  On the
otherhand if $w$ is to be a velocity then $u_w \units d/t$ then it must be
that $u_f \units d^3$ and the $df/dt$ becomes a rate of change of volume.
This needs to be checked against what modelers do.

This system has a conservation law because table \ref{Adjoints} imply that
\[
(a^{-1} \divg )^* = \divg^* (a^{-1})^* = -A \grad a^{-1} a = -A \grad \,,
\]
so this system has the form of the equation discussed in \eqref{Wave System}
and consequently should have a conserved quantity given by 
\[
C = \frac{ \norm{f}_P^2 + \norm{\vec{w}}_S^2 }{2}
\]
This can be checked explicitly:
\begin{align*}
\frac{d C}{d t}
& =
\left< f , \frac{d f}{d t} \right>_P +
\left< \vec{w} , \frac{d \vec{w}}{d t} \right>_S \\
& = \left< f , a^{-1} \divg \, \vec{w} \right>_P +
\left< \vec{w} , {\bf A} \grad  f \right>_S \\
& = \left< (a^{-1} \divg)^* f , \, \vec{w} \right>_S +
\left< \vec{w} , {\bf A} \grad  f \right>_S \\
& = - \left< {\bf A} \, \grad \, f , \vec{w} \right>_S +
\left< \vec{w} , {\bf A} \grad  f \right>_S \\
& = 0 \,.
\end{align*}
As discussed in Section \ref{Harmonic Oscillator}, if $f$ and $\vec{w}$
are solutions of \eqref{System 1} then so are $df/dt$ and $d \vec{w}/dt$
and consequently the classical energy
\[
E  = \frac{\norm{\frac{d f}{d t}}_P^2
+ \norm{\frac{d \vec{w}}{d t}}_S^2}{2} 
   = \frac{\norm{\frac{d f}{d t}}_P^2
+ \norm{{\bf A} \, \grad f }_S^2}{2}
\]
is conserved.

To convert the Maxwell equation \eqref{Maxwell's Equation} to a system assume
that $\vec{v} \in H_C$ in the upper row of \ref{Exact-Sequences} and
that $\vec{w} \in H_C$ in the lower row of \ref{Exact-Sequences} then set
\begin{equation*}
\frac{d \vec{w}}{d t} = {\bf B}^{-1} \, \curl \, \vec{v} \,,\quad 
\frac{d \vec{v}}{d t} = - {\bf A}^{-1} \, \curl \, \vec{w} \,.
\end{equation*}
A conserved quantity is given by
\[
C = \frac{ \norm{\vec{v}}_C^2 + \norm{\vec{w}}_C^2 }{2}
\]
because
\begin{align*}
\frac{d C}{d t} & =
\left< \vec{v} , \frac{d \vec{v}}{d t} \right>_C +
\left< \vec{w} , \frac{d \vec{w}}{d t} \right>_C \\
& =
- \left< \vec{v} , {\bf A}^{-1} \, \curl \, \vec{w} \right>_C +
\left< \vec{w} , {\bf B}^{-1} \, \curl \, \vec{v} \right>_C \\
& = -
\left< {\bf B}^{-1} \curl {\bf A}^{-1} \, {\bf A} \vec{v} , \vec{w} \right>_C +
\left< \vec{w} , {\bf B}^{-1} \, \curl \, \vec{v} \right>_C \\
& = -
\left< {\bf B}^{-1} \curl
\vec{v} , \vec{w} \right>_C +
\left< \vec{w} , {\bf B}^{-1} \, \curl \, \vec{v} \right>_C \\
& = 0
\end{align*}
Additionally,
\[
E = \frac{\norm{\frac{d \vec{w}}{d t}}_C^2 +
          \norm{\frac{d \vec{v}}{d t}}_C^2 }{2}
  = \frac{ \norm{ {\bf B}^{-1} \, \curl \, \vec{v} }_C^2 +
	\norm{ {\bf A}^{-1} \, \curl \, \vec{w} }_C^2 }{2}
  = \frac{ \norm{\curl \, \vec{v} }_S^2 +
	\norm{\curl \, \vec{w} }_S^2 }{2} \,,
\]
is a conserved quantity.
In this case $E$ is essentially the physical energy.

Also first order systems can be made from the second order equations made
from the two term second order operators
${\bf VL}_1$, ${\bf VL}_2$, ${\bf VL}_3$ and ${\bf VL}_4$ given in
\eqref{Linear Combination Operators}. However this requires three first
order equations.  For example, for ${\bf VL}_1$, because
$\vec{v} \in H_C$ let $g \in H_V$ and $\vec{u} \in H_C$ set 
\begin{align*} 
\frac{d g}{d t} & =  \divg {\bf A} \vec{v} \\
\frac{d \vec{u}}{d t} & = {\bf B}^{-1} \curl \vec{v} \\
\frac{d \vec{v}}{d t} & =
	\grad a^{-1} g - {\bf A}^{-1} \curl \vec{u} \,.
\end{align*}
For a conserved quantity set
\[
C = \frac{\norm{\vec{v}}_C^2 + \norm{\vec{u}}_C^2 + \norm{g}_V^2}{2} \\,
\]
so that 
\begin{align*}
\frac{d C}{d t} & = 
  \left< \vec{v} , \frac{d \vec{v}}{d t} \right>_C
+ \left< \vec{u} , \frac{d \vec{u}}{d t} \right>_C
+ \left< g , \frac{d g}{d t} \right>_V \\ 
& =
  \left< \vec{v} , \grad a^{-1} g - {\bf A}^{-1} \curl \vec{u} \right>_C
+ \left< \vec{u} , {\bf B}^{-1} \curl \vec{v} \right>_C
+ \left< g , \divg {\bf A} \vec{v} \right>_V \\ 
& =
\left< \vec{v} , \grad a^{-1} g \right>_C -
\left< \vec{v} , {\bf A}^{-1} \curl \vec{u} \right>_C
+ \left<
{\bf A}^{-1} \curl {\bf B}^{-1} {\bf B} \vec{u} , \vec{v} \right>_C
- \left<
{\bf A}^{-1} {\bf A} \grad a^{-1} g , \vec{v} \right>_V \\ 
& =
\left< \vec{v} , \grad a^{-1} g \right>_C -
\left< \vec{v} , {\bf A}^{-1} \curl \vec{u} \right>_C
+ \left< {\bf A}^{-1} \curl \vec{u} , \vec{v} \right>_C
- \left<
\grad a^{-1} g , \vec{v} \right>_V \\ 
& = 0 \,.
\end{align*}
As before this implies that
\begin{align*}
E & = \frac{
		\norm{\frac{d \vec{v}}{d t}}_C^2 +
		\norm{\frac{d \vec{u}}{d t}}_C^2 +
		\norm{\frac{d g}{d t}}_V^2}{2} \\
  & = \frac{
		\norm{\frac{d \vec{v}}{ d t}}_C^2 +
		\norm{ {\bf B}^{-1} \curl \vec{v}}_C^2 +
		\norm{\divg {\bf A} \vec{v}}_V^2}{2} 
\end{align*}
is conserved.

\subsection{Standard Notation For Wave Equations in Three Dimensions
\label{Section Maxwell Equations}}

\begin{figure}[ht]
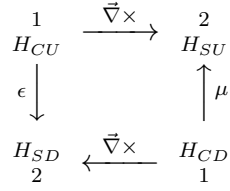

\begin{equation*}
  \begin{CD}
        {\mystack{1}{H_{CU}}}
        @>\curl >> {\mystack{2}{H_{SU}}}
\\
@V{\epsilon }VV @AA{\mu}A @. \\
         {\mystack{H_{SD}}{2}}
        @<\curl << {\mystack{H_{CD}}{1}}
\\
  \end{CD}
\end{equation*}
\caption{Maxwell Exact Sequences}
\label{Maxwell-Sequences}
\end{figure}

\begin{table}[ht]
\begin{center}
\begin{tabular}{|l|c|l|l|}
\hline
quantity & units & name \\
\hline
$\vec{E}$ & $1/d$   & electric field \\
$\epsilon$& $1/d$   & permeability tensor \\
$\vec{D}$ & $1/d^2$ & electric displacement \\
\hline
$\vec{H}$ & $1/d$   & magnetic field \\
$\mu$     & $1/d$   & permittivity tensor \\
$\vec{B}$ & $1/d^2$ & magnetic flux \\
\hline
$\curl$   & $1/d$   & curl operator  \\
\hline
$\vec{J}$ & $1/d^2$     & current \\
\hline
\end{tabular}
\caption{Quantities and their spatial units in the Maxwell
equations \cite{Crain17}.}
\label{Maxwell Units}
\end{center}
\end{table}

The standard notation for the Maxwell equation and the elastic wave equation
are described and the connection to the notation used in this paper are
discussed.  Representing Maxwell's equations \ref{Maxwell Equations} using
diagram chasing uses the center square in Figure \ref{Exact-Sequences} which
is reproduced in Figure \ref{Maxwell-Sequences} using notation appropriate to
Maxwell's equations, that is, by setting ${\bf A} = \epsilon$ and
${\bf B} = \mu$. To make things clearer the space in the upper part of the
diagram are labeled $CU$ and $SU$ with $U$ for upper while the spaces in the
lower row are labeled $SD$ and $CD$ with $D$ for down.
Maxwell's equations can be represented using
$\vec{E} \in H_{CU}$,
$\vec{H} \in H_{CD}$,
$\vec{B} \in H_{SU}$ and
$\vec{D} \in H_{SD}$.
The physical units for Maxwell equations are given in Table \ref{Maxwell Units}.

Diagram chasing with $\vec{E} \in H_{CU}$ and $\vec{H} \in H_{CD}$ gives the
Maxwell equations
\[
\frac{d \vec{H}}{d t} = - \mu^{-1} \curl \vec{E}  \,,\quad
\frac{d \vec{E}}{d t} = \epsilon^{-1} \curl \vec{H} \,.
\]
It is also common to write the Maxwell Equations as
\begin{equation}\label{Maxwell Equations}
\frac{d \vec{B}}{d t} + \curl \vec{E} = 0 \,;\quad
\frac{d \vec{D}}{d t} - \curl \vec{H} = {\vec J} \,.
\end{equation}
\[
\vec{B} = \mu \, \vec{H} \,,\quad \vec{D} = \epsilon \, \vec{E} \,.
\]
Here $\vec{B}$, $\vec{E}$, $\vec{D}$ and $\vec{H}$ are vector functions of
$(x,y,z,t)$ while $\mu$ and $\epsilon$ are symmetric positive definite matrices
that depend only on the spatial variables.  Again the meaning of variables and
their distance units are given in Table \ref{Maxwell Units}.
Set $\vec J = 0$ and then eliminate $\vec{B}$ and $\vec{D}$ from the equation
to get \eqref{Maxwell Equations}.  This system can be written as either of two
second order equations:
\begin{equation}
\frac{d^2 \vec{E}}{d t^2} = 
- \epsilon^{-1} \curl \mu^{-1} \curl \vec{E} \,,\quad
\frac{d^2 \vec{H}}{d t^2} = 
- \mu^{-1} \curl \epsilon^{-1} \curl \vec{H}\,.
\end{equation}

The energy $C$ can be written
\[
C = \int_{\real^3}
\left( \epsilon \, \ \vec{E} \bdot \vec{E} +
       \mu \,        \vec{H}\bdot \vec{H} \right)
\, dx \, dy \, dz \,. \\
\]
The vector identity
\[
\grad \bdot ( \vec{E} \times \vec{H} ) =
(\curl \vec{E}) \bdot \vec{H} - \vec{E} \bdot (\curl \vec{H}) \,
\]
can be used to see that the energy is constant.

The time derivative of $C$ is
\begin{align*}
\frac{d C}{d t}
& = \int_{\real^3}
\left( \epsilon \, \frac{d \vec{E}}{d t} \bdot \vec{E} +
\mu \, \frac{\vec{H}}{d t} \bdot \vec{H} \right)
\, dx \, dy \, dz \,, \\
   & = \int_{\real^3} \left(
	 \curl \vec{H} \bdot \vec{E}
	-\curl \vec{E} \bdot \vec{H}
	  \right) \, dx \, dy \, dz \,, \\ 
   & = - \int_{\real^3} \divg ( \vec{H} \times \vec{E} ) \, dx \, dy \, dz \,, \\
   & = \int_{\real^3} \divg ( \vec{E} \times \vec{H} ) \, dx \, dy \, dz \,, \\
   & = \int_{\real^3} \divg \vec{S} \, dx \, dy \, dz \,, \\
   & = 0 \,.
\end{align*}
The last integral is zero because it was assumed that $\vec{E}$ and $\vec{H}$
are zero far from the origin. Also $\vec{S} = \vec{E} \times \vec{H} $ is
called the Poynting vector which has spatial units $1/d^2$.  The integrand is
the standard energy density confirming $C$ is spatially dimensionless.

\begin{table}[ht]
\begin{center}
\begin{tabular}{|l|c|l|}
\hline
quantity & units & name \\
\hline
$\vec{x}$   & $d$     & spatial position   \\
$\vec{u}$   & $d$     & displacement       \\
$\rho > 0$  & $1/d^3$ & density            \\
$\sigma$    & $1/d^2$ & stress             \\
$e$         & $1$     & strain             \\
$C$         & $1/d^2$ & elastic properties \\
$\lambda>0$ & $1/d$   & Lam\'e parameter   \\
$\mu > 0$   & $1/d$   & Lam\'e parameter   \\
$K$         & $1/d$   & bulk modulus       \\
$Pa$        & $1/d$   & Pascal             \\
\hline
\end{tabular}
\caption{3D Quantities and their spatial units.}
\label{Elastic Units}
\end{center}
\end{table}

For the elastic wave equation for constant material properties note that there
are only two second order differential operators that map vectors to vectors
which are the gradient $\curl \curl $ and $\grad \divg$.  So it is no surprise
that the elastic wave equation \cite{Igel06} is a linear combination of these
operators,
\begin{equation}
\label{Elastic Wave Equation}
\rho \, \frac{d^2 \vec v}{d t^2} =
(\lambda + 2 \, \mu) \grad \divg \vec{v}
- \mu \curl \curl \vec{v} \,,
\end{equation}
where $\mu$ and $\lambda$ are constant scalars with spatial dimension $1/d$.
Note that $\mu = 0$ produces
\begin{equation}
\label{Simple Vector Wave Equation}
\frac{d^2 \vec v}{d t^2} =
\frac{\lambda}{\rho} \grad \divg \vec{v} \,,
\end{equation}
which is a vector wave equation.

If 
\[
a = \sqrt{\frac{(\lambda + 2 \, \mu)}{\rho}} \,,\quad
b = \sqrt{\frac{\mu}{\rho}}
\]
then, like for $c$ in the scalar wave equation, $a$ and $b$ have units $d$,
so the elastic wave equation \eqref{Elastic Wave Equation} can be written in
dimensionless form as
\begin{equation} \label{dimensionless elastic equation}
\frac{d^2 \vec v}{d t^2} =
a^2  \grad \divg \vec{v} - b^2 \curl \curl \vec{v} \,.
\end{equation}
This equation can be converted to a system of three dimensionless first order
equations:
\begin{align*}
\frac{d g}{d t} & = a \, \divg \vec{v} \,; \\
\frac{d \vec{u}}{d t} & = b \, \curl \vec{v} \,; \\
\frac{d \vec{v}}{d t} & = a \, \grad g - b \,\curl \vec{u} \,.
\end{align*}
The quantity
\[
C = \frac{\norm{\vec{v}}^2 + \norm{\vec{u}}^2 + \norm{g}^2}{2} \\,
\]
is conserved because
\begin{align}
\frac{d C}{d t}
& =\left< \frac{d \vec{v}}{d t}, \vec{v} \right> 
+  \left< \frac{d \vec{u}}{d t}, \vec{u} \right> 
+  \left< \frac{d \vec{g}}{d t}, \vec{g} \right> \nonumber \\
& =\left< a \, \grad g - b \,\curl \vec{u}, \vec{v} \right> 
+  \left< b \, \curl \vec{v} , \vec{u} \right> 
+  \left< a \, \divg \vec{v}, \vec{g} \right> \nonumber \\
& = a \, \left< \grad g , \vec{v} \right> 
  - b \, \left< \curl \vec{u}, \vec{v} \right> 
  + b \, \left<  \curl \vec{v} , \vec{u} \right> 
  + a \, \left<  \divg \vec{v}, \vec{g} \right> \nonumber \\
& = - a \, \left< g , \divg \vec{v} \right> 
  - b \, \left<  \vec{u}, \curl \vec{v} \right> 
  + b \, \left<  \curl \vec{v} , \vec{u} \right> 
  + a \, \left<  \divg \vec{v}, \vec{g} \right> \nonumber \\
& = 0 \nonumber
\end{align}

The general elastic wave equations does not fit into the diagram chasing.
This discussion is based on \cite{Etgen87} and a summary of the notation is
given in Table \ref{Elastic Units}.  When there are no external forces the
general elastic wave equation in a material with spatially variable
properties is given by Newton's law applied to the displacements of the
material:
\begin{equation}
\rho \frac{d^2 u_i}{d t^2} =
\sum_{j=1}^3 \frac{ d{\sigma_{i,j}}}{ d {x_j}}
	\,,\quad 1 \leq i \leq 3
\label{General Elastic Wave Equation}
\end{equation} 
where $t$ is time, $\vec{u} = \vec{u}(\vec{x},t)$, are the displacements
of the material, $\rho = \rho(\vec{x})$ is the density of the material,
and $\sigma = \sigma(\vec{x})$ is the symmetric stress tensor:
\[
\sigma_{i,j} = \sum_{k,l = 1}^3 C_{k,l,i,j} \, e_{k,l} \,.
\]
The strains are 
\[
e_{i,j} = \frac{1}{2} 
\left(
\frac{\partial u_i}{\partial x_j} + \frac{\partial u_j}{\partial x_i}
\right) 
\]
which are dimensionless.
The material properties other than density are given by the
$C_{k,l,i,j}$ where $C = C(\vec{x})$ and where $C$ has the
symmetries $C_{k,l,i,j} = C_{l,k,i,j} = C_{k,l,j,i}$.
Consequently $C$ has only 21 independent entries \cite{Etgen87}.  However the
most general second order wave equations generated by diagram chasing
\ref{Linear Combination Operators} have only 13 independent entries 
Even if the matrices in \ref{Linear Combination Operators} are not
assumed symmetric there are only 19 parameters.

\newpage \clearpage
\setcounter{equation}{0}
\section{Mimetic 3D Discretizations \label{Mimetic Discretizations}}

This discussion and notation will follow that in \cite{RobidouxSteinberg2011}.
However, that work was set up to rigorously prove that the discrete operators
in mimetic discretizations have the same properties as the continuum operators
used in vector calculus. Here the focus will be on applying mimetic methods to
physical problems by adding a time variable and its discretization, adding
inhomogeneous and anisotropic material properties and focus on how to use
physical spatial and time units to correctly discretize physical problems.  
Historically, this type of discretization appeared in the Yee grid for
Maxwell's equations \cite{Yee1966} which will be discussed in Section
\ref{E and M section}. This section finishes by describing simulation
programs used to test the results presented here.

A possible  difficulty is that to derive conserved quantities the discrete
versions of 
$a$ and $a^{-1}$,
$b$ and $b^{-1}$,
${\bf A}$ and ${\bf A}^{-1}$ and 
${\bf B}$ and ${\bf B}^{-1}$
must be exact inverses of each other, see
\ref{First Order Systems and Conserved Quantities}.
This is trivial for $a$ and $b$ and also for diagonal ${\bf A}$ 
and ${\bf B}$. The non-diagonal case is important and
has not yet been resolved.

\subsection{Primal and Dual Grids}

\begin{figure}
\begin{center}
\includegraphics[width=4.00in,trim = 0 50 0 150,clip]{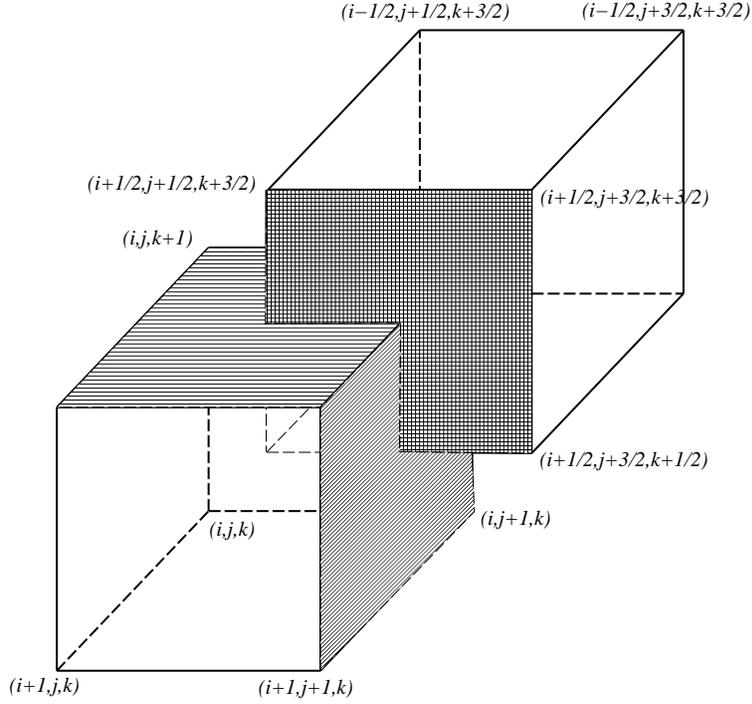}
\end{center}
\caption{The Primal and Dual Grids Taken From \cite{RobidouxSteinberg2011}
\label{Figure-Dual-Grid}}
\end{figure}

\begin{table}
\begin{center}
{\renewcommand{\arraystretch}{1.4}
\begin{tabular}{|c|c|c|}
\hline
primal   & & dual     \\
\hline
nodes   &$\left( i\,\Dx,j\,\Dy,k\,\Dz \right)$   & cells  \\
\hline
        &$\left( \left( i+\half \right) \,\Dx, j \,\Dy, k \,\Dz \right)$          &          \\
edges   &$\left( i \,\Dx, \left( j+\half \right) \,\Dy, k \,\Dz \right)$          & faces    \\
        &$\left( i \,\Dx, j \,\Dy, \left( k+\half \right) \,\Dz \right)$          &          \\
\hline
        &$\left( i \Dx, \left( j+\half \right) \Dy, \left( k+\half \right) \Dz \right) $    &          \\
faces   &$\left( \left( i+\half \right) \Dx, j \Dy, \left( k+\half \right) \Dz \right) $    & edges    \\
        &$\left( \left( i+\half \right) \Dx, \left( j+\half \right) \Dy , k \Dz \right)$    &          \\
\hline
cells   &$\left( \left( i+\half \right) \Dx, \left( j+\half \right) \Dy, \left( k+\half \right) \right) \Dz$ &  nodes \\
\hline
primal  & & dual                                  \\
\hline
\end{tabular}
}
\caption{Notation for the indices of the nodes and the center points of the
edges, faces and cells in the primal and dual grids where
$-\infty < i,j,k < \infty$.
\label{Dual Primal Grids}}
\end{center}
\end{table}

\begin{table}
\begin{center}
{\renewcommand{\arraystretch}{1.4}
\begin{tabular}{|c|c||c|c|}
\hline
units   & primal                         & dual                                 &units  \\
\hline
$1$     &$s_{i,j,k}$                     &$d^\star_{i,j,k}$                     &$1/d^3$\\
\hline
        &$tx_{i+\half, j, k}$             &$nx^\star_{i+\half, j, k}$             &       \\
$1/d  $ &$ty_{i, j+\half, k}$             &$ny^\star_{i, j+\half, k}$             &$1/d^2$\\
        &$tz_{i, j, k+\half}$             &$nz^\star_{i, j, k+\half}$             &       \\
\hline
        &$nx_{i, j+\half, k+\half}$       &$tx^\star_{i, j+\half, k+\half}$       &       \\
$1/d^2$ &$nz_{i+\half, j, k+\half}$       &$ty^\star_{i+\half, j, k+\half}$       &$1/d$  \\
        &$nz_{i+\half, j+\half, k}$       &$tz^\star_{i+\half, j+\half, k}$       &       \\
\hline
$1/d^3$ &$d_{i+\half, j+\half, k+\half}$ &$s^\star_{i+\half, j+\half, k+\half}$ &$1$    \\
\hline
units   & primal                         & dual                                 &units   \\
\hline
\end{tabular}
}
\caption{Notation for the primal and dual, scalar and vector fields where
$-\infty < i,j,k < \infty$.
\label{scalar and vector fields}}
\end{center}
\end{table}

Mimetic discretizations use primal and dual spatial grids as shown in
Figure \ref{Figure-Dual-Grid} and the notation for the nodes, edges, faces
and cells of the grid are given in Table \ref{Dual Primal Grids} while the
notation for scalar and vector fields are given in Table \ref{scalar and
vector fields}.  It is important that the components of vector fields are
not located at the same points in the grid.  All scalar and vector fields
are defined on all of space are assumed to converge to zero far from the
origin. To start all of the grids and all of the scalar and vector functions
will all of three dimensional space.

There two types of scalar fields and also two type of vector fields on
both the primal and dual grids.  On the primal grid there are scalar
fields $s$ that do not have a spatial dimension, vector fields $\vec t$
(for tangent) that have spatial dimension $1/d$, vector fields $\vec n$
(for normal) that have units $1/d^2$, and scalar fields with spatial
dimension $1/d^3$ (as in densities) while the dual grid has the same
types of fields labeled with a superscript star as in $s^\star$.
Note that at each point in the grid there is a value from both the primal and
dual fields.  The spatial dimensions of these fields are different so they
are really different fields.

\subsection{The Discrete Double Exact Sequences}

This section describes the discrete double exact sequences shown in
Figure \ref{Discrete-Exact-Sequences}.  This begins with a description
of the discrete difference operators gradient, curl and divergence on
the primal and dual grids. Next the star or multiplication operators that
describe material properties are discretized.

\begin{figure}
\begin{equation*}
  \begin{CD}
                   \SpaceN
        @>\GRAD >> \SpaceE
        @>\CURL >> \SpaceF
        @>\DIVG >> \SpaceC
\\
@V{a}VV @V{\bf A}VV @A{\bf B}AA @A{b}AA @. \\
                       \StarSpaceC
        @<\DIVGstar << \StarSpaceF
        @<\CURLstar << \StarSpaceE
        @<\GRADstar << \StarSpaceN
\\
  \end{CD}
\end{equation*}
\caption{Discrete Exact Sequences \label{Discrete-Exact-Sequences} }
\end{figure}

\subsubsection{Difference Operators}

The discrete gradient $\GRAD$, curl or rotation $\CURL$ and divergence $\DIVG$
are difference operators on discrete scalar or vector fields.  The formulas
for the dual grid are obtained by making the changes $i \rightarrow i+1/2$,
$j \rightarrow j+1/2$ and $k \rightarrow k+1/2$ in the formulas for the
primal grid.

\noindent {\bf The Gradient:}
If $s \in \SpaceN $ is a discrete scalar field, then its gradient
$\GRAD s = (\GRAD sx, \GRAD sy, \GRAD sz) \in \SpaceE$ is an edge
vector field:
\begin{align}
\GRAD sx_{i+\half, j, k} \equiv &
	\frac{ s_{i+1,j,k}- s_{i,j,k} }{\dx} \,; \nonumber \\
\GRAD sy_{i, j+\half, k} \equiv &
	\frac{s_{i,j+1,k}-s_{i,j,k}}{\dy} \,; \\
\GRAD sz_{i, j, k+\half} \equiv &
	\frac{s_{i,j,k+1}-s_{i,j,k}}{\dz} \,. \nonumber
\end{align}

\noindent {\bf The Curl:}
If $\vec{t} = (tx, ty, tz) \in \SpaceE$ is a discrete edge vector field,
then its curl $\CURL {\vec t} \in \SpaceF$ is a discrete face vector field:
\begin{align}
\CURL tx_{i, j+\half, k+\half} & \equiv
	  \frac{tz_{i, j+1, k+\half}- tz_{i, j, k+\half}}{\dy}
	- \frac{ty_{i, j+\half, k+1}- ty_{i, j+\half, k}}{\dz}
		\,; \nonumber \\
\CURL ty_{i+\half, j, k+\half} & \equiv
	  \frac{tx_{i+\half, j, k+1}- tx_{i+\half, j, k}}{\dz}
	- \frac{tz_{i+1, j, k+\half}- tz_{i, j, k+\half}}{\dx}
		\,; \\
\CURL tz_{i+\half, j+\half, k} & \equiv
	  \frac{ty_{i+1, j+\half, k}- ty_{i, j+\half, k}}{\dx}
	- \frac{tx_{i+\half, j+1, k}- tx_{i+\half, j, k}}{\dy}
		\,. \nonumber 
\end{align}

\noindent {\bf The Divergence:}
If ${\vec n} = (nx, ny, nz ) \in \SpaceF$ is a discrete face vector field,
then its divergence $\DIVG {\vec n} \in \SpaceC $ is a cell scalar field:
\begin{align}
\DIVG {\vec n}_{i+\half, j+\half, k+\half} & \equiv  
	      \frac{nx_{i+1,j+\half,k+\half}- nx_{i,j+\half,k+\half}}{\dx}
			\nonumber \\
	& +  \frac{ny_{i+\half,j+1,k+\half}- ny_{i+\half,j,k+\half}}{\dy}
			\\
	& +  \frac{nz_{i+\half,j+\half,k+1}- nz_{i+\half,j+\half,k}}{\dz}
			\,. \nonumber 
\end{align}

\noindent {\bf The Star Gradient:}
If $s^\star \in \StarSpaceN $ is a discrete star scalar field then
its star gradient $\GRADstar \, s^\star \in \StarSpaceE$ is a star edge
vector field:
\begin{align}
\GRADstar s^\star x_{i, j+\half, k+\half} & \equiv
\frac{s^\star_{i+\half,j+\half,k+\half}-s^\star_{i-\half,j+\half,k+\half} }{\Delta x}; \nonumber \\
\GRADstar s^\star y_{i+\half, j, k+\half} & \equiv
\frac{s^\star_{i+\half,j+\half,k+\half}-s^\star_{i+\half,j-\half,k+\half} }{\Delta x}; \\
\GRADstar s^\star z_{i+\half, j+\half, k} & \equiv
\frac{s^\star_{i+\half,j+\half,k+\half}-s^\star_{i+\half,j+\half,k-\half} }{\Delta x}; \nonumber 
\end{align}

\noindent {\bf The Star Curl:}
If $\vec{t}^\star = (tx^\star, ty^\star, tz^\star ) \in \StarSpaceC$
is a discrete star edge vector field then its curl
$\CURLstar \, \vec{t}^\star \in \StarSpaceF$ is a discrete star face
vector field:
\begin{align}
\CURLstar t^\star x_{i+\half, j, k} & \equiv
          \frac{tz^\star_{i+\half, j+\half, k}-
        tz^\star_{i+\half, j-\half, k}}{\dy}
        - \frac{ty^\star_{i+\half, j, k+\half}-
        ty^\star_{i+\half, j, k-\half}}{\dz}
                \,; \nonumber \\
\CURLstar t^\star y_{i, j+\half, k} & \equiv
          \frac{tx^\star_{i, j+\half, k+\half}-
        tz^\star_{i, j+\half, k-\half}}{\dz}
        - \frac{tz^\star_{i+\half, j+\half, k}-
        tz^\star_{i-\half, j+\half, k}}{\dx}
                \,; \\
\CURLstar t^\star z_{i, j, k+\half} & \equiv
          \frac{ty^\star_{i+\half, j, k+\half}-
        ty^\star_{i-\half, j, k+\half}}{\dx}
        - \frac{tx^\star_{i, j+\half, k+\half}-
        tx^\star_{i, j-\half, k+\half}}{\dy}
                \,. \nonumber
\end{align}

\noindent {\bf The Star Divergence:}
If $\vec{n}^\star = (nx^\star, ny^\star, nz^\star) \in \StarSpaceF$
is a discrete star face vector field then it divergence
$\DIVGstar \vec{n}^\star \in \StarSpaceC$ is a discrete star cell field.
In terms of components
\begin{align}
\DIVGstar {\vec n^\star}_{i, j, k}
                & \equiv
              \frac{nx^\star_{i+\half,j,k}-
        nx^\star_{i-\half,j,k}}{\dx} \nonumber \\
        & +  \frac{ny^\star_{i,j+\half,k}-
        ny^\star_{i,j-\half,k}}{\dy} \\
        & +  \frac{nz^\star_{i,j,k+\half}-
        nz^\star_{i,j,k-\half}}{\dz} \,.  \nonumber
\end{align}

The second order accuracy of the difference operators is confirmed in
{\tt TestAccuracy3.m}.

\subsubsection{Mimetic Properties of Difference Operators}

If $c$ is a constant scalar field then a direct computation 
\cite{RobidouxSteinberg2011} shows that:
\begin{equation}
\GRAD c \equiv 0 \,,\quad
\CURL \GRAD \equiv 0 \,,\quad
\DIVG \CURL \equiv 0 \,,\quad
\GRADstar c \equiv 0 \,,\quad
\CURLstar \GRADstar \equiv 0 \,,\quad
\DIVGstar \CURLstar \equiv 0 \,.
\end{equation}
The code {\tt TestMimetic.m} shows that these relationships are satisfied
up to a small multiple of  {\tt eps}.  These properties are summarized by
saying that the discretization is exact or that the sequences in Figure
\ref{Discrete-Exact-Sequences} are exact, see \cite{RobidouxSteinberg2011}
for a proof that the diagram is exact.

\subsubsection{Discrete Star or Multiplication Operators}

The star operators are multiplication operators that model the material
properties and are given by two positive scalar functions
$a = a(x,y,z)$ and $b=b(x,y,z)$
and two $3 \times 3$ matrix functions
${\bf A} = {\bf A}(x,y,z)$ and ${\bf B} = {\bf B}(x,y,z)$
that are symmetric and positive definite. The spatial dimensions of
$a$ and $b$ must be $1/d^3$ while for ${\bf A}$ and ${\bf B}$ must be $1/d$.

If $s \in \SpaceN$ and $ d^\star = a \, s \in \StarSpaceC$ and
\[
a_{i,j,k} = a(i\,\Dx,j\,\Dy,k\,\Dz )
\]
then 
\[
 d^\star_{i,j,k} = a_{i,j,k} \, s_{i,j,k} \,.
\]
If
if $s^\star \in \StarSpaceN$ and $ d = b \, s^\star \in \StarSpaceC$
and
\[
b_{i,j,k} =
b( i \,\Dx, j \,\Dy, k \,\Dz)
\]
then 
\[
d_{i+\half,j+\half,k+\half} =
b_{i+\half,j+\half,k+\half} \, s^\star_{i+\half,j+\half,k+\half} \,.
\]
The assumption that $a$ and $b$ are not zero implies that these star operators
are invertible.

The discretization of ${\bf A}$, ${\bf B}$, ${\bf A}^{-1}$ and ${\bf B}^{-1}$
is more complicated except for when these matrices are diagonal.
For example let 
\begin{equation}
{\bf A} = 
\left[
\begin{matrix}
Axx & Axy & Axz \\
Ayx & Ayy & Ayz \\
Azx & Azy & Azz
\label{Discrete Material Propreties Matrix}
\end{matrix}
\right] \,,
\end{equation}
where ${\bf A}$ is symmetric, positive definite, and the entries in ${\bf A}$
are functions of $(x,y,z)$.
For ${\vec t} \in \SpaceE$,
${\vec n}^\star = {\bf A} \vec{t} \in \StarSpaceN$ is given by
\begin{align}
nx^\star_{i+\half,j,k} = &
   Axx_{i+\half,j,k} \, tx_{i+\half,j,k} +
   Axy_{i+\half,j,k} \, \overline{ty}_{i+\half,j,k} +
   Axz_{i+\half,j,k} \, \overline{tz}_{i+\half,j,k} \,, \nonumber \\
ny^\star_{i,j+\half,k} = & 
   Ayx_{i,j+\half,k} \, \overline{tx}_{i,j+\half,k} +
   Ayy_{i,j+\half,k} \, ty_{i,j+\half,k} +
   Ayz_{i,j+\half,k} \, \overline{tz}_{i,j+\half,k} \,, \label{A times t} \\
nz^\star_{i,j,k+\half} = & 
   Azx_{i,j,k+\half} \, \overline{tx}_{i,j,k+\half} +
   Azy_{i,j,k+\half} \, \overline{ty}_{i,j,k+\half} +
   Azz_{i,j,k+\half} \, tz_{i,j,k+\half} \,.  \nonumber
\end{align}
To maintain second order accuracy the values with an over bar are
given by average values:
\begin{align*}
\overline{tx}_{i,j+\half,k} & = \frac{
    tx_{i-\half,j,k} + 
    tx_{i+\half,j,k} + 
    tx_{i-\half,j+1,k} + 
    tx_{i+\half,j+1,k}}{4} \,; \\
\overline{tx}_{i,j,k+\half} & = \frac{
    tx_{i-\half,j,k} +
    tx_{i+\half,j,k} +
    tx_{i-\half,j,k+1} +
    tx_{i+\half,j,k+1}}{4} \,; \\
\overline{ty}_{i+\half,j,k} & = \frac{
    ty_{i,  j-\half,k} +
    ty_{i,  j+\half,k} + 
    ty_{i+1,j-\half,k} +
    ty_{i+1,j+\half,k}}{4} \,; \\
\overline{ty}_{i,j,k+\half} & = \frac{
    ty_{i,j-\half,k} +
    ty_{i,j+\half,k} +
    ty_{i,j-\half,k+1} +
    ty_{i,j+\half,k+1}}{4} \,; \\
\overline{tz}_{i,j+\half,k} & = \frac{
    tz_{i,  j,k-\half} + 
    tz_{i,  j,k+\half} + 
    tz_{i,j+1,k-\half} +
    tz_{i,j+1,k+\half}}{4} \,; \\
\overline{tz}_{i+\half,j,k} & = \frac{
    tz_{i,  j,k-\half} +
    tz_{i,  j,k+\half} +
    tz_{i+1,j,k-\half} +
    tz_{i+1,j,k+\half}}{4} \,.
\end{align*}
There are similar formulas for multiplication by
${\bf B}$, ${\bf A}^{-1}$ and ${\bf B}^{-1}$.

\subsection{Discrete Inner Products \label{Discrete Inner Products}}

To study conserved quantities an inner product is needed for each of the eight
linear spaces in the dual exact sequences
Figure \eqref{Discrete-Exact-Sequences}. The inner products will be defined
in terms of four symmetric bilinear forms as in \ref{bilinear form}.
As in the continuum, an important property of the inner products is that
they need to be symmetric, positive definite and importantly dimensionless.

Set $\Delta V = \Delta x \, \Delta y \, \Delta z$.
\noindent
If $s \in \SpaceN$ and $d^\star \in \StarSpaceC$ then
\[
\lbilin s , d \rbilin  = \sum s_{i,j,k} \, d^\star_{i,j,k} \,
\Delta V \,.
\]
\noindent
If $\vec{t} \in \SpaceE$ and $\vec{n}^\star \in \StarSpaceF$ then
\[
\lbilin \vec{t} , \vec{n}^\star \rbilin  = 
\left(
\sum tx_{i+\half,j,k} \, nx^\star_{i+\half,j,k} +
\sum ty_{i,j+\half,k} \, ny^\star_{i,j+\half,k} +
\sum tz_{i,j,k+\half} \, nz^\star_{i,j,k+\half}
\right)
\Delta V \,.
\]

\noindent
If $\vec{n} \in \SpaceE$ and $\vec{t}^\star \in \StarSpaceF$ then
\[
\lbilin \vec{n} , \vec{t}^\star \rbilin  = 
\left(
\sum nx_{i,      j+\half,k+\half} \, tx^\star_{i,      j+\half,k+\half} +
\sum ny_{i+\half,j,      k+\half} \, ty^\star_{i+\half,j,      k+\half} +
\sum nz_{i+\half,j+\half,k}  \,      tz^\star_{i+\half,j+\half,k}
\right)
\Delta V \,.
\]
\noindent
If $g \in \SpaceC$ and $ f^\star \in \StarSpaceN$ then
\[
\lbilin g , f^\star \rbilin = 
\sum g_{i+\half,j+\half,k+\half} \, f^\star_{i+\half,j+\half,k+\half}
\Delta V \,.
\]

The eight inner products are given by the bilinear forms.

\noindent
If $s1, s2 \in \SpaceN$ then
$ \left< s1, s2 \right>_\Nodes =
\lbilin a \, s1 , s2 \rbilin $.

\noindent
If $\vec{t1}, \vec{t2} \in \SpaceE$ then
$\left< \vec{t1}, \vec{t2} \right>_\Edges =
\lbilin {\bf A} \, \vec{t1}, \vec{t2} \rbilin \,. $

\noindent
If $\vec{n1}, \vec{n2} \in \SpaceF$ then
$ \left< \vec{n1}, \vec{n2} \right>_\Faces = 
\lbilin {\bf B}^{-1} \, \vec{n1}, \vec{n2} \rbilin \,.  $

\noindent
If $d1, d2 \in \SpaceC$ then
$
\left< d1, d2 \right>_\Nodes =
\lbilin b^{-1} d1 , d2 \rbilin \,.
$

\noindent
If $s1^\star, s2^\star \in \StarSpaceN$ then
$ \left< s1^\star, s2^\star \right>_\StarNodes =
\lbilin b \, s1^\star, s2^\star \rbilin \,. $

\noindent
If $\vec{t1}^\star, \vec{t2}^\star \in \StarSpaceE$ then
$ \left< \vec{t1}^\star, \vec{t2}^\star \right>_\StarEdges =
\lbilin {\bf B} \, \vec{t1}^\star, \vec{t2}^\star \rbilin \,.  $

\noindent
If $\vec{n1}^\star, \vec{n2}^\star \in \StarSpaceF$ then
$
\left< \vec{n1}^\star, \vec{n2}^\star \right>_\StarFaces = 
\lbilin {\bf A}^{-1} \, \vec{n1}^\star, \vec{n2}^\star \rbilin \,. 
$

\noindent
If $d1^\star, d2^\star \in \StarSpaceC$ then
$
\left< d1^\star, d2^\star \right>_\StarCells = 
\lbilin a^{-1} d1^\star , d2^\star \rbilin \,.
$

\subsection{Adjoint Operators}

The derivation of the adjoints of the discrete operators in the discrete
exact sequences shown in Figure \ref{Discrete-Exact-Sequences} are 
similar to the derivation for the continuum the adjoint operators
defined in Section \ref{Continuum Adjoint Operators}. Again note that the
discrete operators are not mapping of  a space into itself. Some of the
adjoints can easily be guessed by diagram chasing using
\ref{Discrete-Exact-Sequences}:

\begin{align}
\GRAD^*& = - a^{-1} \, \DIVGstar \, {\bf A} &
\GRADstar^* & = - b^{-1} \, \DIVG \, {\bf B}
\,,\nonumber\\
\CURL^* & = + {\bf A}^{-1} \, \CURLstar \, {\bf B}^{-1} &
\CURLstar^* & = + {\bf B}^{-1} \, \CURL \, {\bf A}^{-1} 
\,,\nonumber\\
\DIVG^* & = - {\bf B} \GRADstar  b^{-1} &
\DIVGstar^* & = - {\bf A} \GRAD  a^{-1}
\,,\nonumber\\
{\bf A}^* & = {\bf A}^{-1} &
{\bf B}^* & = {\bf B}^{-1}
\,, \label{DiscreteAdjoints} \\
a^* & = a^{-1} &
b^* & = b^{-1}
\nonumber\,.
\end{align}

When working with systems of first order wave equation adjoints of
products of operators will be needed. For example \eqref{Product Adjoints}
gives 
\[
({\bf A} \GRAD)^*
= \GRAD^* {\bf A}^*
= - a^{-1} \DIVGstar {\bf A} {\bf A}^{-1}
= - a^{-1} \DIVGstar \,.
\]
Also 
\[
(a^{-1} \DIVGstar)^* =
- {\bf A} \GRAD  a^{-1} \, a =
- {\bf A} \GRAD 
\]
Consequently 
\[
(a^{-1} \DIVGstar {\bf A} \GRAD)^*
= ({\bf A} \GRAD)^* (a^{-1} \DIVGstar)^*
= a^{-1} \DIVGstar {\bf A} \GRAD
\]
so this operator is self adjoint.
In terms of inner products, 
if $s \in \in \SpaceN$ and $\vec{n}^\star \in \StarSpaceF$ then
\[
\left< {\bf A} \, \GRAD \, s , \vec{n}^\star \right>_\StarFaces =
- \left< s  , a^{-1} \, \DIVGstar \, \vec{n}^\star \right>_\Nodes \,.
\]
In addition
\[
\left< (a^{-1} \DIVGstar {\bf A} \GRAD) f, f \right>_\Nodes = 
- \left< {\bf A} \GRAD f, {\bf A} \GRAD f \right>_\StarFaces \leq 0 \,,
\]
so this operator is negative.

\subsection{Testing the Differential Operators}

Test show that the differential operators gradient, curl and divergence are
second order accurate and that the curl of the gradient and the divergence
of the curl are zero up to a small multiple of machine epsilon.
The exactness is tested using {\tt TestZero3.m} and the accuracy is tested
using {\tt TestAccuracy3.m} which depends on {\tt ConvergenceRate.m}.
The grids are generated using {\tt Grids3.m} and the differential operators
are given by {\tt Grad3.m}, {\tt Curl3.m} and {\tt Div3.m }. The scalar
and vector fields used for the tests are given by

\begin{tabular}{ll}
{\tt FieldScalarPrimal3.m},  &  \\
{\tt FieldTangentPrimal3.m}, & {\tt FieldGradientPrimal3.m}, \\
{\tt FieldNormalPrimal3.m},  & {\tt FieldCurlPrimal3.m}, \\
{\tt FieldDensityPrimal3.m}, & {\tt FieldDivergenceDual3.m}. 
\end{tabular}

\noindent
Note that there are no star differential operators listed. This is because
only uniform grids are used so the code for the star operators is the
same as for the non-star operators.

\newpage \clearpage
\setcounter{equation}{0}
\section{Discretizing Wave Equations in 3D
\label{Discretize Wave}}
The results in Section Section \ref{Mimetic Discretizations} will be used to
discretize the three dimensional scalar wave equation and Maxwell's wave
equation in three dimensions.

A technical detail in the simulation code is given the spatial discretization
how can one choose a reasonable time step. The estimates in \cite{LinesPB99}
were helpful.

\subsection{The Scalar Wave Equation \label{Scalar Wave 3d}}

The second order scalar wave equation \eqref{Second Order Scalar} can be
written as as first order system as in \eqref{System 1}. But here the notation
will be changed to match that in Section \ref{Mimetic Discretizations}:
\begin{equation*}
\frac{\partial s}{\partial t} = a^{-1} \divg \, \vec{v} \,,\quad
\frac{\partial \vec{v}}{\partial t} = {\bf A} \grad  s  \,,
\end{equation*}
with $s \in H_P$ and $\vec{v} \in H_S$.
This system will be discretized using the operators described in Section
\ref{Mimetic Discretizations} so now let $s \in \SpaceN$ and
$v \in \StarSpaceF$ and then the leapfrog discretization is
\begin{equation}
\frac{s^{n+1} - s^{n}}{\dt} =  a^{-1}\, \DIVGstar v^{n+\half} \,,\quad
\frac{v^{n+\half} - v^{n-\half}}{\dt} = {\bf A} \GRAD s^n \,.
\label{Discrete 3D Wave}
\end{equation}
If $s^0$ and $v^{\half}$ are given then the leapfrog scheme for $n \geq 0$ is
\[
s^{n+1} = s^{n} + \dt \, a^{-1}\, \DIVGstar v^{n+\half} \,,\quad
v^{n+\thalf} = v^{n+\half} + \dt \, {\bf A} \GRAD s^{n+1} \,.
\] 
If $s^0$ and $v^0$ are given then a second order accurate value
for $v^\half$ is given by the Taylor series
\begin{align}
v^{\half} & =  v(\Delta t /2) \nonumber \\
               & = v(0) 
                 + v'(0) \, \frac{\Delta t}{2}
                 + \frac{v''(0)}{2}  \, 
               \left(\frac{\Delta t}{2}\right)^2 \nonumber\\
               & = v(0)
                 + {\bf A} \grad s(0) \frac{\Delta t}{2}
                 + \frac{{\bf A} \grad \, s'(0)}{2}
                   \left(\frac{\Delta t}{2}\right)^2 \nonumber\\
               & = v(0)
                 + {\bf A} \grad s(0) \frac{\Delta t}{2}
                 + \frac{1}{2} {\bf A} \grad \, a^{-1} \divg \vec{v}(0)
                   \left(\frac{\Delta t}{2}\right)^2
\label{Taylor Series Initial}
\end{align}

The discretization of the first order system gives a discretization of the
second order scalar wave equation \eqref{Second Order Scalar} as
\begin{equation}
\frac{s^{n+1} - 2\, s^{n} + s^{n-1}}{\dt^2}
= a^{-1}\, \DIVGstar {\bf A} \GRAD s^n \,, 
\end{equation}
and as discussed in Section \ref{3D Wave Equations} on continuum
wave equations.
This is also a discretization of the vector wave equation
\begin{equation}
\frac{v^{n+\thalf} - 2\, v^{n+\half} + v^{n-\half}}{\dt^2}
= {\bf A} \GRAD a^{-1}\, \DIVGstar v^{n+\half} \,. 
\end{equation}

The results in Section \ref{ODEs} give two conserved quantities
for the discretization. To see this using \eqref{Conserved A* half}
and \eqref{Conserved A full}  set
$A$ to $a^{-1} \, \DIVGstar$,
$A^xx$ to $-{\bf A} \GRAD$,
$g^{n+\half}$ to $v^{n+\half}$ and
$f^n$ to $s^n$
to get the conserved conserved quantities
\begin{equation}
C^n =
 \norm{s^{n}}_\Nodes^2
 + \norm{\frac{v^{n+1/2} + v^{n-1/2}}{2}}_\StarFaces^2
 - \frac{\Delta t^2}{4} \norm{ {\bf A} \GRAD \, s^{n} }_\StarFaces^2 \,,
\end{equation}
\begin{equation}
C^{n+1/2} =
\norm{v^{n+1/2}}_\StarFaces^2
+ \norm{\frac{s^{n+1} + s^{n}}{2}}_\Nodes^2
- \frac{\Delta t^2}{4} \norm{ a^{-1} \, \DIVGstar \, v^{n+1/2} }_\Nodes^2 \,.
\end{equation}
These formulas agree with those derived in detail in Appendix
\ref{Appendix Details}.

\subsection{Testing the Wave Equation Codes}

The test were done in stages, first for trivial material properties with
{\tt Wave3DTMP.m}. For this code the energies $C_n$ and $C_{n+\half}$ are
constant to a small multiple of {\tt eps} and the solutions are again
forth order accurate for the example tested.

Currently working on {\tt Wave3DGMP.m} to model more complex material
properties. Can do constant $a$ and constant diagonal ${\bf A}$.

\subsection{Maxwell's Equations \label{E and M section}}

This is to be redone. xxx

Assume that ${\vec E} \in \SpaceE$ and $\vec{H} \in \StarSpaceE$ so that,
using the notation in the previous sections, the Maxwell system
\ref{Maxwell Equations} will be discretized as 
\[
\frac{{\vec E}^{n+1} - {\vec E}^n}{\dt} = \epsilon^{-1} \, \CURLstar \vec{H}^{n+\half} \,,\quad
\frac{\vec{H}^{n+\half} - \vec{H}^{n-\half}}{\dt} = - \mu^{-1} \, \CURL {\vec E}^n \,.
\]
where in Exact Sequence diagram \eqref{Discrete-Exact-Sequences}
$A = \epsilon$ and $B = \mu$. Here $\epsilon$ and $\mu$
can be symmetric positive definite matrices.  If
$\vec{E} = (Ex, ER,Ez)$ and  $\vec{H} = (Hx, Hy,Hz)$
then Table \ref{Dual Primal Grids} shows that $Ex$ and $Hx$ are indexed as
\[
Ex^n_{i+\half,j,k} \,,\quad Hx^{n+\half}_{i,j+\half,k+\half}
\]
just as in Yee's paper \cite{Yee1966}.  If $\vec{E}^0$ and $\vec{H}^\half$
are given then the leapfrog scheme for $n \geq 0$ is
\[
{\vec E}^{n+1} = {\vec E}^n
+ \dt \, \epsilon^{-1} \, \CURLstar \vec{H}^{n+\half} \,,\quad
\vec{H}^{n+3/2} = \vec{H}^{n+\half}
- \dt \, \mu^{-1} \, \CURL {\vec E}^{n+1} \,.
\]

Using a similar argument, it is easy to see that
\[
C_{n+1/2} =
\norm{\frac{\vec{E}^{n+1} + \vec{E}^{n}}{2}}_\Edges^2
+ \norm{ \vec{H}^{n+1/2}}_\StarEdges^2 
- \frac{\Delta t^2}{4}
\norm{ \epsilon^{-1} \, \CURLstar \, \vec{H}^{n+1/2} }_\Edges^2 
\]
is a conserved quantity and that
\[
C_{n+1/2} \geq
\norm{\frac{\vec{E}^{n+1} + \vec{E}^{n}}{2}}_\Edges^2 +
\left( 1 - \frac{\Delta t^2}{4} \norm{\epsilon^{-1} \, \CURLstar}^2 \right)
\norm{ \vec{H}^{n+1/2}}_\StarEdges^2  \,.
\]
So $C_{n+1/2} \geq 0$ for $\Delta t$ sufficiently small provided
$\norm{\epsilon^{-1} \, \CURLstar}$ is finite.

Also
\[
C_n =  \norm{\vec{E}^{n}}_\Edges^2 
 -\frac{\Delta t^2}{4} \norm{ \mu^{-1} \CURL \, \vec{E}^{n} }_\StarFaces^2 
 +\norm{\frac{\vec{H}^{n+1/2} + \vec{H}^{n-1/2}}{2}}_\StarFaces^2 \,.
\]
is a conserved quantity and 
\[
\norm{C_n} \ge	
\left(1 - \frac{\Delta t^2}{4}
\norm{\mu^{-1} \CURL}^2\right)
\norm{\vec{E}^n}_\Edges^2
+ \norm{\frac{\vec{H}^{n+1/2} + \vec{H}^{n-1/2}}{2}}_\StarFaces^2 \,.
\]
So $\norm{C_n}$ is positive for sufficiently small $\Delta t$ if
$\norm{ \mu^{-1} \CURL \, \vec{E}^{n} }$ is finite.
Also, these formulas agree with those derived in detail in Appendix
\ref{Appendix Details}.

The codes {\tt Maxwell.m} and {\tt MaxwellStar.m} confirm that our algorithms
conserve $C_{n+1/2}$ and $C_n$ to two parts in $10^{16}$. Additionally,
the divergence of the curl of the electric and magnetic fields are constant
to one part in $10^{14}$ when there are no sources.

\newpage \clearpage
\setcounter{equation}{0}
\section{Implementation in 2D \label{Implementaion 2D}}

The goal is to illustrate how scalar and vector functions are discretized
in 2D and then show how the gradient and divergence are discretized.
This will provide some intuition on how to discretize functions and
differential operators in 3D.  For 3D the primal and dual grids are shown
in Figure \ref{Figure-Dual-Grid} and described in Table \ref{Dual Primal Grids}.
In some ways the 2D discretization is more complicated than the 3D 
because in 3D cells have edges and faces while in 2D the edge of a
cell can correspond to either an edge or a face in 3D.

The implementation of mimetic finite difference discretizations in a bounded
region can be annoying due to the two staggered grids and indices that contain
a half. Additionally, the components of vector fields are discretized
at different points. There are also problems at the boundary that the
star operators will be used fix.  Initially the star operators will be
trivial, except at the boundary. The simulation region will be the unit square.

For 3D a detailed description of the mimetic discretization is given
in Section \ref{Mimetic Discretizations}. In particular the left most
box in Figure \ref{Discrete-Exact-Sequences} will be used for the 2D
discretization.

To facilitate programming, all indices are integers greater than zero,
for example $i \geq 1$ and $j \geq 1$.

\subsection{2D Exact solution}
\begin{align*}
s & = \sqrt{m^2+n^2} \\
u(x,y,t) & = \cos(c \, s \, \pi \, t )
\sin(m \pi x) \sin(n \pi y) \\
v(x,y,t) & = \frac{1}{s} \, \sin(c \, s \, \pi \, t)
     (m \, \cos(m \pi x) \, \sin(n \pi y), n \, \sin(m \pi x), \, \cos(n \pi y))
\end{align*}

\subsection{The 2D primal and dual grids}

\begin{figure}
\begin{center}
\includegraphics[width=5.00in,trim=100 247 100 250,clip]{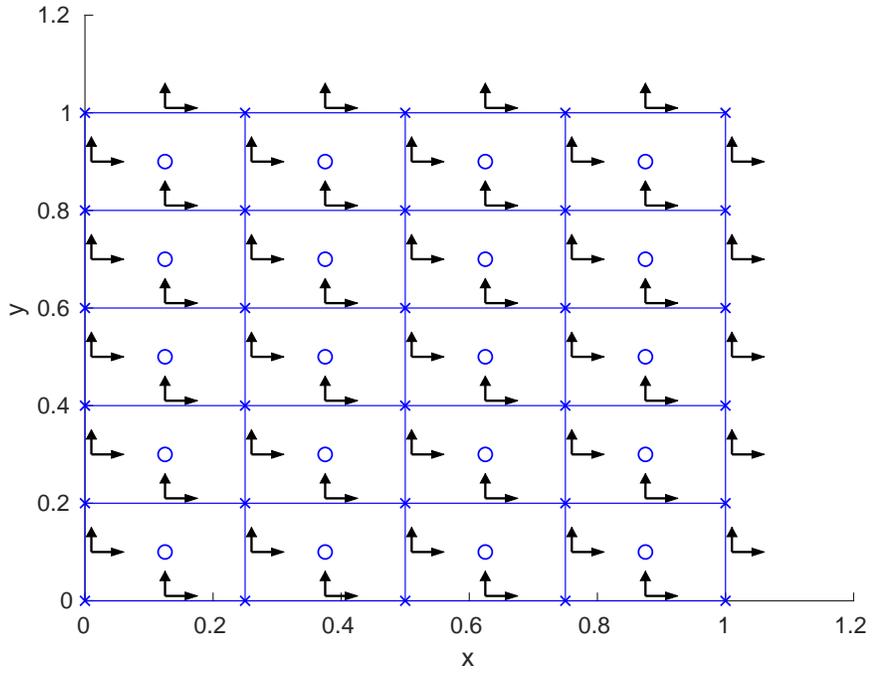}
\end{center}
\caption{2D primal grid tangent and normal vector fields, $Nx = 4$, $Ny = 5$
\label{Tangent Normal Primal}}
\end{figure}
\begin{figure}
\begin{center}
\includegraphics[width=5.00in,trim=100 247 100 250,clip]{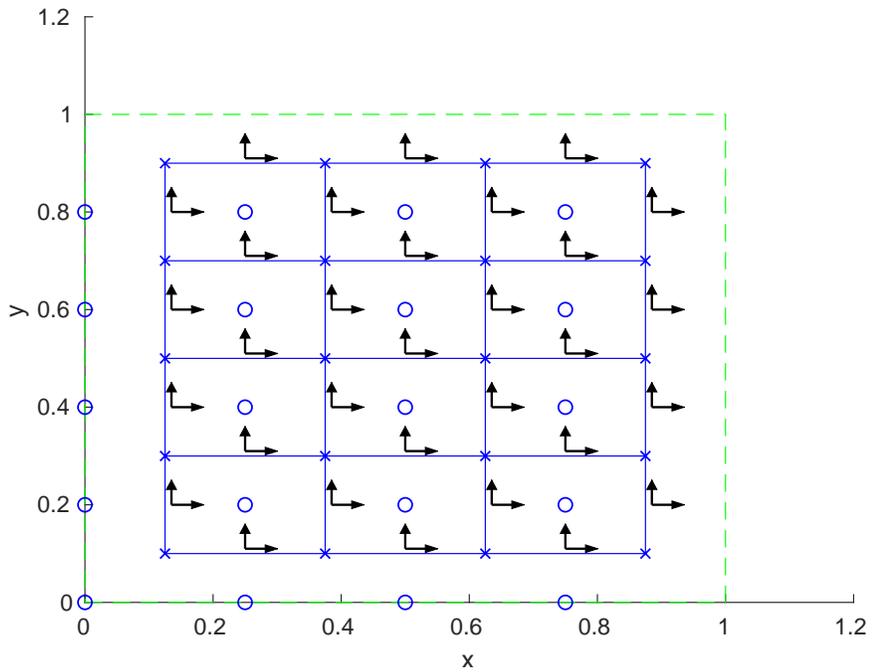}
\end{center}
\caption{2D dual grid tangent and normal vector fields, $Nx = 4$, $Ny = 5$
\label{Tangent Normal Dual}}
\end{figure}

Let $Nx,Ny$ be positive integers and then set $dx = 1/Nx$, $dy = 1/Ny$.
The primal grid nodes (cell corners) are given by
\begin{align*}
(xp(i,j), yp(i,j)) & = ((i-1) \, dx, (j-1) \, dy) \,,\quad
1 \leq i \leq Nx+1 \,,
1 \leq j \leq Ny+1 \,,
\end{align*}
and the dual grid nodes are given by
\begin{align*}
(xd(i,j), yd(i,j)) & = ( (i-1/2) \, dx, (j-1/2) \, dy) \,, \quad
1 \leq i \leq Nx \,,
1 \leq j \leq Ny \,.
\end{align*}
Note that in the primal grid $Nx$ is the number of cells in
$x$ direction and $Ny$ is the number of cells in the $y$ direction.  For
$Nx = 4$ and $Ny = 5$ the positions of scalar and vector function on the
primal and dual grids are illustrated in Figures \ref{Tangent Normal Primal}
and \ref{Tangent Normal Dual}, see {\tt FigurePrimalDual2.m}.
Compare these to Figure \ref{Figure-Dual-Grid} for a 3D grid.

Important points are that a primal grid cell center is given by a dual
grid node and the dual grid centers are given by given by the {\em interior}
primal grid nodes. Additionally the location of dual grid tangent vectors
are the same as the location of the primal grid interior normal vectors
and the position of the dual grid normal vectors are the same as the position
of the interior grid tangent vectors. On the boundary of the primal grid the
positions of points and vectors do not correspond to anything in the dual grid.
As will be seen this is important for representing boundary conditions for
partial differential equations.

\subsection{Discretizing continuum functions on the 2D grids}

\begin{figure}
\begin{center}
\begin{tabular}{c}
\includegraphics[width=6.00in,height=8.00in,trim=100 240 110 215]{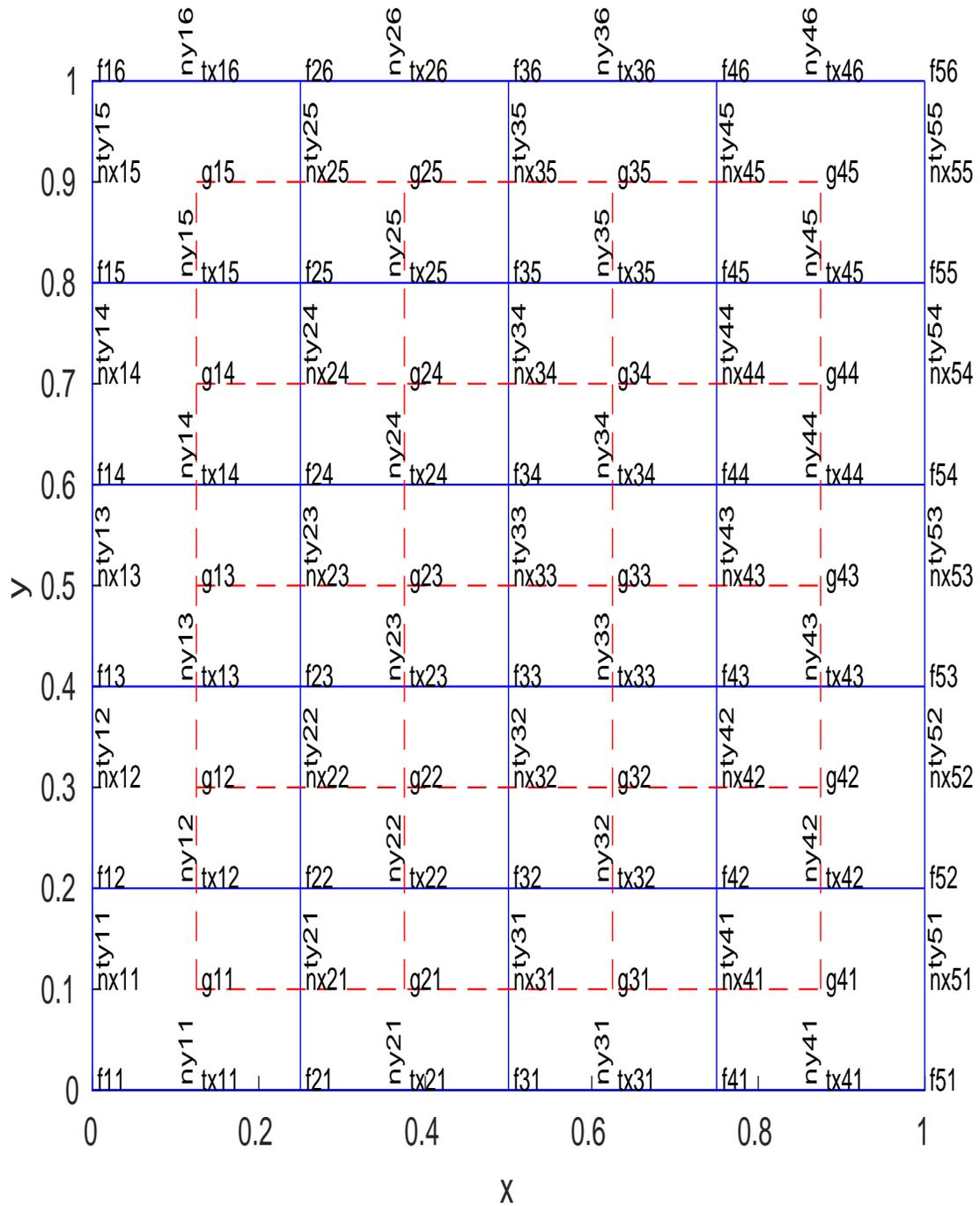}
\end{tabular}
\end{center}
\caption{Scalar and vector fields on the primal grid.
\label{Scalar Vector Primal}}
\end{figure}

\begin{figure}
\begin{center}
\begin{tabular}{c}
\includegraphics[width=6.00in,height=8.00in,trim=100 240 110 215]{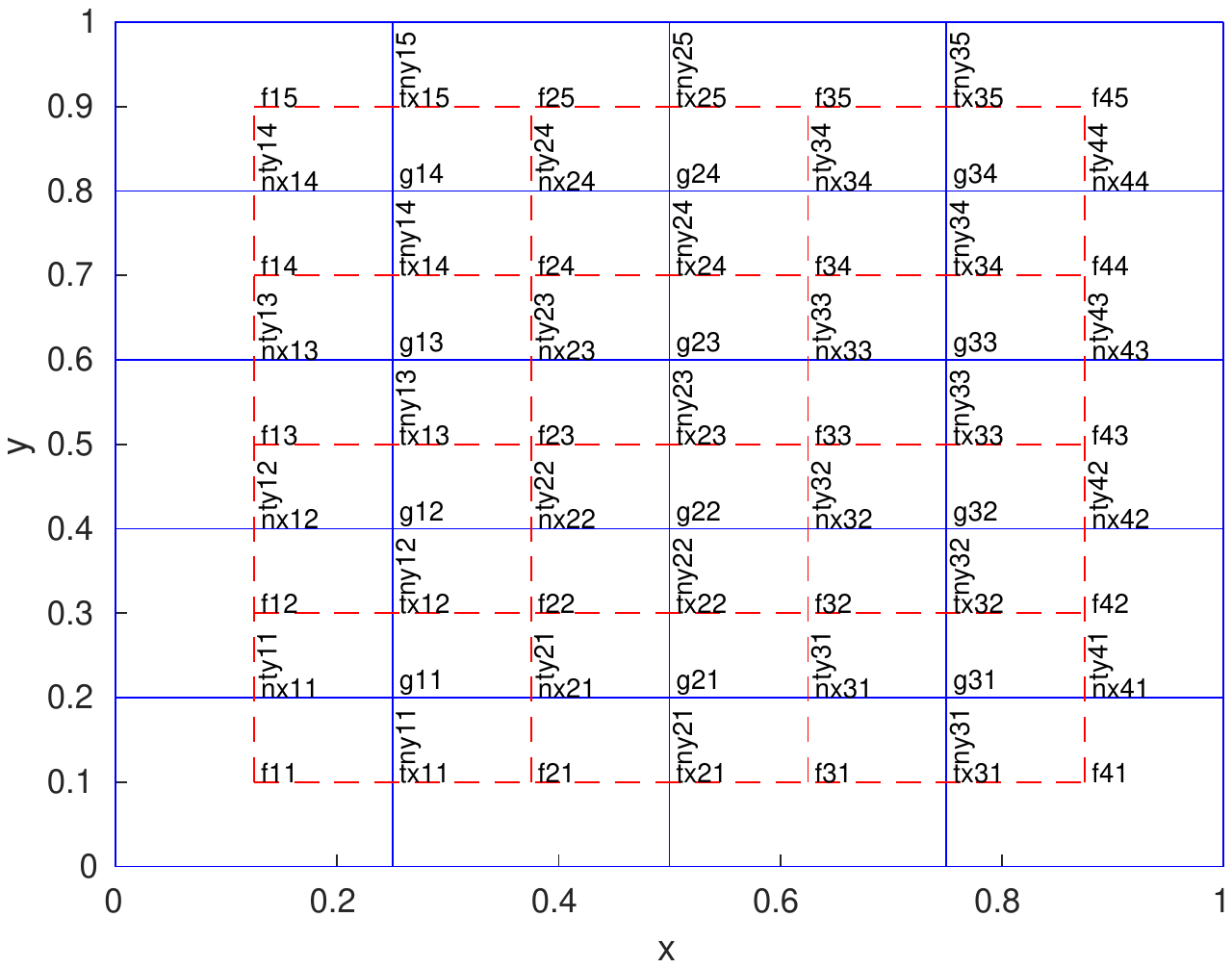}
\end{tabular}
\end{center}
\caption{Scalar and vector fields on the dual grid,
\label{Scalar Vector Dual}}
\end{figure}

Scalar functions will be discretized at either the nodes (corners) or
cell centers of the primal and dual grids while vector fields will be
discretized at the centers of the  of the edges of the cells as illustrated
in Figures \ref{Tangent Normal Primal} and \ref{Tangent Normal Dual} for grids
in a region that is a unit square (see {\tt FigureDetailsPrimalDual2.m}).
There are a total of eight types of discretized functions.
In the text, function names on the primal grid have a $p$ appended as in
$fp$ while functions on the dual grid have a $d$ appended as in $fd$.
In the figures which grid the functions are on is clear.

For the primal grid there are four cases.  A scalar function $f(x,y)$
with spatial weight $1$ is discretized at primal cell nodes:
\begin{equation}
fp(i,j) = f(xp(i),yp(j)) \,,\quad i \leq Nx+1 \,,\quad j \leq Ny+1 \,.
\end{equation}
A scalar function $g(x,y) $ with spatial weights $1/d^3$ is discretized at
primal grid cell centers:
\begin{equation}
gp(i,j) = g(xd(i),yd(j)) \,,\quad i \leq Nx \,,\quad j \leq Ny \,.
\end{equation}
A vector function $(tx(x,y),ty(x,y)$ with spatial weight $1/d$ (tangent on
primal grid) is discretized at cell edge centers:
\begin{align}
txp(i,j) = & tx(xd(i),yp(j)) \,,\quad i \leq Nx   \,,\quad j \leq Ny+1 \,;\\
typ(i,j) = & ty(xp(i),yd(j)) \,,\quad i \leq Nx+1 \,,\quad j \leq Ny \,. 
\end{align}
A vector function $(nx(x,y),ny(x,y)$ with spatial weight $1/d^2$ (normal on
primal grid) is also discretized at cell edge centers:
\begin{align}
nxp(i,j) = & nx(xp(i),yd(i)) \,,\quad i \leq Nx+1 \,,\quad j \leq Ny \,;\\
nyp(i,j) = & ny(xd(i),yp(j)) \,,\quad i \leq Nx   \,,\quad j \leq Ny+1 \,. 
\end{align}

For the dual grid there are also four cases.  A scalar function $f(x,y)$ with
spatial weight $1$ is discretized at the dual cell nodes:
\begin{equation}
fd(i,j) = f(xd(i),yd(j)) \,,\quad i \leq Nx \,,\quad j \leq Ny \,.
\end{equation}
A scalar function $g(x,y) $ with spatial weights $1/d^3$ is discretized at
dual grid cell centers nodes:
\begin{equation}
gd(i,j) = g(xp(i+1),yp(j+1)) \,,\quad
	i \leq Nx-1 \,,\quad j \leq Ny-1 \,.
\end{equation}
A vector function $(tx(x,y),ty(x,y)$ with spatial weight $1/d$ (tangent on
dual grid) is discretized at cell edge centers:
\begin{align}
txd(i,j) = & tx(xp(i+1),yd(j)) \,,\quad i \leq Nx-1 \,,\quad j \leq Ny \,;\\
tyd(i,j) = & ty(xd(i),yp(j+1)) \,,\quad i \leq Nx \,,\quad j \leq Ny \,. 
\end{align}
A vector function $(nx(x,y),ny(x,y)$ with spatial weight $1/d^2$ (normal
on the dual grid) is discretized at cell face centers:
\begin{align}
nxd(i,j) = & nx(xd(i),yp(j+1)) \,,\quad i \leq Nx \,,\quad j \leq Ny-1 \,;\\
nyd(i,j) = & ny(xp(i+1),yd(j)) \,,\quad i \leq Nx-1 \,,\quad j \leq Ny \,. 
\end{align}

Examples of 2D scalar and vector fields can be generated using
{\tt ScalarField2p.m},
{\tt ScalarField2d.m},
{\tt TangentField2p.m},
{\tt TangentField2d.m},
{\tt NormalField2p.m},
{\tt NormalField2d.m},
{\tt DensityField2p.m},
{\tt DensityField2d.m}
and tested using {\tt TestFields2pd.m}

\subsection{The Star Operators}

For simplicity the discussion of star operators will begin for material
properties that are constant. These star operators will map quantities
defined on the primal grid to the dual grid while their inverses do the
opposite.  For constant and isotropic materials the star operators are
multiplication by a constant for quantities that are defined at the same
point in the grids.  For scalar variables this is done by multiplication
by a constant $a>0$ with spatial dimension $1/d^3$. For vectors a constant
diagonal matrix
\[
{\bf A} =
\left[
\begin{matrix}
A11 & 0  \\
0 & A22 \\
\end{matrix}
\right]
\]
with $A11>0$ and $A22>0$ spatial dimension $1/d$ will be used.
If $A11 \neq A22$ then there is a simple anisotropy.  Away from the
boundaries of the region, the two star operators are inverses of each
other.  For boundary value problems, the mismatch between the sizes of the
primal and dual grids will be used to represent the boundary conditions.
The star operators will be given in pairs that are essentially inverse of
each other, first the mapping from the primal grid to the dual grid and then
the inverse.

If $gd = \star fp$ then
\begin{equation}
gd(i,j) = a \,fp(i+1,j+1)
\,,\quad 1 \leq i \leq N_x-1
\,,\quad 1 \leq j \leq N_y -1
\,.
\end{equation}
If $fp = \star gd$ then
\begin{equation}
fp(i,j) = \frac{1}{a} \, gd(i-1,j-1)
\,,\quad 2 \leq i \leq N_x
\,,\quad 2 \leq j \leq N_y
\,.
\end{equation}
In this case $fp$ is not defined on the boundary of the primal grid.

If $fd = \star gp$ then
\begin{equation}
gp(i,j) = a\,fd(i,j)
\,,\quad 1 \leq i \leq N_x
\,,\quad 1 \leq j \leq N_y 
\,.
\end{equation}
If $gp = \star fd$ then
\begin{equation}
gp(i,j) = \frac{1}{a}\,fd(i,j)
\,,\quad 1 \leq i \leq N_x
\,,\quad 1 \leq j \leq N_y 
\,.
\end{equation}

If $(txd,tyd) = \star (nxp,nyp)$ then
\begin{align}
txd(i,j) &= \frac{1}{A11} \, nxp(i+1,j)
\,,\quad 1 \leq i \leq N_x-1
\,,\quad 1 \leq j \leq N_y
\,. \nonumber \\
tyd(i,j) &= \frac{1}{A22} \, nyp(i,j+1)
\,,\quad 1 \leq i \leq N_x
\,,\quad 1 \leq j \leq N_y-1
\,.
\end{align}

If $(nxp,nyp) = \star (txd,tyd)$ then
\begin{align}
nxp(i,j) &= A11 \, txd(i-1,j)
\,,\quad 2 \leq i \leq N_x
\,,\quad 1 \leq j \leq N_y
\,. \nonumber \\
nyp(i,j) &= A22 \, tyd(i,j-1)
\,,\quad 1 \leq i \leq N_x
\,,\quad 2 \leq j \leq N_y
\,.
\end{align}

If $(nxd, nyd) = \star(txp,typ)$ then
\begin{align}
nxd(i,j) &= A11 \, txp(i,j+1)
\,,\quad 1 \leq i \leq N_x
\,,\quad 1 \leq j \leq N_y-1
\,. \nonumber \\
nyd(i,j) &= A22 \, typ(i+1,j)
\,,\quad 1 \leq i \leq N_x-1
\,,\quad 1 \leq j \leq N_y
\,. \label{nxd star txp}
\end{align}

If $(txp,typ) = \star (nxd, nyd)$ then
\begin{align}
txp(i,j) &= \frac{1}{A11} \, nxd(i,j-1)
\,,\quad 1 \leq i \leq N_x
\,,\quad 2 \leq j \leq N_y
\,. \nonumber \\
typ(i,j) &= \frac{1}{A22} \, nyd(i-1,j)
\,,\quad 2 \leq i \leq N_x
\,,\quad 1 \leq j \leq N_y
\,.
\end{align}

\subsection{2D Differential Operators}

The 2D discrete differential operators are the gradient and divergence
which are given by {\tt Grad2p.m} and {\tt Grad2d.m} and {\tt Div2p.m}
and {\tt Div2d.m}. All the differential operators use centered differences.
The operators on the primal and dual grids differ in their indexing.
The programs {\tt TestGradDiv2p.m} and {\tt TestGradDiv2d.m} show that
the gradient and divergence operators are second order accurate on the
primal and dual grids.

\subsection{Scalar Wave Equation}

The discretization of the 2D scalar wave equation is the same as the
discretization of 3D scalar wave equation given in \eqref{Discrete 3D Wave}
and has been implemented in {\tt Wave2D.m}. The wave equation is written
as a first order system of differential equations that are discretized
using staggered space-time grids and {\tt Grad2p} and {\tt Div2d}.
In general the approximate solutions of the discrete wave equation are
second order accurate. For some cases the solutions are forth order accurate
and there is at least one example where the solution is exact
(set $m1 = n1 = 1$ and eliminate the $m2$, $n2$ part of the test
solution to see this).
Both conservation laws are constant to at least 1 part in $10^{15}$
when the star operator is trivial.

What star was used in Wave2D.m? See \eqref{nxd star txp}.

\subsection{Boundary Conditions}

The needs a rewrite.

Figures \ref{Scalar Vector Primal} and \ref{Scalar Vector Dual} illustrates
the position of the scalar field $f$ and the scalar function $g$ which is
the Laplacian of $f$, that is the divergence of the gradient of $f$.
It also illustrates the positions of the boundary conditions which must
specify $g$ on the boundary:
\begin{align*}
g_{1,j}     \,,\, & 1 \leq j \leq Ny+1 \,\;\\
g_{Nx+1,j}  \,,\, & 1 \leq j \leq Ny+1 \,\;\\
g_{i,1}     \,,\, & 1 \leq i \leq Nx+1 \,\;\\
g_{i,Ny+1,} \,,\, & 1 \leq i \leq Nx+1 \,.
\end{align*}
Note that the values of $g$ are defined twice at the corner points
of the region:
$(1,1)$; $(1,Ny+1)$; $(Nx+1,1)$ and $(Nx+1,Ny+1)$.
In fact in standard discrete boundary value problems these values of $g$
are not needed and can be assigned any value. On the other hand creating
a data structure that doesn't have these values creates a programming mess.
The figures were generated using {\tt Figure2DDiv.m}, {\tt Figure2DGrad.m}
and {\tt Figure2DLap.m}.

The typical boundary condition is of mixed or Robin type, that is,
\[
\alpha \, \vec{n} \bdot \vec{v} + \beta f = \gamma \,,
\]
which in the discrete setting becomes
\begin{align*}
\text{at  } y = 0 \quad &
\alpha_{1,j} \, vy_{1,j} +\beta_{1,j} \, g_{1,j} =
\gamma_{1,j} \,; \\
\text{at  } y = 1 \quad &
\alpha_{Nx+1,j} \, vy_{Nx+1,j} +\beta_{Nx+1,j} \, g_{Nx+1,j} =
\gamma_{Nx+1,j} \,; \\
\text{at  } x = 0 \quad &
\alpha_{i,1} \, vy_{i,1} +\beta_{i,1} \, g_{i,1} =
\gamma_{i,1} \,; \\
\text{at  } x = 1 \quad &
\alpha_{i,Ny+1} \, vy_{i,Ny+1} +\beta_{i,Ny+1} \, g_{i,Ny+1} =
\gamma_{i,Ny+1} \,. 
\end{align*}
These equations can be trivially solved for
$ g_{i,1} $, $ g_{i,Ny+1} $, $ g_{1,j} $, $ g_{Nx+1,j} $ 
which could give two different values of $f$ for the corner points.
In the typical explicit time stepping algorithms for wave equations,
these values are never used.

Check this, probably not correct.
For Dirichlet boundary conditions
the program {\tt Wave2D.m} confirms that the solution of the 2D wave
equation is second order accurate while the program {\tt Wave2DExact.m}
illustrates some cases where the solutions are accurate up to some
small multiple of {\tt eps}.

%
\newpage \clearpage
\bibliography{citations}
\begin{appendix}
\newpage \clearpage
\setcounter{equation}{0}
\section{Energy Preserving Discretizations of the Harmonic Oscillator
\label{C-N} }

Here the well-known fact that the Crank-Nicholson discretization conserves
the discrete analog of the energy for the harmonic oscillator is shown.
It is also shown that the methods introduced in \cite{WanBihloNave2015}
produce a discretization that is equivalent to the Crank-Nicholson
discretization.

\subsection{Conserving the Simple Energy}

The Crank-Nicholson discretization does preserve the simple
energy \eqref{Simple Conserved}:
\[
\frac{ u_{n+1}-u_n }{\Delta t} = - \omega \, \frac{v_{n+1}+v_n}{2} \,,\quad
\frac{ v_{n+1}-v_n }{\Delta t} =   \omega \, \frac{u_{n+1}+u_n}{2} \,.
\]
This gives a discretization of the second order differential equation:
\begin{equation}\label{2-CN}
\frac{u_{n+2} - 2 \, u_{n+1} +u_n }{\Delta t^2} + 
\omega^2 \frac{u_{n+2} + 2 \, u_{n+1} +u_n }{4} = 0  \,.
\end{equation}
Then
\begin{align*}
C_{n+1}^2 - C_n^2 & = 
\frac{1}{2} \left(v_{n+1} + v_n\right) \, \left(v_{n+1} - v_n\right) + 
\frac{1}{2} \left(u_{n+1} + u_n\right) \, \left(u_{n+1} - u_n\right) \\
& =
  \frac{\Delta t \, \omega}{4} \,
  \left(v_{n+1} + v_n\right) \,
  \left(u_{n+1} + u_n\right) -
  \frac{\Delta t \, \omega}{4} \,
  \left(u_{n+1} + u_n\right) \,
  \left(v_{n+1} + v_n\right) \equiv 0 \,,
\end{align*}
so $C_n$ is conserved.  Write the system as
\begin{align*}
u_{n+1} + \frac{\Delta t \, \omega}{2} v_{n+1} & =
   u_n  - \frac{\Delta t \, \omega}{2} v_n \\
v_{n+1} - \frac{\Delta t \, \omega}{2} u_{n+1} & =
   v_n  + \frac{\Delta t \, \omega}{2} u_n
\,,
\end{align*}
so that the scheme is implicit, that is it involves the inversion of a
$2\times2$ matrix. The coefficient matrix is always invertible, so there
is no restriction on the size of $\Delta t$, that is, the scheme is
unconditionally stable.

\subsection{The Conservation Law First}

Following the discussion in \cite{WanBihloNave2015} it is easy to show that
the only reasonable discretization that conserves the simple conservation
law \eqref{Simple Conserved} is equivalent to the Crank-Nicholson
discretization. First compute using \eqref{Simple Conserved} that
\[
C_{n+1}^2 - C_n^2 = 
\left(u_{n+1}-u_n\right)
\left(u_n+u_{n+1}\right) +
\left( v_{n+\frac{3}{2}} - v_{n-\frac{1}{2}} \right)
\frac{v_{n-\frac{1}{2}}+2 v_{n+\frac{1}{2}}+v_{n+\frac{3}{2}}}{4} \,.
\]
Choosing
\[
\frac{u_{n+1}-u_n}{\Delta t} =
- \omega \frac{v_{n-\frac{1}{2}}+2 v_{n+\frac{1}{2}}+v_{n+\frac{3}{2}}}{4}
\]
and
\[
\frac{v_{n+\frac{3}{2}} - v_{n-\frac{1}{2}}}{2 \, \Delta t} =
\omega \frac{u_n+u_{n+1}}{2}
\]
will make the $C_n$ constant.
If $\alpha =  \Delta t \, \omega / 2 $ then these equations can be written
\begin{align*}
u_{n+1} + \frac{\alpha}{2} v_{n+3/2} & =
	u_n - \alpha v_{n+1/2} -\frac{\alpha}{2} v_{n-1/2} \,, \\
- 2 \, \alpha \, u_{n+1}  + v_{n+3/2} & = 2 \, \alpha u_{n} +v_{n-1/2} \,.
\end{align*}
So the difference equations are implicit.

It is easy to check that $u_n$ satisfies the second order difference
equation \eqref{2-CN}. Unfortunately, this discretization produces the
same $u_n$ values as the Crank-Nicholson scheme but with a greater
computational cost.  Setting
\[
v_n = \frac{v_{n+1/2}+v_{n-1/2}}{2} \,,
\]
converts this scheme along with it's conserved quantity to the
Crank-Nicholson scheme along with it's conserved quantity.  

\newpage \clearpage
\setcounter{equation}{0}

\section{Details for Discrete Conserved Quantities \label{Appendix Details}}

\subsection{Scalar Wave}

As before a second order discrete equation and a second order average
will be needed
\begin{align*}
\frac{u^{n+1} - 2\, u^{n} + u^{n-1}}{\dt^2}
& = a^{-1}\,\frac{\DIVGstar v^{n+\half}- \DIVGstar v^{n-\half}}{\dt} \\
& = a^{-1}\, \DIVGstar \frac{v^{n+\half}- v^{n-\half}}{\dt} \\
& = a^{-1}\, \DIVGstar {\bf A} \GRAD u^n \\
\frac{u^{n+1} + 2\, u^{n} + u^{n-1}}{4}
& =   u^n + \frac{u^{n+1} - 2\, u^{n} + u^{n-1}}{4} \\
& =   u^n + \frac{\dt^2}{4}\frac{u^{n+1} - 2\, u^{n} + u^{n-1}}{\dt^2} \\
& =   u^n + \frac{\dt^2}{4} a^{-1}\, \DIVGstar {\bf A} \GRAD u^n 
\\
\end{align*}

To find a conserved quantity let
\begin{align*}
C1_{n+1/2} & = \norm{\frac{u^{n+1} + u^{n}}{2}}_\Nodes^2 \,, \\
C2_{n+1/2} & = \norm{ v^{n+1/2}}_\StarFaces^2 \,, \\
C3_{n+1/2} & = \Delta t^2 \norm{ a^{-1} \, \DIVGstar \, v^{n+1/2} }_\Nodes^2 \,.
\end{align*}
As before compute:
\begin{align*}
C1_{n+1/2}- C1_{n-1/2}
& = \left< \frac{u^{n+1} + 2 \, u^{n}  + u^{n-1}}{4} , u^{n+1} - u^{n-1} \right>_\Nodes \\
& = \left< u^{n} + \frac{\Delta t^2}{4} a^{-1} \, \DIVGstar \, {\bf A} \GRAD u^{n} , u^{n+1} - u^{n-1} \right>_\Nodes \,;  \\
& =
\left< u^{n} , u^{n+1} - u^{n-1} \right>_\Nodes 
+
\frac{\Delta t^2}{4}
\left< a^{-1} \, \DIVGstar \, {\bf A} \GRAD u^{n}
, u^{n+1} - u^{n-1} \right>_\Nodes \,;
\end{align*}
Using the adjoint equation xxx gives
\begin{align*}
C2_{n+1/2}- C1_{n-1/2}
& = \left< v^{n+1/2} + v^{n-1/2} , v^{n+1/2} - v^{n-1/2} \right>_\StarFaces \\
& = \left< v^{n+1/2} + v^{n-1/2} , \Delta t \, {\bf A} \GRAD u^n \right>_\StarFaces \\
& = - \Delta t \, \left< a^{-1} \, \DIVGstar \, v^{n+1/2} + a^{-1} \, \DIVGstar \, v^{n-1/2} , u^n \right>_\Nodes \\
& = -\Delta t \, \left< \frac{u^{n+1}-u^{n-1}}{\Delta t} , u^n \right>_\Nodes \\
& = - \left< u^n , u^{n+1}-u^{n-1} \right>_\Nodes \,;
\end{align*}
Also
\begin{align*}
C3_{n+1/2}- C1_{n-1/2} & = \Delta t^2
\left< a^{-1} \, \DIVGstar \, v^{n+1/2} - a^{-1} \, \DIVGstar \, v^{n-1/2} \,,\, a^{-1} \, \DIVGstar \, v^{n+1/2} + a^{-1} \, \DIVGstar \, v^{n-1/2} \right>_\Nodes \\
& = \Delta t^ 2 \left< a^{-1} \, \DIVGstar \left( v^{n+1/2} - v^{n+1/2} \right) \,,\, \frac{u^{n+1}-u^{n-1}}{\Delta t} \right>_\Nodes \\
& = \Delta t^ 2 \left< - \Delta t \, a^{-1} \, \DIVGstar \, -{\bf A} \GRAD u^{n} \,,\,
	\frac{u^{n+1}-u^{n-1}}{\Delta t} \right>_\Nodes \\
& = \Delta t^ 2 \left< a^{-1} \, \DIVGstar \, {\bf A} \GRAD u^{n} \,,\, u^{n+1}-u^{n-1} \right>_\Nodes \,.
\end{align*}
Consequently $C = C1 + C2 - C3/4$ is a conserved quantity:
\[
C_{n+1/2} = \norm{\frac{u^{n+1} + u^{n}}{2}}^2 + \norm{ v^{n+1/2}}^2 
 - \frac{\Delta t^2}{4} \norm{ a^{-1} \, \DIVGstar \, v^{n+1/2} }^2 \,.
\]
This implies that 
\[
C_{n+1/2} \geq \norm{\frac{u^{n+1} + u^{n}}{2}}^2 +
\left( 1 - \frac{\Delta t^2}{4} \norm{a^{-1} \, \DIVGstar}^2 \right)
\norm{ v^{n+1/2}}^2  \,.
\]
So $C_{n+1/2} \geq 0$ for $\Delta t$ sufficiently small provided
$\norm{a^{-1} \, \DIVGstar}$ is finite.

Next look for an analog $C_{n}$ of the scalar conserved quantity
\begin{align*}
C1_{n} & = \norm{\frac{v^{n+1/2} + v^{n-1/2}}{2}}_\StarFaces^2 \,, \\
C2_{n} & = \norm{u^{n}}_\Nodes^2 \,, \\
C3_{n} & = \Delta t^2 \norm{{\bf A} \GRAD \, u^{n} }_\StarFaces^2 \,. 
\end{align*}
First compute
\begin{align*}
C1_{n+1}- C1_{n}
& = \left<
\frac{v^{n+3/2} + 2 \, v^{n 1/2} + v^{n-1/2}}{4} \,,\,
   v^{n+3/2} -v^{n-1/2} \right>_\StarFaces \\
& =
\left< v^{n+1/2} , v^{n+3/2} -v^{n-1/2} \right>_\StarFaces
+ \frac{\Delta t^2}{4}
\left<
{\bf A} \GRAD \, a^{-1} \, \DIVGstar \, v^{n+1/2} \,,\,
v^{n+3/2} -v^{n-1/2} \right>_\StarFaces
\,.
\end{align*}
Using the adjoint equation xxx gives
\begin{align*}
C2_{n+1} - C2_{n} 
& = \left< u^{n+1}-u^{n} \,,\, u^{n+1}+u^{n} \right>_\Nodes \\
& = \left< \Delta t \,  a^{-1} \, \DIVGstar \, v^{n+1/2} \,,\, u^{n+1}+u^{n} \right>_\Nodes \\
& = \Delta t \left< v^{n+1/2} \,,\, -{\bf A} \GRAD u^{n+1}-{\bf A} \GRAD u^{n} \right> \\
& = \Delta t \left< v^{n+1/2} \,,\,
	-\frac{v^{n+3/2}- v^{n-1/2}}{\Delta t} \right> \\
& = - \left< v^{n+1/2} \,,\,
	v^{n+3/2}- v^{n-1/2} \right>  \,.
\end{align*}
Also
\begin{align*}
C3_{n+1}- C3_{n} & = \Delta t^2 
\left<
   {\bf A} \GRAD u^{n+1}-{\bf A} \GRAD u^n  \,,\,  {\bf A} \GRAD u^{n+1}+{\bf A} \GRAD u^n \right>_\StarFaces \\
& = \Delta t^2 \left<  {\bf A} \GRAD u^{n+1} - {\bf A} \GRAD u^n  \,,\, \frac{v^{n+3/2}-v^{n-1/2}}{\Delta t}
\right>_\StarFaces \\
& = \Delta t^2 \left<
   \Delta t \, {\bf A} \GRAD a^{-1} \, \DIVGstar \, v^{n+1/2} \,,\, \frac{v^{n+3/2}-v^{n-1/2}}{\Delta t}
\right> \\
& = \Delta t^2 \left<
   {\bf A} \GRAD a^{-1} \, \DIVGstar \, v^{n+1/2} \,,\, v^{n+3/2}-v^{n-1/2} \right> \,.
\end{align*}
Consequently $C_n = C1_n + C2_n - C3_n/4$ is a conserved quantity: 
\[
C_n =  \norm{u^{n}}^2 
 -\frac{\Delta t^2}{4} \norm{ {\bf A} \GRAD \, u^{n} }^2 
 +\norm{\frac{v^{n+1/2} + v^{n-1/2}}{2}}^2 \,.
\]
This implies that
\[
\norm{C_n} \ge	 \left(1 - \frac{\Delta t^2}{4} \norm{{\bf A} \GRAD}^2\right) \norm{u^n}^2
	+ \norm{\frac{v^{n+1/2} + v^{n-1/2}}{2}}^2 \,,
\]
so $\norm{C_n}$ is positive for sufficiently small $\Delta t$ if
$\norm{ {\bf A} \GRAD \, u^{n} }$ is finite.

\subsection{Maxwell}

To study conserved quantities for Maxwell's equations the second order
discrete difference and average will be needed:
\begin{align*}
\frac{\vec{E}^{n+1} - 2\, \vec{E}^{n} + \vec{E}^{n-1}}{\dt^2}
& = - \epsilon^{-1}\, \CURLstar \mu^{-1} \CURL \vec{E}^n \\
\frac{\vec{E}^{n+1} + 2\, \vec{E}^{n} + \vec{E}^{n-1}}{4}
& =   \vec{E}^n -
\frac{\dt^2}{4} \epsilon^{-1}\, \CURLstar \mu^{-1} \CURL \vec{E}^n 
\end{align*}

To find a conserved quantity $C_{n+1/2}$ let
\begin{align*}
C1_{n+1/2} & = \norm{\frac{\vec{E}^{n+1} + \vec{E}^{n}}{2}}_\Edges^2 \,, \\
C2_{n+1/2} & = \norm{ \vec{H}^{n+1/2}}_\StarEdges^2 \,, \\
C3_{n+1/2} & = \Delta t^2 \norm{\epsilon^{-1} \, \CURLstar \, \vec{H}^{n+1/2} }_\Edges^2  \,.
\end{align*}

As before compute:
\begin{align*}
C1_{n+1/2}- C1_{n-1/2}
& = \left< \frac{\vec{E}^{n+1} + 2 \, \vec{E}^{n}  + \vec{E}^{n-1}}{4} , \vec{E}^{n+1} - \vec{E}^{n-1} \right>_\Edges \\
& = \left< \vec{E}^{n} - \frac{\Delta t^2}{4} \epsilon^{-1} \, \CURLstar \, \mu^{-1}
\CURL \vec{E}^{n} , \vec{E}^{n+1} - \vec{E}^{n-1} \right>_\Edges \,;  \\
& =
\left< \vec{E}^{n} , \vec{E}^{n+1} - \vec{E}^{n-1} \right>_\Edges 
- \frac{\Delta t^2}{4}
\left< \epsilon^{-1} \, \CURLstar \, \mu^{-1} \CURL \vec{E}^{n}
, \vec{E}^{n+1} - \vec{E}^{n-1} \right>_\Edges \,;
\end{align*}

Using the adjoint equation xxx gives
\begin{align*}
C2_{n+1/2}- C1_{n-1/2}
& = \left< \vec{H}^{n+1/2} + \vec{H}^{n-1/2} \,, \vec{H}^{n+1/2} - \vec{H}^{n-1/2} \right>_\StarFaces \\
& = \left< \vec{H}^{n+1/2} + \vec{H}^{n-1/2} \,, - \Delta t \, \mu^{-1} \CURL \vec{E}^n \right>_\StarFaces \\
& = - \Delta t \, \left< \epsilon^{-1} \, \CURLstar \, \vec{H}^{n+1/2} + \epsilon^{-1} \, \CURLstar \, \vec{H}^{n-1/2} , \vec{E}^n \right>_\Edges \\
& = -\Delta t \, \left< \frac{\vec{E}^{n+1}-\vec{E}^{n-1}}{\Delta t} , \vec{E}^n \right>_\Edges \\
& = - \left< \vec{E}^n , \vec{E}^{n+1}-\vec{E}^{n-1} \right>_\Edges \,;
\end{align*}
Also
\begin{align*}
C3_{n+1/2}- C1_{n-1/2} & = \Delta t^2
\left< \epsilon^{-1} \, \CURLstar \, \vec{H}^{n+1/2} - \epsilon^{-1} \, \CURLstar \, \vec{H}^{n-1/2} \,,\, \epsilon^{-1} \, \CURLstar \, \vec{H}^{n+1/2} + \epsilon^{-1} \, \CURLstar \, \vec{H}^{n-1/2} \right>_\Edges \\
& = \Delta t^ 2 \left< \epsilon^{-1} \, \CURLstar \left( \vec{H}^{n+1/2} - \vec{H}^{n+1/2} \right) \,,\, \frac{\vec{E}^{n+1}-\vec{E}^{n-1}}{\Delta t} \right>_\Edges \\
& =\Delta t^ 2 \left< - \Delta t \, \epsilon^{-1} \, \CURLstar \, \mu^{-1} \CURL \vec{E}^{n} \,,\,
	\frac{\vec{E}^{n+1}-\vec{E}^{n-1}}{\Delta t} \right>_\Edges \\
& = - \Delta t^ 2 \left< \epsilon^{-1} \, \CURLstar \, \mu^{-1} \CURL \vec{E}^{n} \,,\, \vec{E}^{n+1}-\vec{E}^{n-1} \right>_\Edges \,.
\end{align*}
Consequently $C = C1 + C2 - C3/4$ is a conserved quantity:
\[
C_{n+1/2} =
\norm{\frac{\vec{E}^{n+1} + \vec{E}^{n}}{2}}_\Edges^2
+ \norm{ \vec{H}^{n+1/2}}_\StarEdges^2 
- \frac{\Delta t^2}{4}
\norm{ \epsilon^{-1} \, \CURLstar \, \vec{H}^{n+1/2} }_\Edges^2 \,.
\]
This implies that 
\[
C_{n+1/2} \geq
\norm{\frac{\vec{E}^{n+1} + \vec{E}^{n}}{2}}_\Edges^2 +
\left( 1 - \frac{\Delta t^2}{4} \norm{\epsilon^{-1} \, \CURLstar}^2 \right)
\norm{ \vec{H}^{n+1/2}}_\StarEdges^2  \,.
\]
So $C_{n+1/2} \geq 0$ for $\Delta t$ sufficiently small provided
$\norm{\epsilon^{-1} \, \CURLstar}$ is finite.

Next look for a conserved quantity $C_{n}$:
\begin{align*}
C1_{n} & = \norm{\frac{\vec{H}^{n+1/2} + \vec{H}^{n-1/2}}{2}}_\StarFaces^2 \,, \\
C2_{n} & = \norm{\vec{E}^{n}}_\Edges^2 \,, \\
C3_{n} & = \Delta t^2 \norm{\mu^{-1} \CURL \, \vec{E}^{n} }_\StarFaces^2 \,. 
\end{align*}
First compute
\begin{align*}
C1_{n+1}- C1_{n}
& = \left<
\frac{\vec{H}^{n+3/2} + 2 \, \vec{H}^{n 1/2} + \vec{H}^{n-1/2}}{4} \,,\,
   \vec{H}^{n+3/2} -\vec{H}^{n-1/2} \right>_\StarFaces \\
& =
\left< \vec{H}^{n+1/2} , \vec{H}^{n+3/2} -\vec{H}^{n-1/2} \right>_\StarFaces
+ \frac{\Delta t^2}{4}
\left<
\mu^{-1} \CURL \, \epsilon^{-1} \, \CURLstar \, \vec{H}^{n+1/2} \,,\,
\vec{H}^{n+3/2} -\vec{H}^{n-1/2} \right>_\StarFaces
\,.
\end{align*}
Using the adjoint equation xxx gives
\begin{align*}
C2_{n+1} - C2_{n} 
& = \left< \vec{E}^{n+1}-\vec{E}^{n} \,,\, \vec{E}^{n+1}+\vec{E}^{n} \right>_\Edges \\
& = \left< \Delta t \,  \epsilon^{-1} \, \CURLstar \, \vec{H}^{n+1/2} \,,\, \vec{E}^{n+1}+\vec{E}^{n} \right>_\Edges \\
& = \Delta t \left< \vec{H}^{n+1/2} \,,\, -\mu^{-1} \CURL \vec{E}^{n+1}-\mu^{-1} \CURL \vec{E}^{n} \right>_\StarFaces \quad\text{(adjoint)} \\
& = \Delta t \left< \vec{H}^{n+1/2} \,,\,
	-\frac{\vec{H}^{n+3/2}- \vec{H}^{n-1/2}}{\Delta t} \right>_\StarFaces \\
& = - \left< \vec{H}^{n+1/2} \,,\,
	\vec{H}^{n+3/2}- \vec{H}^{n-1/2} \right>_\StarFaces  \,.
\end{align*}
Also
\begin{align*}
C3_{n+1}- C3_{n} & = \Delta t^2 
\left<
   \mu^{-1} \CURL \vec{E}^{n+1}-\mu^{-1} \CURL \vec{E}^n  \,,\,  \mu^{-1} \CURL \vec{E}^{n+1}+\mu^{-1} \CURL \vec{E}^n \right>_\StarFaces \\
& = \Delta t^2 \left<  \mu^{-1} \CURL \vec{E}^{n+1} - \mu^{-1} \CURL \vec{E}^n  \,,\, \frac{\vec{H}^{n+3/2}-\vec{H}^{n-1/2}}{\Delta t}
\right>_\StarFaces \\
& = \Delta t^2 \left<
   \Delta t \, \mu^{-1} \CURL \epsilon^{-1} \, \CURLstar \, \vec{H}^{n+1/2} \,,\, \frac{\vec{H}^{n+3/2}-\vec{H}^{n-1/2}}{\Delta t}
\right>_\StarFaces \\
& = \Delta t^2 \left<
   \mu^{-1} \CURL \epsilon^{-1} \, \CURLstar \, \vec{H}^{n+1/2} \,,\, \vec{H}^{n+3/2}-\vec{H}^{n-1/2} \right>_\StarFaces \,.
\end{align*}
Consequently $C_n = C1_n + C2_n - C3_n/4$ is a conserved quantity: 
\[
C_n =  \norm{\vec{E}^{n}}_\Edges^2 
 -\frac{\Delta t^2}{4} \norm{ \mu^{-1} \CURL \, \vec{E}^{n} }_\StarFaces^2 
 +\norm{\frac{\vec{H}^{n+1/2} + \vec{H}^{n-1/2}}{2}}_\StarFaces^2 \,.
\]
This implies that
\[
\norm{C_n} \ge	
\left(1 - \frac{\Delta t^2}{4}
\norm{\mu^{-1} \CURL}^2\right)
\norm{\vec{E}^n}_\Edges^2
+ \norm{\frac{\vec{H}^{n+1/2} + \vec{H}^{n-1/2}}{2}}_\StarFaces^2 \,,
\]
so $\norm{C_n}$ is positive for sufficiently small $\Delta t$ if
$\norm{ \mu^{-1} \CURL \, \vec{E}^{n} }$ is finite.

The codes {\tt Maxwell.m} and {\tt MaxwellStar.m} confirm that our algorithms
conserve $C_{n+1/2}$ and $C_n$ to two parts in $10^{16}$. Additionally,
the divergence of the curl of the electric and magnetic fields are constant
to one part in $10^{14}$ when there are no sources.

\newpage \clearpage
\setcounter{equation}{0}
\section{Conservation Laws and Positive Solutions}

Do a search on "Positivity-Preserving" for diffusion equations. xxx

Conservation laws that say the total amount of some positive substance
is conserved play an important role in modeling using partial differential
equations, for example the Navier-Stokes equations
\cite{PeyretTaylor83}(equations 1.5, 1.6 and 1.7) can be put into this form.
To provided some insight into discretizing such conservation laws,
two important but simple cases will be considered.  For a similar discussion
see Chapter 11 in \cite{Oliver16}. Also \cite{GerritsmaPJZ18} develop
mimetic methods for anisotropic 2D diffusion. For positivity preserving
for the 3D heat equation see \cite{GaoWU15}

\subsection{Transport}

The transport equation in one dimension is given by
\[
\frac{\partial \rho}{\partial t} +
\frac{\partial \, v \, \rho}{\partial x} = 0 \,,
\]
where $\rho = \rho(x,t)$ is a density and $v = v(x)$ is the velocity of
transport. An important assumption is that $\rho \geq 0$ as it typically
represents the density of some substance.  The general solution of this
equation is
\[
\rho (x,t) = w(x-vt) \,,
\]
where $w(x) = \rho(x,0)$ is the initial data. This solution is
a right translation of $w(x)$. This equation also has an important
conservation law:
\[
\int_{-\infty}^{\infty} \rho (x,t) \, d x =
\int_{-\infty}^{\infty} w(x) \, d x  \,.
\]
The conserved quantity is the total amount of material being transported.
Also note that if $w(x) \geq 0$ then $\rho(x,t) \geq 0$ for all $t$. These
two properties are central to this discussion.
Our interest is in finite difference discretizations of equations that
have a similar conservation law and maintain the positivity of the solution.

We assume that $\Delta x > 0 $ and use two grids: a primal grid
$x_i = i \, \Delta x$ that has cells $[x_i,x_{i+1}]$ and a grid of
cell centers $x_{i+\half} = (i+\half) \, \Delta x$ where $-\infty < i < \infty$.
Note that if $\rho$ is a density then it has spatial dimension $1/d^k$ in
a space of dimension $k$ suggesting that $\rho$ should be in a cells. If
a primal grid is chosen then the discretization of $\rho$ is
\[
\rho^{n+\half}_{i+\half} \,.
\]
We will use the conservation of material
\[
\Delta x \, \rho^{n+\half}_{i+\half} 
\]
in a cell to discretize this equation as
\[
\Delta x \, \rho^{n+3/2}_{i+\half} = 
	\Delta x \, \rho^{n+\half}_{i+\half} 
	+ \Delta t \, v_i \, \rho^{n+\half}_{i-\half}
	- \Delta t \, v_{i+1} \, \rho^{n+\half}_{i+\half} \,.
\]
Rewrite this as
\[
\frac{\rho^{n+3/2}_{i+\half} - \rho^{n+\half}_{i+\half}}{\Delta t}
+ \frac{v_{i+1} \, \rho^{n+\half}_{i+\half} -v_i \, \rho^{n+\half}_{i-\half}}
{\Delta x} =  0 \,,
\]
to see that the discretization is a first order approximation of the
differential equation. As an update of the density the equation becomes
\[
\rho^{n+3/2}_{i+\half} =
\rho^{n+\half}_{i+\half} 
	+ \frac{\Delta t}{\Delta x} \, v_i \, \rho^{n+\half}_{i-\half}
	- \frac{\Delta t}{\Delta x} \, v_{i+1} \, \rho^{n+\half}_{i+\half} \,.
\]
Now if 
\[
\frac{\Delta t}{\Delta x} \, v_i \geq 0  \,,\quad
1 - \frac{\Delta t}{\Delta x} \, v_{i+1} \geq 0 \,,
\]
that is if 
\[
v_i \geq 0 \,,\quad
\frac{\Delta t}{\Delta x} \, v_{i+1} \leq 1 \,,
\]
then the discretization preserves the positivity of the discrete solution
and is the well known upwind scheme. This scheme is not useful if 
the velocity $v = v(x)$ has both negative and positive values.
To fix this consider $v$ rather than $\rho$.

\begin{figure}
\begin{center}
\begin{tabular}{cc}
   \includegraphics[width=2.75in,trim = 0 0 0 350]{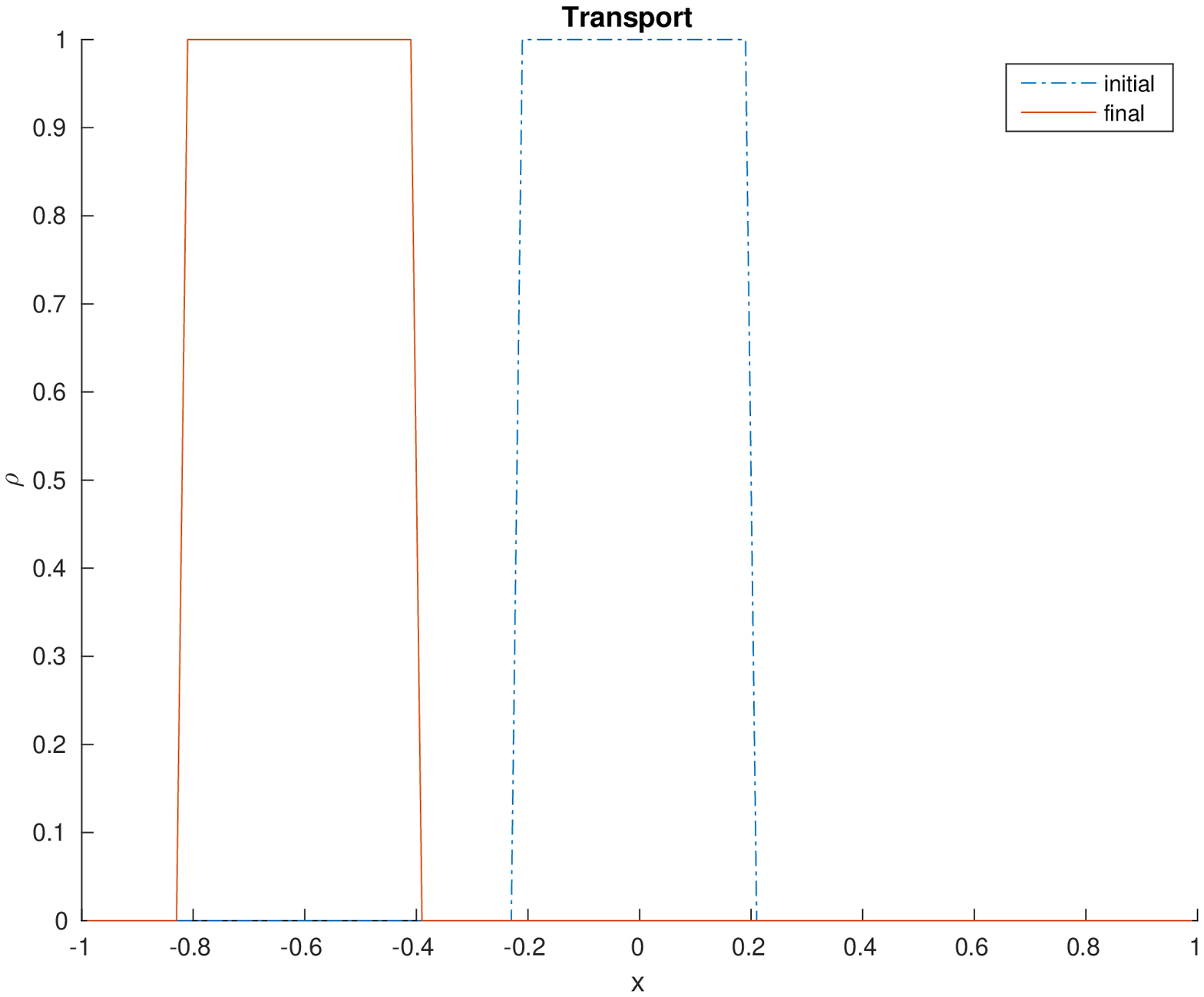} &
   \includegraphics[width=2.75in,trim = 0 0 0 350]{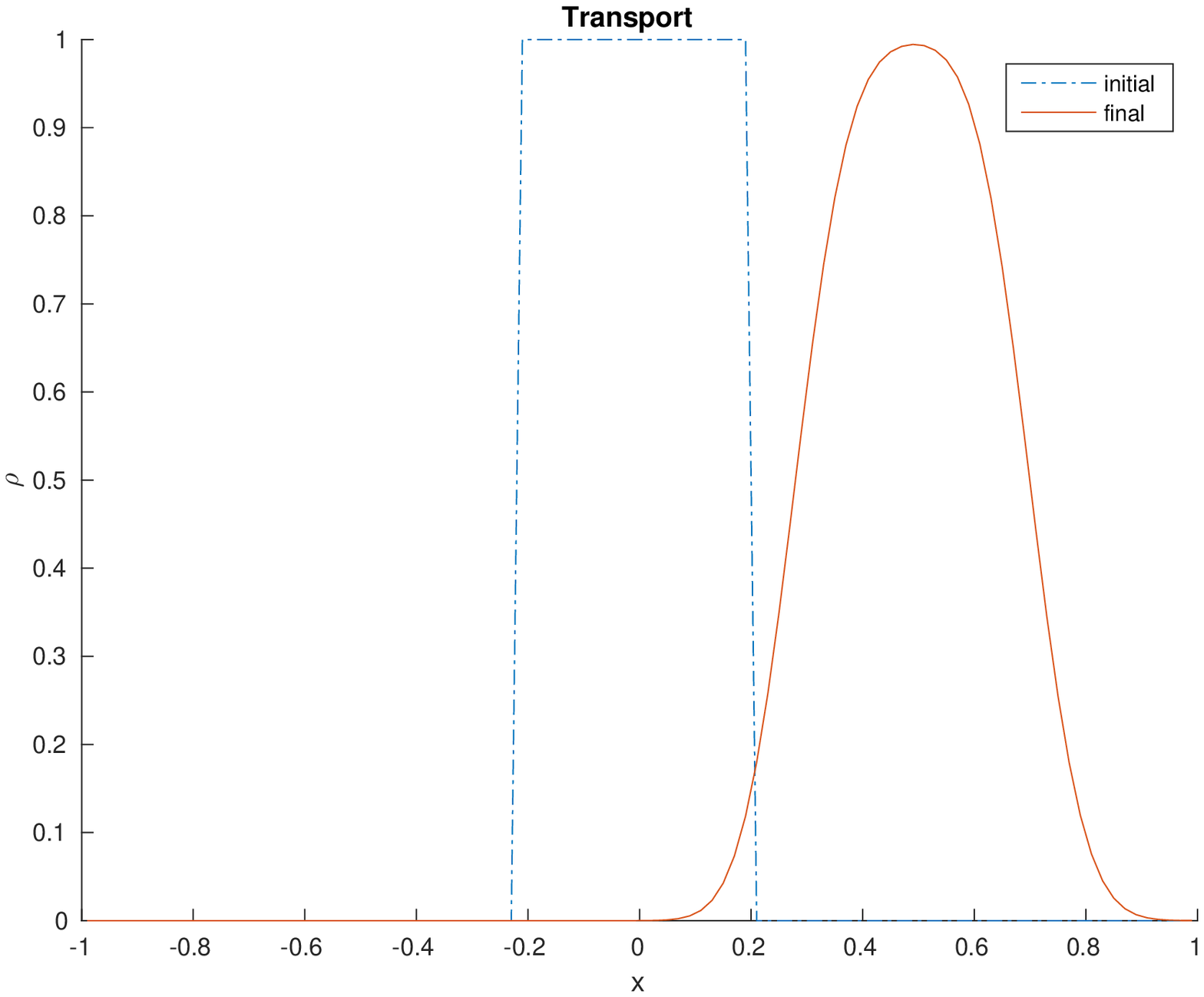} \\
A & B \\
\end{tabular}
\caption{A: Left transport of a square wave $v \, \Delta t / \Delta x = -1$.
B: Right transport of a square wave with $v = 0.4167$. (See {\tt Transport.m})}
\label{Transport Figure}
\end{center}
\end{figure}

So consider the edges of the cells and compute the amount of material
being transferred between the neighboring cells, that is for each time
step $n$, for all $i$ compute the discrete solution as follows:
\begin{align*}
\text{if } v_i \geq 0 \text{ then} \quad & 
\rho^{n+3/2}_{i-\half} = \rho^{n+3/2}_{i-\half}
- v_i \frac{\Delta t}{\Delta x} \, \rho^{n+\half}_{i-\half} \,;\\
& \rho^{n+3/2}_{i+\half} = \rho^{n+3/2}_{i+\half}
+ v_i \frac{\Delta t}{\Delta x} \, \rho^{n+\half}_{i-\half} \,; \\
\text{if } v_i \leq 0 \text{ then}\quad &
\rho^{n+3/2}_{i-\half} = \rho^{n+3/2}_{i-\half}
-  v_i \frac{\Delta t}{\Delta x} \, \rho^{n+\half}_{i+\half} \,;\\
& \rho^{n+3/2}_{i+\half} = \rho^{n+3/2}_{i+\half}
+ v_i \frac{\Delta t}{\Delta x} \, \rho^{n+\half}_{i+\half} \,.
\end{align*}
If $v_i$ is positive then this removes some material from cell
$i-\half$ and put it into cell $i+\half$ and conversely if $v_i$
is negative. If $V = \max(|v_i|)$ then the most material that can
be removed from cell $i-\half$ is 
\[
V \,  \frac{\Delta t}{\Delta x} \rho^{n+\half}_{i-\half} \,,
\]
so to keep $\rho \geq 0$ it must be that
\[
V \,  \frac{\Delta t}{\Delta x}  \leq 1 \,.
\]
An interesting feature of this algorithm is that for 
$v_i \,  \Delta t / \Delta x = \pm 1$ it gives an exact solution
of solution as shown in Figure \ref{Transport Figure}.
This is an upwind scheme for velocities that change direction that
keeps that preserves $\rho \geq 0$ and conserves the amount material
being transported. As done in {\tt Transport.m} this scheme can be
implemented with out the conditional in the update loop.

\begin{figure}
\begin{center}
\begin{tabular}{cc}
   \includegraphics[width=2.75in,trim = 0 0 0 350]{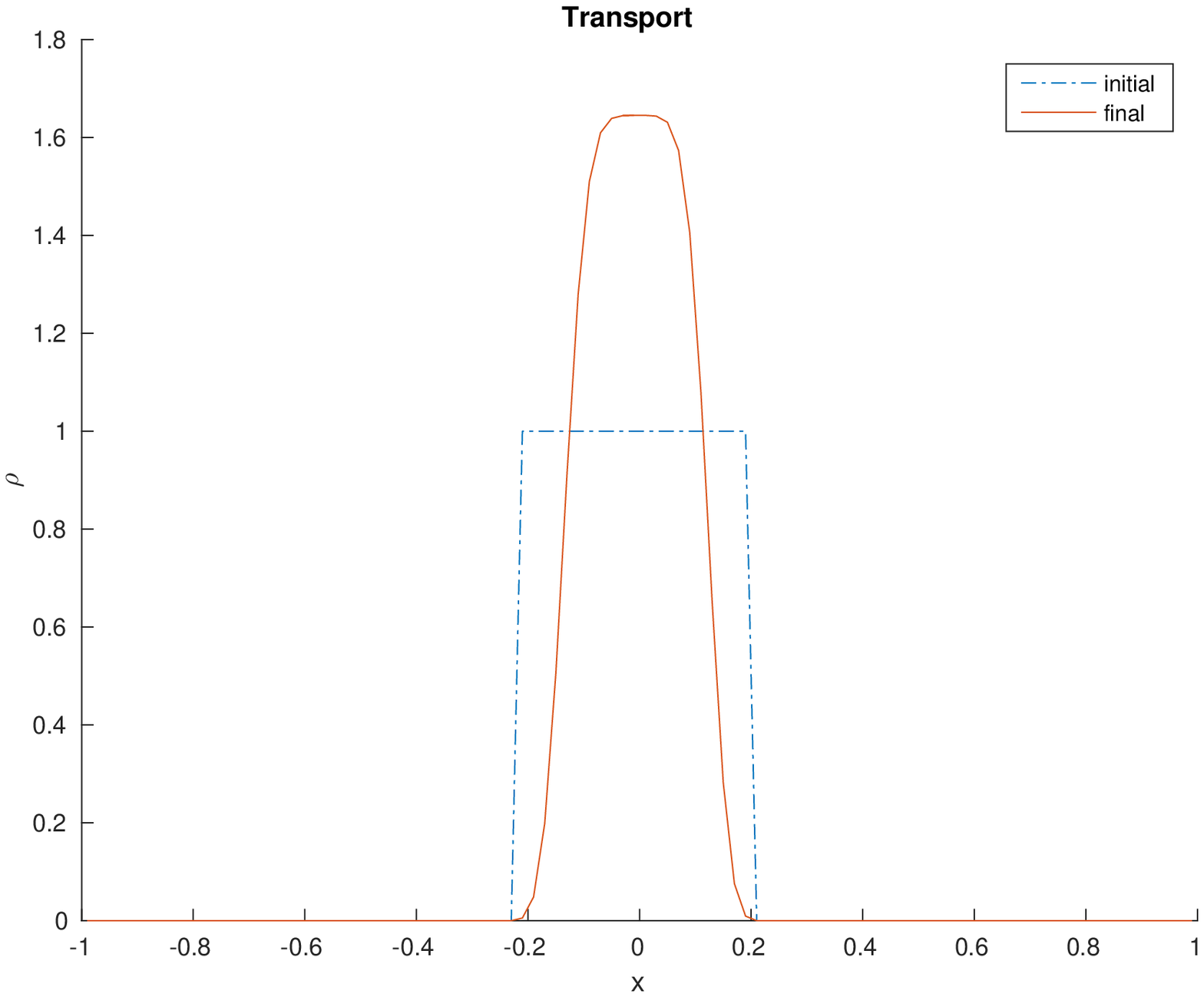} &
   \includegraphics[width=2.75in,trim = 0 0 0 350]{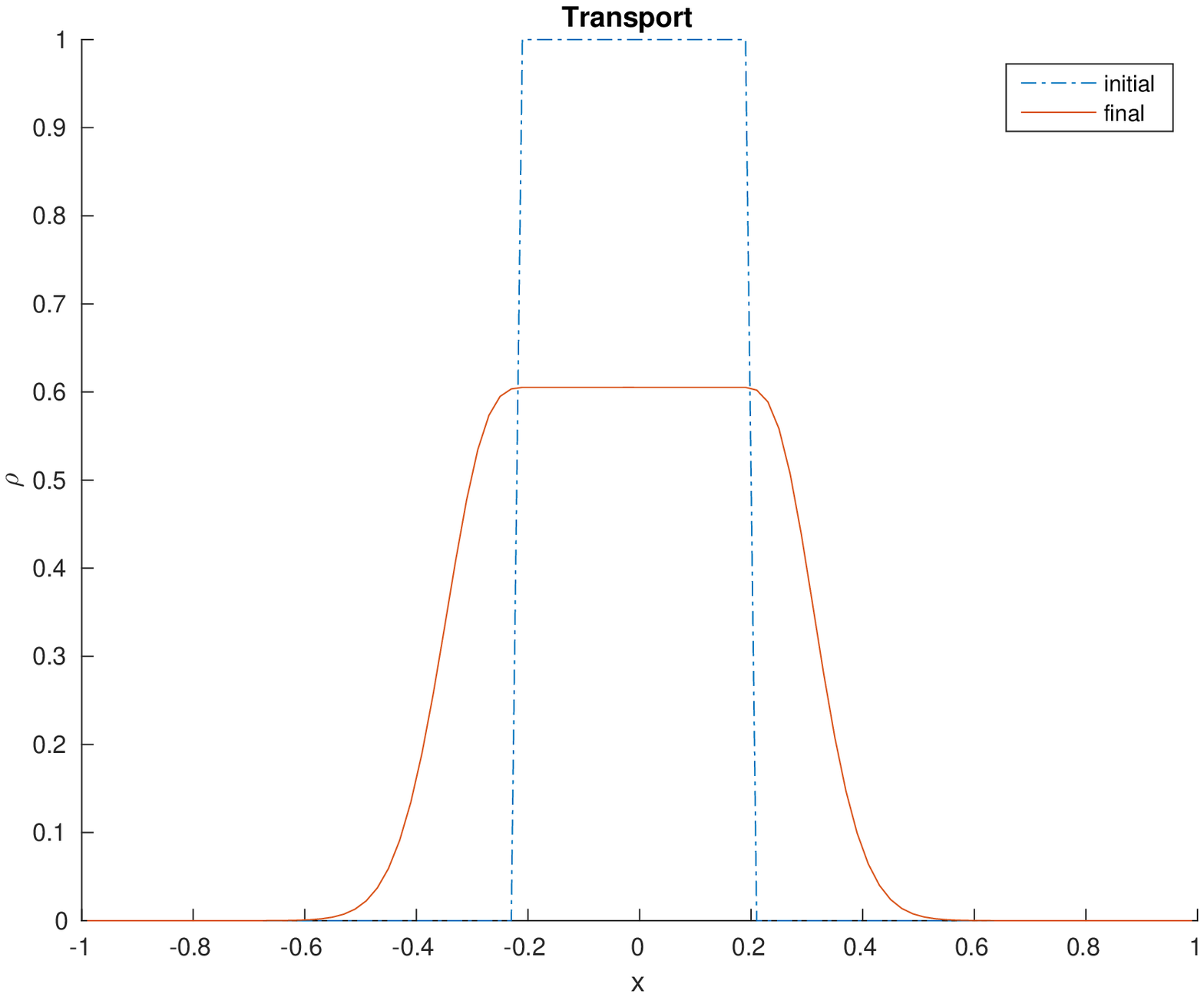} \\
A & B \\
\end{tabular}
\caption{A: Collapse with $v = -x$. B: Expand with $v = x$.
(See {\tt Transport.m})}
\label{Colapse Expand Figure}
\end{center}
\end{figure}

Not all discretizations preserve positive solutions, for example
the Lax-Wendroff, Richtmyer, and MacCormac schemes do not for linear
equations (see {\tt Lax-Wendroff-Positive.nb}).
This can also be seen by by choosing initial data $f_i$ that are
all zero except for one $i$ where $f_i = 1$. For linear equations
the Richtmyer and MacCormac schemes produce the same solution as
the Lax-Wendroff scheme.

\subsection{Diffusion}

The diffusion equation in one dimension is given by
\[
\frac{\partial \rho}{\partial t} =
\frac{\partial}{\partial x}  D \, \frac{\partial \rho}{\partial x} \,,
\]
where $\rho = \rho(x,t)$ is the heat density $D = D(x) \geq 0$ is the
diffusion coefficient. For this discussion $t \geq 0$ and $\rho$ is smooth
and zero for large values of $|x|$. Then integrating the differential
equation gives
\[
\int_{-\infty}^{\infty} \rho(x,t) \, dx = 0 \,.
\]
If $\rho(x,0) \geq 0$ then the solution of the equation is given by
convolution with a Gaussian so then $\rho(x,t)\geq 0 $ for $t \geq 0$.

The standard forward time center space finite difference discretization
of this equation is given by
\[
\frac{\rho_{i+\half}^{n+\half} - \rho_{i+\half}^{n-\half}}{\Delta t} =
\frac{1}{\Delta x} \left(
D_{i+1} \frac{\rho_{i+\thalf}^{n-\half} - \rho_{i+\half}^{n-\half}}{\Delta x}-
D_{i} \frac{\rho_{i+\half}^{n-\half} - \rho_{i-\half}^{n-\half}}{\Delta x}
\right)
\]
or in computational form
\[
\rho_{i+\half}^{n+\half} =
\rho_{i+\half}^{n-\half} +
\frac{\Delta t}{\Delta x^2}
\left(
 D_{i+1} \, \rho_{i+\thalf}^{n-\half}
-\left( D_{i+1} + D_{i} \right)  \, \rho_{i+\half}^{n-\half}
+D_{i}   \, \rho_{i-\half}^{n-\half}
\right).
\]
This algorithm will preserve positive solutions for
\[
\left( D_{i+1} + D_{i} \right) \, \frac{\Delta t}{\Delta x^2} \leq 1 \,,
\]
which is the standard stability constraint for this discretization.

\newpage \clearpage
\setcounter{equation}{0}
\section{Other Stuff}

To be revised.

The paper \cite{VantHofVuik2019} extends mimetic methods in
2D to curvilinear grids for several linear wave equations.

Section \ref{3D Wave Equations} reviews some continuum wave equations in 3
dimensions when the material properties are constant.  The main issue is
understanding the role of the spatial dimension for distance $d$ plays in the
partial differential equations.
The scalar and vector wave equations, the
elastic wave equation \cite{EtgenObrian07}
and Maxwell's equations are introduced and a conserved
quantity is given for each equation. See \cite{Yee1966}
\cite{Collino2006} for a related explicit scheme for Maxwell equation.

In Section \ref{Second Order D0s} continuum second order differential
operators for anisotropic and inhomogeneous materials are introduced.  The
main idea is to use the notion of a double exact sequence and diagram
chasing to define a large class of second order spatial differential operators
that can be used to define wave equations.  This idea is motivated by the
exact sequences used in differential geometry.  A knowledge of differential
forms is not required for understanding this material but can be helpful
\cite{Arapura99}.  A discrete double exact sequence is critical for the
discussion of discretizations using staggered grids.  Importantly, the
discrete double exact sequence cannot be reduced to a single
sequence.  The paper \cite{Palha20141394}, Figure 9, uses a double exact
sequence that is called a De Rham complex.

Additionally, weighted inner products for scalar and vector functions are
introduced and used to define adjoint operators and to show that the second
order operators are either positive or negative.  The main difference between
the discussion here and that in \cite{RobidouxSteinberg2011} is the
introduction of variable material properties.  This discussion depends
heavily on the spatial units of the dependent variables, the differential
operators, and the material properties.

In Section \ref{3D Wave Equations} the second order
differential operators defined in the previous section are combined with a
second time derivative to define several types of wave equations with variable
material properties.  The second order equations are written as a first
order system that has properties similar to the systems studied earlier.
This then gives an automatic definition of a conserved quantity.
At the end of the section Maxwell's equations and the general elastic wave
equations examples are studied.

The material below to be revised soon!
In Section \ref{Mimetic Discretizations} primal and dual staggered grids in
3D are introduced. These grids are the same as those introduce by Yee 
\cite{Yee1966} in 1966 to discretize Maxwell's equations. Consequently there
are two types of discrete scalar fields and two types of discrete vector fields.
The differential operators divergence, gradient and curl are discretized as
in \cite{RobidouxSteinberg2011}.  Because two grids are used there are two
discrete version of each of the first order discrete operators divergence
$\divg$, curl $\curl$ and gradient $\grad$.
Additionally it is shown how to discretize the material properties.  This
section continues by defining discrete inner products and adjoint operators
critical for understanding important properties of the discrete operators.
Note that the paper \cite{MohamedHiraniSataney16} also used a dual grid
differential form method to discretize the Navier-Stokes Equations.
For an introduction to the relationship of vector calculus to
differential forms see the notes \cite{Arapura99}.

In Section \ref{E and M section} we show how to create a conserved quantity
for Maxwell's equations that is constant to one in $10^{15}$, and
also show that the divergence of the electric and magnetic fields
are constant to one part in $10^{13}$, see {\tt Maxwell.m}.

For isotropic and homogeneous materials, simulations show that the three
dimensional scalar wave equation and Maxwell's equations without sources
the approximate energy is constant to less than one part in $10^{15}$.
Additionally, for the scalar wave equation the curl of the velocity is
constant to less than one part in $10^{13}$ and the divergence of the
electric and magnetic fields are constant to less than one part in
$10^{13}$, see {\tt ScalarWave.m} and {\tt Maxwell.m}.

Should section 10 be an appendix?

\subsection{Notes}

This will be revised.

For the latest, see the minisymposium at a recent SIAM meeting \cite{SIAM}.
For more information on steady state problems see \cite{LipnakovBVM14}
For an idea of the difficulties encountered in discretizing Maxwell's
equations see \cite{ChristliebTang14}.
Others have used approximate quantities to study time discretizations 
\cite{HairerLW05,NGI,EngleSD05,GansS00,Quispel08}.
It appears that most of the energy preserving methods are implicit, but by
introducing additional variables, explicit methods that conserve a modified
energy are discussed in \cite{Tao2016}.

One complexity of mimetic spatial discretizations is caused by having
primal and dual grids.  This is leads to there being a primal gradient,
curl, and divergence and dual gradient, curl and divergence. The dual
operators are labeled with a star $\star$. This complexity was already
present in the paper by Yee \cite{Yee1966} which has evolved into the FDTD
discretization method \cite{wikiFDTD}.

There are several minor problems caused by writing wave equations as a
second order differential equations or as a system of two first order
equations.  For example second order equations are not exactly equivalent
to first order system.  Additionally, for the discrete equations there are
problems in converting the initial data for the second order equations to
data for the first order equations and vice versa.  Additionally, because
the equations studied are linear, if  they conserve some quantity, they will
conserve infinitely many quantities and thus there are choices in what
conserved quantity to study. For the first order system there is a natural
{\em primitive} conserved quantity.

This paper was inspired by the papers \cite{WanBihloNave2015} and
\cite{Yee1966}.  We note that in \cite{SchuhmannWeiland2001}
(see equation (45)) the same stability constraint was found as the
one in this paper for conserving the classical energy by modifying
the discretization of Maxwell's equations.  In \cite{GaoZang13} an
implicit (ADI) method is developed that has a modified energy that
is similar to the one used here but the added term is positive while
the added term here is negative.  For a finite element approach that
produce many of the same results that as in this paper see
\cite{TaylorFournier10,BrezziBM14}.

The paper \cite{Sanderse13} gives an overview of energy conserving methods
for Navier-Stokes equations and develops some implicit Runge-Kutta methods
for doing this. The thesis \cite{Capuano15} addresses energy conservation for
turbulent flows.  For a differential forms approach to discretization
see \cite{PerotZusi14,Teixeira13} and additionally for multisymplectic time
integration approach to Maxwell's equations see \cite{StToDeMa2015}.
For two dimensional problems see
\cite{ChenLL08,Perot2011,PalhaGerritsma16,Salmon07,VeigaLV15,MorinishiLCM98,AdjointHymanShashkov97, NaturalHymanShashkov97}.
The papers \cite{Tonti14, WanBihloNave2015} take a novel approach to finding
discrete models.  For a finite-element approach to vector wave equations
see Section 2.3.2 of \cite{arnoldelastic2010}.

For higher order mimetic methods, see \cite{2015SanchezETALAlgorithms}.

Others have used summation by parts to obtain mimetic like discretizations
\cite{WangKreiss2017,FernandezHZ14,NordstromL13}.

It would be very interesting to extend the methods described here to
nonlinear partial differential equations \cite{CaiJY18}.

\newpage \clearpage

\end{appendix}
\newpage \clearpage
\end{document}